\documentclass[final,11pt,times]{elsarticle}
\usepackage[top=0.8in, bottom=0.8in,left=0.5in, right=0.5in]{geometry}
\usepackage{enumitem}
\usepackage{lineno}
\usepackage{graphicx,subfigure,color}
\usepackage{ctex}
\usepackage{float}
\usepackage{amsmath}
\graphicspath{{FIGURE/}}
\usepackage{subfigure}
\usepackage{epstopdf}
\usepackage{algorithm,algorithmic}
\usepackage{amssymb}
\usepackage{setspace}
\usepackage{booktabs,multirow}
\usepackage[english]{babel}
\usepackage{amsthm}
\newtheorem{remark}{Remark}
\newtheorem{theorem}{Theorem}[section]

\newtheorem{assumption}{Assumption}
\newcommand\dlmu[2][4cm]{\hskip1pt\underline{\hbox to #1{\hss#2\hss}}\hskip3pt}
\usepackage{amsmath}
\let\oldequation\equation
\let\oldendequation\endequation
\renewenvironment{equation}{\linenomathNonumbers\oldequation}{\oldendequation\endlinenomath}

\let\oldequation\equation
\let\oldendequation\endequation
\renewenvironment{equation}{\linenomathNonumbers\oldequation}{\oldendequation\endlinenomath}

\let\oldalign\align
\let\oldendalign\endalign
\renewenvironment{align}{\linenomathNonumbers\oldalign}{\oldendalign\endlinenomath}

\let\oldgather\gather
\let\oldendgather\endgather
\renewenvironment{gather}{\linenomathNonumbers\oldgather}{\oldendgather\endlinenomath}

\let\oldgather\gather
\let\oldendgather\endgather
\renewenvironment{gather}{\linenomathNonumbers\oldgather}{\oldendgather\endlinenomath}

\usepackage{pdfpages}
\usepackage{hyperref}
\usepackage{cleveref}
\usepackage[numbers]{natbib}
\biboptions{sort&compress}

\graphicspath{{./figure/}}
\hypersetup{colorlinks = true,
	linkcolor = black,
	urlcolor = blue,
	citecolor = blue}

\makeatletter
\def\@author#1{\g@addto@macro\elsauthors{\normalsize%
		\def\baselinestretch{1}%
		\upshape\authorsep#1\unskip\textsuperscript{%
			\ifx\@fnmark\@empty\else\unskip\sep\@fnmark\let\sep=,\fi
			\ifx\@corref\@empty\else\unskip\sep\@corref\let\sep=,\fi
		}%
		\def\authorsep{\unskip,\space}%
		\global\let\@fnmark\@empty
		\global\let\@corref\@empty  %% Added
		\global\let\sep\@empty}%
	\@eadauthor={#1}
}

\makeatletter
\def\ps@pprintTitle{%
	\let\@oddhead\@empty
	\let\@evenhead\@empty
	\def\@oddfoot{\hfil\thepage\hfil}
	\let\@evenfoot\@oddfoot
}
\makeatother
\begin{document}
%	\begin{linenumbers}
		\begin{frontmatter}
	
	\title{Model and higher-order structure-preserving scheme for the chemotaxis equations modeling the growth of bacterial colonies \\}
	\tnotetext[fund]{This work is supported in part by the Natural Science Foundation of China (Nos. 12431014, 12571417), the 111 Project (No. D23017), and the Postgraduate Scientific Research Innovation Project of Xiangtan University (No. XDCX2025Y192), Project No. 2026CDJZKPT-09 supported by the Fundamental Research Funds for the Central Universities.}
	\author[XTU]{Shuyan Ren}
	\ead{rensymath0820@163.com}
	
	
	\author[CQU]{Kun Wang\corref{cor3}}
	\ead{kunwang@cqu.edu.cn}\cortext[cor3]{Corresponding author.}
	
	\author[HNAM]{Nianyu Yi}
	\ead{yinianyu@xtu.edu.cn }

	\address[XTU]{School of Mathematics and Computational Science, Xiangtan University, Xiangtan, 411105, PR China}
	\address[CQU]{College of Mathematics and Statistics, Key Laboratory of Nonlinear Analysis and Its Applications (Ministry of Education), Chongqing University, Chongqing 401331, PR China}
	\address[HNAM]{National Center for Applied Mathematics in Hunan, School of Mathematics and Computational Science, Xiangtan University, Xiangtan 411105, PR China}
	
	\begin{abstract}
In this paper, we study the model and higher-order structure-preserving numerical scheme for the chemotaxis equations modeling the growth of bacterial colonies. By using a new variational analysis method,  we first rigidly derive the nonnegativity for two solutions of the chemotaxis equations. Then, utilizing the finite element method with  mass-lumping in space and the BDF scheme with  the positivity-preserving extrapolation to deal with the nonlinear terms in  time,   we propose a high-order fully discrete framework for the chemotaxis equations by combining with the Karush-Kuhn-Tucker (KKT) shift and the Lagrange multiplier method, which preserves the mass conservation of the discrete system and the nonnegativity of two approximation solutions. Moreover, apart from confirming the existence and uniqueness of the nonlinear equation generated by the Lagrange multiplier method,  by classifying its structure in detail, we prove that the secant iteration method for this nonlinear equation is locally convergent. Finally, we show some experiments to verify the theoretical prediction and efficiency of the proposed scheme, including consistent simulations of the growth of the colonies under different nutrient concentration and agar hardness.
	\end{abstract}

%\journal{Mathematical Models and Methods in Applied Sciences}	
\end{frontmatter}

\section{Introduction}
\label{1_introduction}
\begin{figure}[htbp]
	\centering
	\includegraphics[width=0.7\textwidth, height=0.5\textwidth]{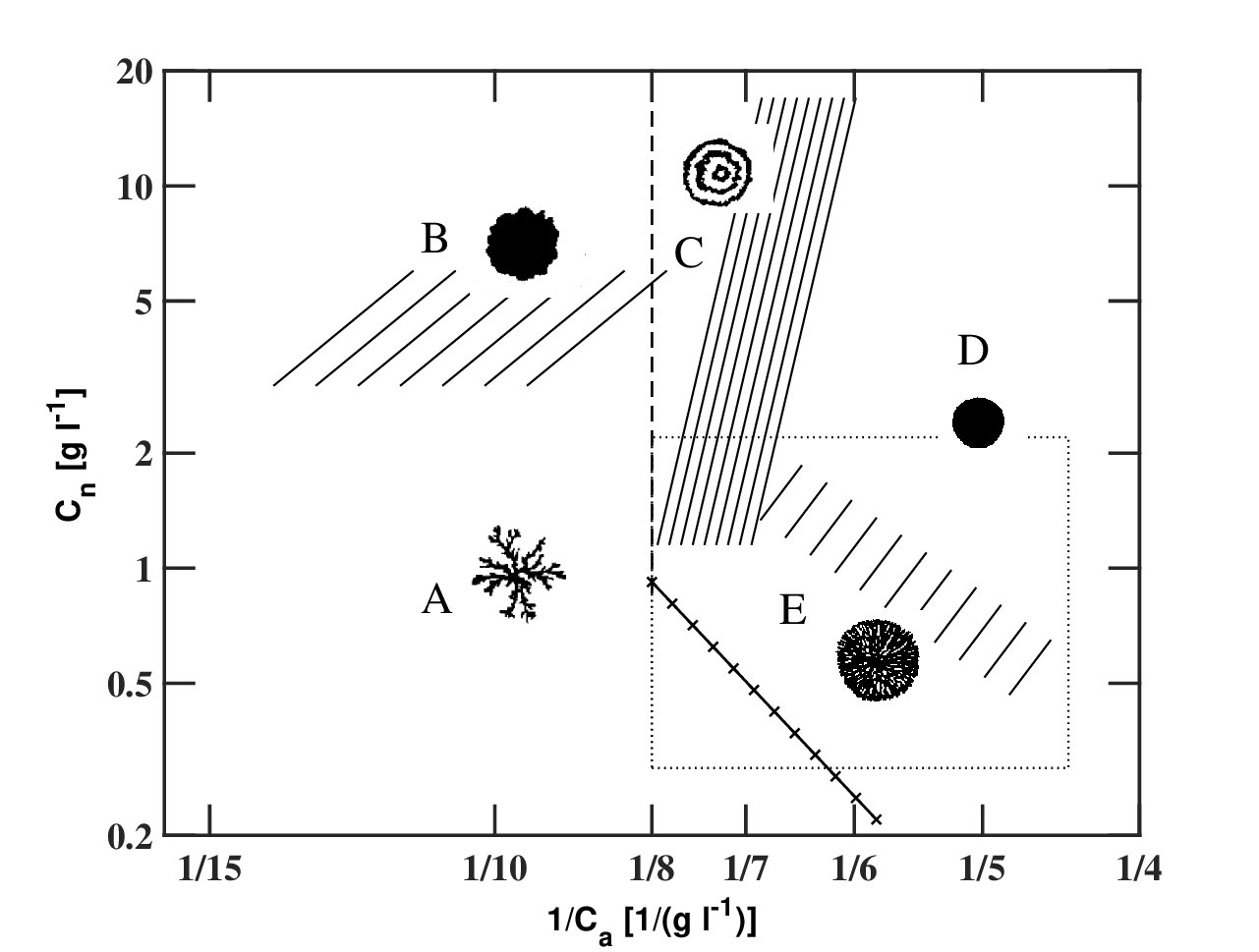}
	\caption{Phase diagram of colony patterns in the $(C_n,1/C_a)$ parameter
		plane, where $C_n$ denotes the nutrient concentration and $C_a$
		the agar concentration, so that $1/C_a$ serves as a measure of
		the softness of the agar medium \cite{KAMT1997}.
	}
	\label{bacterial_behavior}
\end{figure}

Studying bacterial collective behavior, surface motility and colony morphogenesis is always an interesting topic (see, e.g., \cite{Vlamakis2013,Allard2016}),  especially for the cases of spatially heterogeneous environments such as the plant rhizosphere, nutrients and root exudates generate chemical gradients that can be sensed by bacteria, thereby inducing chemotactic migration toward favorable regions and facilitating nutrient acquisition and early root colonization \cite{Allard2016}. Effective migration and surface spreading are important for bacterial populations to occupy the rhizosphere and establish stable communities, while successful root colonization provides an essential basis for the plant growth-promoting and biocontrol functions of Bacillus subtilis. Therefore, bacterial responses to nutrient gradients and population spreading in spatially constrained environments are closely related to plant-microbe interactions, rhizosphere ecology, and agricultural biocontrol. Thin agar plates provide a typical and well-controlled experimental system for investigating the effects of nutrient availability, the physical properties of the medium, and bacterial motility on population spreading. Experimental observations have shown that the morphology of Bacillus subtilis colonies varies significantly with the nutrient and agar concentrations. As illustrated in Figure \ref{bacterial_behavior}, under low nutrient and high agar concentrations, namely, on a nutrient-poor and relatively hard medium, the colony typically develops a diffusion-limited aggregation (DLA)-like pattern characterized by pronounced slender branches (region A) \cite{MSM2000}. When the agar concentration remains high while the nutrient concentration increases, the branches gradually thicken and tend to merge, eventually forming an Eden-like pattern with an approximately circular shape and a rough outer boundary (region B) \cite{MMT1992}. In region C, alternating phases of active bacterial motion and temporary arrest associated with cell division lead, through repeated cycles, to the formation of pronounced concentric ring-like structures \cite{SMW2017}. On a nutrient-rich and relatively soft medium with a low agar concentration, the constraints on bacterial motility are weakened, allowing the colony to spread more uniformly and develop a macroscopically homogeneous disk-like structure (region D) \cite{MXZ2005}. By contrast, when both the nutrient and agar concentrations are relatively low, the colony exhibits a dense-branching morphology (DBM), characterized by rich internal branching structures and a comparatively smooth outward-propagating envelope (region E) \cite{KAMT1997}. These distinct colony morphologies result from the interplay among nutrient diffusion and depletion, bacterial proliferation and motility, chemotaxis, and the constraints imposed by the surrounding medium. In particular, nutrient availability regulates bacterial growth and resource acquisition, whereas the agar concentration affects colony spreading by altering the hardness of the medium and the resulting constraints on bacterial motility. Chemotaxis, on the other hand, drives directed bacterial migration along nutrient gradients. The coupling of these mechanisms collectively determines bacterial aggregation and spatial organization, the propagation speed of the colony front, and the resulting morphological evolution. Therefore, developing mathematical models that capture the coupled effects of nutrient transport, bacterial proliferation and motility, chemotactic responses, and environmental constraints is important for understanding and predicting bacterial migration, colonization, and colony spreading in nutrient-limited and spatially heterogeneous environments, as well as for gaining further insight into their roles in rhizosphere ecology and biological control.

Mesibov et al. experimentally demonstrated that bacterial chemotactic responses exhibit distinct threshold and saturation effects and interpreted these phenomena in terms of attractant binding to chemoreceptors \cite{RI1976}. Subsequently, Lapidus and Schiller incorporated this concentration-dependent receptor response into the Keller-Segel model to describe bacterial chemotaxis over different attractant concentration ranges \cite{RG1973}. To characterize the observation that bacteria remain highly motile near the nutrient-rich advancing front of a colony but become nearly immotile in the nutrient-depleted interior, Kawasaki et al. \cite{KAMT1997} proposed a coupled diffusion-reaction model for bacterial density and nutrient concentration, in which nonlinear degenerate diffusion is used to describe bacterial motility and a nutrient-dependent reaction term accounts for bacterial proliferation. The model is capable of reproducing the morphological transition from branching colonies to homogeneous disk-like colonies and of capturing the dependence of colony propagation speed on nutrient concentration. Golding et al. subsequently compared several reaction-diffusion models for bacterial colonies and emphasized the important roles of nonlinear degenerate diffusion and chemotaxis in branching growth and morphological evolution \cite{EI2000}. Ben-Jacob et al. further provided a comprehensive discussion of cooperative self-organization and the formation of complex spatiotemporal patterns in microorganisms under adverse environmental conditions, with particular emphasis on cell motility, chemotaxis, chemical communication, and environmental responses \cite{IYIE1998}. Building on these studies, Leyva et al. \cite{FG2013} introduced a chemotactic flux driven by the nutrient concentration gradient and compatible with the nonlinear degenerate diffusion, thereby obtaining a coupled reaction-diffusion-chemotaxis system for investigating the effects of nutrient-driven chemotaxis on colony-front propagation, branch formation, and morphological stability.

On the other hand, considerable effort has been devoted to the theoretical analysis of related bacterium-nutrient coupling systems. Winkler \cite{WM2022} analyzed a doubly degenerate nutrient-taxis model in two spatial dimensions, showing global existence of weak solutions with nonnegative bacterial density and strictly positive nutrient concentration. Li and Winkler \cite{LM2022} subsequently considered the corresponding one-dimensional problem and obtained global weak solutions with the same positivity properties, together with boundedness and large-time behavior. Lankeit and Winkler \cite{LJ2023} provided a comprehensive survey of chemotaxis-consumption systems, summarizing existing results on the existence, boundedness, and large-time behavior of solutions, as well as the formation of spatial structures. For a two-dimensional doubly degenerate nutrient-taxis system with a general cross-diffusive flux, Winkler \cite{WM2024} established the existence of global bounded weak solutions under suitable growth conditions on the bacterial response function. Tran and Yang \cite{TY2025} further investigated the three-dimensional doubly degenerate nutrient-taxis system and established the global existence of weak solutions under broad conditions on the bacterial response. More recently, Zhang and Li \cite{ZL2026} analyzed a planar doubly degenerate nutrient-taxis model supplemented with logistic growth and established global weak solvability under suitable regularity assumptions on the initial data. These studies have substantially advanced the theoretical understanding of such models, particularly with respect to global existence, boundedness, weak solvability, and large-time behavior. However, investigations of the system concerning structural properties such as total mass conservation and the nonnegativity of solutions remain relatively limited. These properties are important for understanding the dynamics of the model and for developing reliable numerical methods.

For the numerical method perspective, a number of numerical discretizations and simulations have been developed for bacterium-nutrient coupling models. Early studies mainly focused on direct numerical simulations of colony propagation and morphological evolution. For example, Leyva et al. \cite{FG2013} employed a second-order explicit Runge-Kutta method combined with a finite difference spatial discretization and GPU parallel computing to perform high-resolution simulations of the reaction-diffusion-chemotaxis system, with particular emphasis on the effects of nutrient chemotaxis on colony propagation speed and aggregation patterns. More recently, increasing attention has been paid to structure-preserving discretizations. Barrera-Cardozo and Rueda-G\'{o}mez proposed a fully discrete mass-lumping finite element scheme for a doubly degenerate nutrient-taxis system \cite{BR2026}. By introducing an auxiliary variable and a regularization technique, their scheme preserves discrete total mass conservation, pointwise estimates for the nutrient concentration, approximate nonnegativity of the bacterial density, and unconditional stability with respect to a modified energy. However, their method employs a first-order semi-implicit Euler discretization in time, provides only regularization-dependent approximate nonnegativity for the bacterial density, and still leads to a nonlinear and semi-coupled discrete system. Despite these advances, the construction of high-accuracy structure-preserving schemes still faces several challenges. First, the nonlinear terms are strongly coupled; a fully implicit treatment leads to complicated nonlinear algebraic systems, whereas linearization may destroy the intrinsic structural properties of the continuous system. Second, total mass conservation constitutes a global equality constraint, while nonnegativity represents a local inequality constraint, making it difficult to preserve both properties simultaneously at the discrete level. Finally, mass correction, which is generally used in the references usually introduces a nonlinear algebraic equation, whose solvability and the convergence of the corresponding iterative method also require careful analysis.

Designing efficient numerical methods that inherit these structural properties at the discrete level is crucial for obtaining reliable numerical simulations. The mass-lumping finite element method (FEM) diagonalizes the mass matrix arising from the standard finite element discretization through appropriate numerical quadrature, thereby reducing computational complexity. Owing to its nodal quadrature structure, it is also particularly suitable for imposing pointwise constraints. For higher-order finite elements, a direct application of standard mass lumping may lead to zero or negative quadrature weights. To overcome this difficulty, suitable bubble functions are usually introduced to construct higher-order mass-lumping finite element spaces with positive quadrature weights while retaining the desired high-order approximation accuracy \cite{GPJ2001}. On the other hand, the Karush-Kuhn-Tucker (KKT) conditions and Lagrange multiplier methods provide a natural framework for handling inequality and conservation constraints. In recent years, Cheng and Shen employed this idea to construct structure-preserving numerical schemes, in which the predicted solution is corrected through the KKT conditions and Lagrange multipliers within a predictor-corrector framework, so that structural properties such as positivity, bound preservation, and mass conservation can be maintained \cite{CS2022I, CS2022II}. Based on these observations, we study the model and structure-preserving simulation for the chemotaxis equations modeling the growth of bacterial colonies in this paper, which has the following main contributions:

\begin{itemize}
  \item Besides the mass conservation, through using a new variational  analysis method, we  establish  the nonnegativity for two  solutions of  the chemotaxis system modeling the bacterial colony growth.
  \item We propose a high-order structure-preserving fully discrete framework for the chemotaxis equations, which consists of utilizing the FEM with mass-lumping in space and the BDF scheme in time in the predictor step,  then  recovering the mass conservation and the nonnegativity that may be lost during the linearized prediction by the KKT shift and the Lagrange multiplier method  in the corrector step.
  \item Apart from confirming the existence and uniqueness of the nonlinear equation generated by the Lagrange multiplier method, by classifying its structure in detail, we prove that the secant iteration method for this nonlinear equation is locally convergent.
  \item  Consistent simulations for the growth of the  bacterial colonies  are presented under different nutrient concentration and agar hardness.
\end{itemize}

The rest of the paper is structured as follows. In Section \ref{2_model}, we recall the development of the chemotaxis equations modeling bacterial colonies and derive the model's structures, especially the nonnegativity of two  solutions for  the chemotaxis system. In Section \ref{3_structure_preserving}, using the mass-lumping finite element method and BDF scheme in the predictor step, and the KKT conditions and the Lagrange multiplier technique in the corrector step, we develop the high-order structure-preserving fully discrete scheme, which can achieve second- and third-order in time and space, respectively. Then, some numerical experiments are presented in Section \ref{4_test}, where we examine the temporal and spatial convergence orders, verify their structure-preserving properties, and investigate bacterial pattern formation. Finally, we end with conclusions in Section \ref{5_conclusion}.

\section{Model and its structures}
\label{2_model}		
In this section, we first introduce the model to describe the growth of bacterial colonies. Then, we will deduce the structures of this model, which is very important in developing high-fidelity numerical methods.	
\subsection{Model}
We consider a chemotaxis model for the growth of bacterial colonies on the surface of a thin agar medium. Mesibov et al. identified threshold and saturation effects in bacterial chemotaxis \cite{RG1973}, which Lapidus and Schiller subsequently incorporated into a concentration-dependent extension of the Keller-Segel model \cite{RI1976}. In the original framework of Kawasaki et al. \cite{KAMT1997}, nutrients diffuse through the medium and are consumed by bacteria, while the bacteria proliferate through nutrient uptake and migrate over the agar surface by diffusion. Accordingly, the nutrient concentration and the bacterial density are governed by the following reaction-diffusion system:
\begin{subequations}
	\label{system_1}
	\begin{alignat}{2}
		n_t & = D_n \Delta n - f(n,b),
		&\quad &(\mathbf{x},t)\in\Omega\times(0,T], \label{eq:n} \\
		b_t & = \nabla \cdot (D_b \nabla b) + \tilde{\theta} f(n,b),
		&\quad &(\mathbf{x},t)\in\Omega\times(0,T], \label{eq:b}
	\end{alignat}
\end{subequations}
where $\Omega\subset\mathbb{R}^2$ is a bounded domain and $T>0$ is the final time. Moreover, $n=n(\mathbf{x},t)=n(x,y,t)$ denotes the nutrient concentration, $b=b(\mathbf{x},t)=b(x,y,t)$ denotes the bacterial density, $D_n>0$ is the diffusion coefficient of the nutrient, $D_b$ is the bacterial diffusion coefficient, $\tilde{\theta}>0$ is the conversion rate from nutrient consumption to the bacterial growth, and $f(n,b)$ denotes the nutrient consumption rate by bacteria.

A central modeling ingredient is the choice of the bacterial diffusion coefficient $D_b$. When the nutrient supply is insufficient or the bacterial density is very low, bacterial cells are nearly immotile. In contrast, the bacterial motility is substantially enhanced as both the nutrient concentration and the bacterial density increase. To represent this dependence, Kawasaki et al. \cite{KAMT1997} modeled the bacterial diffusivity by the degenerate nonlinear coefficient
$$
D_b=\sigma nb.
$$
The parameter $\sigma$ is written as $\sigma=\sigma_0(1+\zeta)$. Here, $\sigma_0>0$ denotes the mean bacterial motility coefficient associated with the softness of the agar medium, while $\zeta$ represents a random fluctuation. Specifically, $\zeta$ is assumed to follow a triangular distribution on $[-1,1]$ and is generated by the inverse transform method \cite{LN1986} from pseudorandom variables uniformly distributed on $[0,1]$. Larger values of $\sigma_0$ correspond to lower agar concentrations and hence softer media.

For the $f(n,b)$, we assume that it follows Michaelis-Menten-type kinetics of the form \cite{FG2013}
$$
f(n,b)=\frac{\kappa nb}{1+\gamma n}.
$$
Here, $\kappa>0$ and $\gamma>0$ are prescribed constants. For sufficiently small nutrient concentrations, $\gamma n\ll1$, and hence $1+\gamma n\approx 1$. The reaction term therefore reduces to
$$
f(n,b)\approx \kappa nb.
$$

Accordingly, system \eqref{system_1} can be rewritten as
\begin{subequations}
	\label{system_2}
	\begin{alignat}{2}
		n_t & = D_n \Delta n - \kappa nb,
		&\quad &(\mathbf{x},t)\in\Omega\times(0,T],  \\
		b_t & = \nabla \cdot \bigl(\sigma nb \nabla b\bigl) + \tilde{\theta} \kappa nb,
		&\quad &(\mathbf{x},t)\in\Omega\times(0,T].
	\end{alignat}
\end{subequations}
Although this model already captures colony growth driven by the combined effects of degenerate diffusion and nutrient consumption, Ben-Jacob et al.~\cite{IYIE1998,EI2000} pointed out that the development of bacterial colonies is a self-organized process involving several chemotactic mechanisms. In particular, nutrient-directed chemotaxis plays an important role in accelerating the spreading of the colony front and shaping the branching morphology of the resulting patterns. Accordingly, a nutrient-dependent chemotactic term is added to system \eqref{system_2} to describe this mechanism. We define the chemotactic flux as \cite{RI1976, RG1973}
$$
J_c=\sigma nb^2\frac{\chi_0 K_d}{(K_d+n)^2}\nabla n,
$$
where $\chi_0>0$ is a constant measuring the strength of the chemotactic response, and $K_d>0$ denotes the receptor-ligand dissociation constant. Subtracting the divergence of this flux $J_c$ from the bacterial density equation in \eqref{system_2}, we arrive at the following reaction-diffusion chemotaxis model:
\begin{subequations}
	\label{system_3}
	\begin{alignat}{2}
		n_t &= D_n \Delta n - \kappa nb, &\quad &(\mathbf{x},t)\in\Omega\times(0,T],  \\
		b_t &= \nabla \cdot \bigl(\sigma nb \nabla b\bigr) + \theta \kappa nb
		- \nabla \cdot \left( \sigma nb^2 \frac{\chi_0 K_d}{(K_d+n)^2} \nabla n \right), &\quad &(\mathbf{x},t)\in\Omega\times(0,T].
   \end{alignat}
\end{subequations}

We first introduce the normalized space-time variables
\begin{align}
	\tilde{x} &= \left( \frac{\tilde{\theta} K_d \kappa}{D_n} \right)^{1/2} x, \qquad
	\tilde{y} = \left( \frac{\tilde{\theta} K_d \kappa}{D_n} \right)^{1/2} y, \qquad
	\tilde{t} = \tilde{\theta} K_d \kappa\, t, \nonumber
\end{align}
while $n$, $b$, and $\sigma$ are scaled as
\begin{align}
	\tilde{n} &= \frac{n}{K_d}, \qquad
	\tilde{b} = \frac{b}{\tilde{\theta} K_d}, \qquad
	\tilde{\sigma} = \frac{\tilde{\theta} K_d^2}{D_n}\sigma. \nonumber
\end{align}
With these definitions, system \eqref{system_3} can be written in the following dimensionless form. In what follows, the tildes are omitted for simplicity:
\begin{subequations}
	\label{system_4}
	\begin{alignat}{2}
		n_t &= \Delta n - nb, &\quad &(\mathbf{x},t)\in\Omega\times(0,T], \label{n}\\
		b_t &= \nabla \cdot \bigl(\sigma nb \nabla b\bigl) + nb
		- \chi_0 \nabla \cdot \left( \frac{\sigma nb^2}{(1+n)^2}\nabla n \right), &\quad &(\mathbf{x},t)\in\Omega\times(0,T] \label{b}.
    \end{alignat}
\end{subequations}

To complete \eqref{system_4}, we impose it with the initial conditions
\begin{align}\label{initial_condition}
	n(\mathbf{x},0)=n_0(\mathbf{x}), \qquad
	b(\mathbf{x},0)=b_0(\mathbf{x}), \qquad
	\mathbf{x}\in \Omega,
\end{align}
and the homogeneous Neumann boundary conditions
\begin{align}\label{boundary_condition}
	\frac{\partial n}{\partial \boldsymbol{\nu}}=0, \qquad
	\frac{\partial b}{\partial \boldsymbol{\nu}}=0, \qquad
	(\mathbf{x},t)\in \partial\Omega\times(0,T].
\end{align}
Here, $\partial\Omega$ is the boundary of the domain $\Omega$, and
$\boldsymbol{\nu}$ denotes the unit outward normal vector on $\partial\Omega$.
\subsection{Structures of the continuous model}
System \eqref{system_4}-\eqref{boundary_condition} possesses several fundamental structures, including total mass conservation and nonnegativity of the nutrient concentration and bacterial density. When constructing numerical schemes, it is crucial to preserve these intrinsic features, since they have a direct impact on the reliability of the simulation results. Before proceeding to derive the structure, we first state the following regularity assumption on the exact solution which will be used in the subsequent analysis.
 \begin{assumption}\label{regularity}
 Assume that the initial data $n_0(\mathbf{x})$, $b_0(\mathbf{x})$ are nonnegative, i.e., $n_0(\mathbf{x})\geq 0$ and $b_0(\mathbf{x})\geq 0$.
 \end{assumption}

The chemotaxis equations \eqref{system_4}-\eqref{boundary_condition} have the following structures,
\begin{theorem}[Total mass conservation]
\label{continue_mass}
The total mass of the system \eqref{system_4}-\eqref{boundary_condition} is conserved in time, that is,
	\begin{align}
		\frac{d}{dt} \int_\Omega (n + b) \,d\mathbf{x} = 0.
	\end{align}
\end{theorem}
\begin{proof}
	Integrating the sum of \eqref{n} and \eqref{b} over $\Omega$
	and using the boundary conditions \eqref{boundary_condition}, we obtain
	\begin{align*}
		\frac{d}{dt}\int_\Omega (n+b)\,d\mathbf{x}=0.
	\end{align*}
	This completes the proof.
\end{proof}
 \begin{theorem}[Nonnegativity preservation]  \label{positivity_continous}
	Under Assumption \ref{regularity}, the solution $(n,b)$ of
	\eqref{system_4}-\eqref{boundary_condition} satisfies
	 \begin{align}
	 	n(\mathbf{x},t)\geq 0,\qquad b(\mathbf{x},t)\geq 0,
	 	\qquad (\mathbf{x},t)\in\overline{\Omega}\times[0,T].
	 \end{align}
 \end{theorem}
 \begin{proof}
 	We first prove the nonnegativity of $n(\mathbf{x},t)$. Let
 	$$
 	n_-(\mathbf{x},t):=\max\{-n(\mathbf{x},t),0\}\ge 0
 	$$
 	denote the negative part of $n(\mathbf{x},t)$. Testing \eqref{n} with $n_-(\mathbf{x},t)$ and integrating over $\Omega$ yields
 	$$
 	\int_\Omega n_t n_-\,d\mathbf{x}
 	=
 	\int_\Omega \Delta n\, n_-\,d\mathbf{x}
 	-
 	\int_\Omega nb\,n_-\,d\mathbf{x}.
 	$$
 	For $n(\mathbf{x},t)<0$, the definition of the negative part gives
 	$
 	n(\mathbf{x},t)=-n_-(\mathbf{x},t).
 	$
 	Thus,
 	$$
 	\int_\Omega n_t n_-\,d\mathbf{x}
 	=
 	-\frac12\frac{d}{dt}\|n_-\|_{L^2(\Omega)}^2.
 	$$
 	Moreover, using the homogeneous Neumann boundary condition
 	$\frac{\partial n}{\partial \boldsymbol{\nu}}=0$ in \eqref{boundary_condition}
 	and integration by parts, we have
 	$$
 	\int_\Omega \Delta n\, n_-\,d\mathbf{x}
 	=
 	-\int_\Omega \nabla n\cdot \nabla n_-\,d\mathbf{x}
 	=
 	\|\nabla n_-\|_{L^2(\Omega)}^2.
 	$$
 	Similarly,
 	$$
 	-\int_\Omega nb\,n_-\,d\mathbf{x}
 	=
 	\int_\Omega b\,n_-^2\,d\mathbf{x}.
 	$$
 	Therefore,
 	$$
 	-\frac12\frac{d}{dt}\|n_-\|_{L^2(\Omega)}^2
 	=
 	\|\nabla n_-\|_{L^2(\Omega)}^2
 	+
 	\int_\Omega b\,n_-^2\,d\mathbf{x},
 	$$
 	which implies
 	$$
 	\frac12\frac{d}{dt}\|n_-\|_{L^2(\Omega)}^2
 	\le
 	\|b(\cdot,t)\|_{L^\infty(\Omega)}
 	\|n_-\|_{L^2(\Omega)}^2.
 	$$
 	Since $n_0\ge 0$ according to Assumption \ref{regularity}, we have
 	$n_-(0)=0$. By Gr\"onwall's inequality, we get
 	$$
 	\|n_-\|_{L^2(\Omega)}^2 \le 0.
 	$$
 	Since $\|n_-\|_{L^2(\Omega)}^2$ is nonnegative, we conclude that
 	$$
 	\|n_-\|_{L^2(\Omega)}^2 = 0
 	\qquad \forall t\in[0,T],
 	$$
 	and thus
 	$$
 	n(\mathbf{x},t)\ge 0
 	\quad \text{in } \overline{\Omega}\times[0,T].
 	$$
 	
 	We next prove the nonnegativity of $b(\mathbf{x},t)$.
 	The chemotactic drift term
 	$
 	-\chi_0\sigma\nabla\cdot
 	\left(
 	\frac{nb^2}{(1+n)^2}\nabla n
 	\right)
 	$
 	prevents a direct negative-part test with $b_-$, since it leads to mixed terms that are not readily closed. To overcome this difficulty, we introduce the transformation
 	\begin{align}
 		\label{b_equivalent}
 		u:=b\,e^{-\frac{\chi_0 n}{1+n}}.
 	\end{align}
 	Since $n\ge 0$, we have $1+n>0$. Moreover, the exponential factor is strictly positive, and hence proving $u\ge 0$ is equivalent to proving $b\ge 0$.
 	
 	Obviously, \eqref{b_equivalent} implies
 	\begin{align*}
 		\nabla b
 		=
 		\nabla\left(e^{\frac{\chi_0 n}{1+n}}u\right)
 		=
 		e^{\frac{\chi_0 n}{1+n}}\nabla u
 		+
 		\frac{\chi_0}{(1+n)^2}
 		e^{\frac{\chi_0 n}{1+n}}u\nabla n
 		=
 		e^{\frac{\chi_0 n}{1+n}}
 		\left(
 		\nabla u
 		+
 		\frac{\chi_0}{(1+n)^2}u\nabla n
 		\right),
 	\end{align*}
 	which gives
 	$$
 	\nabla b
 	-
 	\frac{\chi_0 b}{(1+n)^2}\nabla n
 	=
 	e^{\frac{\chi_0 n}{1+n}}\nabla u.
 	$$
 	Since $\sigma$ is spatially constant for each realization, we have
 	\begin{align*}
 		&\sigma nb\nabla b
 		-
 		\frac{\chi_0\sigma nb^2}{(1+n)^2}\nabla n
 		=
 		\sigma nb
 		\left(
 		\nabla b
 		-
 		\frac{\chi_0 b}{(1+n)^2}\nabla n
 		\right)
 		=
 		\sigma n
 		e^{\frac{2\chi_0 n}{1+n}}
 		u\nabla u.
 	\end{align*}
 	Hence, equation \eqref{b} can be rewritten as
 	\begin{align}
 		\label{b_t_equal}
 		b_t
 		=
 		\nabla\cdot
 		\left(
 		\sigma n
 		e^{\frac{2\chi_0 n}{1+n}}
 		u\nabla u
 		\right)
 		+
 		ne^{\frac{\chi_0 n}{1+n}}u.
 	\end{align}
 	On the other hand, from \eqref{b_equivalent}, we have
 	\begin{align}
 		\label{b_positivity_1}
 		b_t
 		&=
 		\partial_t
 		\left(
 		e^{\frac{\chi_0 n}{1+n}}u
 		\right)
 		=
 		e^{\frac{\chi_0 n}{1+n}}
 		\left(
 		u_t
 		+
 		\frac{\chi_0 n_t}{(1+n)^2}u
 		\right),
 	\end{align}
 	which, combined with \eqref{b_t_equal}, yields
 	\begin{align}
 		\label{b_positivity_2}
 		e^{\frac{\chi_0 n}{1+n}}
 		\left(
 		u_t
 		+
 		\frac{\chi_0 n_t}{(1+n)^2}u
 		\right)
 		=
 		\nabla\cdot
 		\left(
 		\sigma n
 		e^{\frac{2\chi_0 n}{1+n}}
 		u\nabla u
 		\right)
 		+
 		ne^{\frac{\chi_0 n}{1+n}}u.
 	\end{align}
 	
 	Moreover, by the homogeneous Neumann boundary conditions
 	\eqref{boundary_condition}, we have
 	\begin{align*}
 		\frac{\partial u}{\partial \boldsymbol{\nu}}
 		&=
 		\frac{\partial}{\partial \boldsymbol{\nu}}
 		\left(
 		b\,e^{-\frac{\chi_0 n}{1+n}}
 		\right)
 		=
 		e^{-\frac{\chi_0 n}{1+n}}
 		\left(
 		\frac{\partial b}{\partial \boldsymbol{\nu}}
 		-
 		\frac{\chi_0 b}{(1+n)^2}
 		\frac{\partial n}{\partial \boldsymbol{\nu}}
 		\right)
 		=0
 	\end{align*}
 	on $\partial\Omega\times(0,T]$.
 	
 	Now we define the negative part of $u(\mathbf{x},t)$ by
 	\begin{align}
 		\label{u_negative_part}
 		u_-(\mathbf{x},t):=\max\{-u(\mathbf{x},t),0\}\geq 0,
 	\end{align}
 	and denote
 	\begin{align}
 		\Omega_-(t):=
 		\{\mathbf{x}\in\Omega:\ u(\mathbf{x},t)<0\}.
 	\end{align}
 	Then $u_-\geq 0$ in $\Omega$. Moreover, $u_-$ vanishes outside
 	$\Omega_-(t)$, and all integral terms involving $u_-$ are therefore supported on $\Omega_-(t)$.
 	
 	On $\Omega_-(t)$, we have
 	\begin{align}
 		u=-u_-,
 		\qquad
 		\nabla u=-\nabla u_-,
 		\qquad
 		u_t=-(u_-)_t .
 	\end{align}
 	Hence,
 	\begin{align}
 		\label{b_positivity_3}
 		\partial_t
 		\left(
 		e^{\frac{\chi_0 n}{1+n}}u
 		\right)
 		=
 		-\partial_t
 		\left(
 		e^{\frac{\chi_0 n}{1+n}}u_-
 		\right),
 		\quad
 		(\mathbf{x},t)\in\Omega_-(t)\times(0,T].
 	\end{align}
 	Since $u$ is continuous, the free boundary of the negative set satisfies
 	\begin{align}
 		\label{zero_boundary}
 		u=0,
 		\qquad
 		u_-=0,
 		\qquad
 		\text{on } \partial\Omega_-(t)\cap\Omega .
 	\end{align}
 	
 	By integrating equation \eqref{b_positivity_3} over $\Omega_-(t)$
 	and using equations \eqref{b_positivity_1} and \eqref{b_positivity_2},
 	we get
 	\begin{align*}
 		\frac{d}{dt}
 		\int_{\Omega_-}
 		e^{\frac{\chi_0 n}{1+n}}u_-\,d\mathbf{x}
 		&=
 		-\int_{\Omega_-}
 		\partial_t
 		\left(
 		e^{\frac{\chi_0 n}{1+n}}u
 		\right)
 		\,d\mathbf{x}
 		\\
 		&=
 		-\int_{\Omega_-}
 		\Biggl[
 		\nabla\cdot
 		\left(
 		\sigma n
 		e^{\frac{2\chi_0 n}{1+n}}
 		u\nabla u
 		\right)
 		+
 		ne^{\frac{\chi_0 n}{1+n}}u
 		\Biggr]\,d\mathbf{x}
 		\\
 		&=
 		-\int_{\Omega_-}
 		\nabla\cdot
 		\left(
 		\sigma n
 		e^{\frac{2\chi_0 n}{1+n}}
 		u\nabla u
 		\right)
 		\,d\mathbf{x}
 		+
 		\int_{\Omega_-}
 		ne^{\frac{\chi_0 n}{1+n}}u_-\,d\mathbf{x}.
 	\end{align*}
 	The boundary contribution associated with the divergence term vanishes under the zero-flux condition \eqref{zero_boundary}. Therefore,
 	\begin{align*}
 		&-\int_{\Omega_-}
 		\nabla\cdot
 		\left(
 		\sigma n
 		e^{\frac{2\chi_0 n}{1+n}}
 		u\nabla u
 		\right)
 		\,d\mathbf{x}
 		=
 		-\int_{\partial\Omega_-}
 		\sigma n
 		e^{\frac{2\chi_0 n}{1+n}}
 		u\,
 		\frac{\partial u}{\partial \boldsymbol{\nu}}
 		\,ds
 		=0.
 	\end{align*}
   Thus, by noting \eqref{u_negative_part}, we arrive at
 	\begin{align}
 		\label{dt_u_negative_part}
 		\frac{d}{dt}
 		\int_{\Omega_-}
 		e^{\frac{\chi_0 n}{1+n}}u_-\,d\mathbf{x}
 		=
 		\int_{\Omega_-}
 		ne^{\frac{\chi_0 n}{1+n}}u_-\,d\mathbf{x}
 		&\le
 		\|n(\cdot,t)\|_{L^\infty(\Omega)}
 		\int_{\Omega_-}
 		e^{\frac{\chi_0 n}{1+n}}u_-\,d\mathbf{x}.
 	\end{align}
 	
 	Since $b_0\ge 0$, it follows from
 	$$
 	u(\mathbf{x},0)
 	=
 	b_0
 	e^{-\frac{\chi_0 n_0}{1+n_0}}
 	\ge 0
 	$$
 	that
 	$$
 	u_-(\mathbf{x},0)\equiv0.
 	$$
 	Hence,
 	$$
 	\int_{\Omega_-(0)}
 	e^{\frac{\chi_0 n_0(\mathbf{x})}
 		{1+n_0(\mathbf{x})}}
 	u_-(\mathbf{x},0)\,d\mathbf{x}
 	=0.
 	$$
 	By Gr\"onwall's inequality, we conclude from
 	\eqref{dt_u_negative_part} that
 	$$
 	\int_{\Omega_-(t)}
 	e^{\frac{\chi_0 n}{1+n}}u_-\,d\mathbf{x}
 	=0,
 	\qquad \forall  t\in[0,T].
 	$$
 	Since
 	$$
 	e^{\frac{\chi_0 n}{1+n}}>0,
 	$$
 	this implies
 	$$
 	u_-(\mathbf{x},t)\equiv 0
 	\qquad \text{in } \Omega\times[0,T].
 	$$
 	Therefore,
 	$$
 	u(\mathbf{x},t)\ge 0
 	\qquad \text{in } \Omega\times[0,T],
 	$$
 	which yields
 	$$
 	b(\mathbf{x},t)\ge 0
 	\qquad \text{in } \Omega\times[0,T].
 	$$
 	This completes the proof.
 \end{proof}
\section{Structure-preserving schemes}	
\label{3_structure_preserving}	
In this section, we study the structure-preserving scheme for the equations
\eqref{system_4}-\eqref{boundary_condition}. After introducing the high-order finite element method with mass-lumping  in space, we propose the  fully discrete schemes with first- and second-order convergence in time, which are proved to be unconditional total mass conservation and nonnegativity-preservation.
\subsection{Finite element with mass-lumping}
Mass-lumping finite element methods have been developed over several decades for the numerical approximation of time-dependent partial differential equations \cite{GPJ2001}. For linear finite elements, mass lumping can be realized by employing a suitable nodal quadrature rule without modifying the underlying finite element space. The situation is more subtle for quadratic elements. The standard $\mathbb{P}^2$ finite element with exact integration leads to a symmetric positive definite but non-diagonal mass matrix, and hence does not provide a nodal discrete inner product  form required by the pointwise KKT correction which will be used to preserve the positivity in this paper. On the other hand, under the required quadrature accuracy,  directly applying mass lumping at the six nodes of the standard $\mathbb{P}^2$  element results in zero weights at the vertices, yielding a singular lumped mass matrix. To overcome this difficulty, the standard quadratic polynomial space is enriched by suitable bubble functions so that a mass-lumped element with strictly positive quadrature weights and sufficient accuracy can be constructed \cite{GPJ2001}. Such mass-lumped constructions have also been incorporated into structure-preserving finite element methods for nonlinear parabolic equations \cite{yj2023}. In the present work, mass lumping is employed not only to diagonalize the mass matrix, but, more importantly, to provide a nodal discrete inner product with strictly positive quadrature weights. This structure is well suited to the subsequent KKT-Lagrange multiplier correction and preserves the monotonicity of the scalar mass-correction equation, which is essential for enforcing nodal nonnegativity and mass conservation simultaneously.

To derive the fully discrete scheme for equations
\eqref{system_4}-\eqref{boundary_condition}, we first introduce the spatial discretization subject to homogeneous Neumann boundary conditions on the computational domain $\Omega$. Let $\mathcal{T}_h$ be a shape-regular and quasi-uniform simplicial triangulation of $\Omega$. Given a simplex $K\in\mathcal{T}_h$, we use $\mathbb{P}^r(K)$ to represent the polynomial space on $K$ with total degree not exceeding $r$. For the linear mass-lumping element, we take
$
\widetilde{\mathbb{P}}^1(K)=\mathbb{P}^1(K),
$
whereas for the quadratic case, the standard space $\mathbb{P}^2(K)$ is enriched by suitable bubble functions to construct $\widetilde{\mathbb{P}}^2(K)$. Accordingly, we define
\begin{align}
	V_h^r
	=
	\left\{
	v_h\in H^1(\Omega):
	v_h|_K\in\widetilde{\mathbb{P}}^r(K),
	\quad \forall K\in\mathcal{T}_h
	\right\}, \qquad r=1,2.
\end{align}
The associated nodal quadrature rules have strictly positive weights and sufficient algebraic accuracy, leading to diagonal positive definite mass matrices. For completeness, Table \ref{tab:poly_space_details} summarizes the two-dimensional local spaces $\widetilde{\mathbb{P}}^1$ and $\widetilde{\mathbb{P}}^2$ used in this work.

\begin{table}[htbp]
	\centering
	\caption{Basis functions, interpolation nodes, and lumped quadrature weights associated with the finite element spaces $\widetilde{\mathbb{P}}^1$ and $\widetilde{\mathbb{P}}^2$ \cite{GPJ2001, yj2023}.}
	\label{tab:poly_space_details}
	\renewcommand{\arraystretch}{1.15}
	\begin{tabular}{c p{5.2cm} p{3.6cm} l}
		\toprule
		Finite element space & Shape functions & Barycentric node & Lumped quadrature weight \\
		\midrule
		$\widetilde{\mathbb{P}}^1$
		& $\mu_1,\ \mu_2,\ \mu_3$
		& $(1,0,0),\ (0,1,0),\ (0,0,1)$
		&\ \ \  $\frac{1}{3},\ \frac{1}{3},\ \frac{1}{3}$ \\
		\cline{1-4}
		\multirow{7}{*}{$\widetilde{\mathbb{P}}^2$}
		& $\mu_1(2\mu_1-1)+3\mu_1\mu_2\mu_3$
		& $(1,0,0)$
		&\ \ \ $\frac{1}{20}$ \\
		& $\mu_2(2\mu_2-1)+3\mu_1\mu_2\mu_3$
		& $(0,1,0)$
		& \ \ \ $\frac{1}{20}$ \\
		& $\mu_3(2\mu_3-1)+3\mu_1\mu_2\mu_3$
		& $(0,0,1)$
		&\ \ \ $\frac{1}{20}$ \\
		& $4\mu_2\mu_3-12\mu_1\mu_2\mu_3$
		& $(0,\frac{1}{2},\frac{1}{2})$
		& \ \ \ $\frac{2}{15}$ \\
		& $4\mu_3\mu_1-12\mu_1\mu_2\mu_3$
		& $(\frac{1}{2},0,\frac{1}{2})$
		& \ \ \ $\frac{2}{15}$ \\
		& $4\mu_1\mu_2-12\mu_1\mu_2\mu_3$
		& $(\frac{1}{2},\frac{1}{2},0)$
		&\ \ \ $\frac{2}{15}$ \\
		& $27\mu_1\mu_2\mu_3$
		& $(\frac{1}{3},\frac{1}{3},\frac{1}{3})$
		&\ \ \ $\frac{9}{20}$ \\
		\bottomrule
	\end{tabular}
	
	\vspace{0.3cm}
	\parbox{0.92\textwidth}{\small
		\textbf{Note:} Here, $\mu_i$ denotes the $i$-th area coordinate on a triangular element, which equals 1 at the $i$-th vertex and 0 at the other two vertices.
	}
\end{table}

Let $\Sigma_h=\{\mathbf{x}_i\}_{i=1}^N$ denote the set of quadrature points associated with the numerical integration rule, where $\mathbf{x}_i\in \bar{\Omega}$ is the $i$-th quadrature point and $N$ is the grid point. We define the discrete inner product by
$$
(f,g)_h:=\sum_{i=1}^N \omega_i f(\mathbf{x}_i)g(\mathbf{x}_i),
$$
and the induced discrete norm by
$$
\|f_h\|_h:=\sqrt{(f_h,f_h)_h}, \qquad \forall f_h\in V_h^r, \quad r=1, 2.
$$
where $\omega_i$ is the corresponding quadrature weight.

\subsection{First-order fully discrete scheme}
Based on the mass-lumping finite element method and the corresponding discrete inner product introduced in the last subsection, we next construct the first-order (BDF1) structure-preserving fully discrete scheme.

Let $\tau$ denote the temporal step size, and set
$t^m=m\tau$ for $m=0,1,\ldots,M$, where $M=T/\tau$.
The spatial approximation is based on a mass-lumping finite element method, whereas a predictor-corrector (IMEX) strategy is adopted for time integration. To enforce the desired structural properties, the KKT shift and Lagrange multiplier correction are incorporated into the resulting scheme \cite{CS2022I,CS2022II,FQ2024}. More precisely, the first-order fully discrete scheme is given as follows.

\medskip
\noindent\textbf{Scheme 1.}
Set $\lambda_h^0=\eta_h^0=0$. Given $n_h^0=\Pi_h n_0$ and $b_h^0=\Pi_h b_0$, where $\Pi_h$ denotes the interpolation operator, we compute $n_h^{m+1}$ and $b_h^{m+1}$ for $m=0,1,\ldots,M-1$ as follows:

\noindent\textbf{Step 1.} Find $\tilde n_h^{m+1}\in V_h^r$ such that,
\begin{align}\label{first_order_predictor_n}
	\left(
	\frac{\tilde n_h^{m+1}-n_h^m}{\tau}, \phi_h
	\right)_h
	+
	(\nabla \tilde n_h^{m+1}, \nabla\phi_h)_h
	+
	(b_h^m\tilde n_h^{m+1}, \phi_h)_h
	=0, \quad \forall \phi_h\in V_h^r.
\end{align}

\noindent\textbf{Step 2.} Find $\tilde b_h^{m+1}\in V_h^r$ such that,
\begin{align}\label{first_order_predictor_b}
	\left(
	\frac{\tilde b_h^{m+1}-b_h^m}{\tau}, \psi_h
	\right)_h
	+
	(\sigma n_h^m b_h^m\nabla\tilde b_h^{m+1},\nabla\psi_h)_h
	=
	(\tilde n_h^{m+1}b_h^m, \psi_h)_h
	+
	\chi_0\left(\frac{\sigma n_h^m(b_h^m)^2}{(1+n_h^m)^2}\nabla n_h^m, \nabla\psi_h\right)_h, \quad \forall \psi_h\in V_h^r.
\end{align}

\noindent\textbf{Step 3.} Find $\xi_h^{m+1}\in\mathbb{R}$ from the scalar nonlinear algebraic equation:
\begin{align}\label{first_order_xi_equation}
	F(\xi_h^{m+1})&:=
	\sum_{i=1}^{N}
	\left[
	\left(\tilde n_h^{m+1}(\mathbf{x}_i)+\tau\theta\xi_h^{m+1}\right)_+\right]\omega_i
	+\sum_{i=1}^{N}\left[
	\left(\tilde b_h^{m+1}(\mathbf{x}_i)+\tau(1-\theta)\xi_h^{m+1}\right)_+
	\right]\omega_i
	-
	(n_h^m+b_h^m,1)_h
	=0,
\end{align}
where $(s)_+=\max\{s,0\}$ and $\theta \in (0, 1)$.

\noindent\textbf{Step 4.} Find $n_h^{m+1}$ and $\lambda_h^{m+1}$ from
\begin{subequations}\label{first_order_correction_n_system}
	\begin{align}
		&\frac{n_h^{m+1}(\mathbf{x}_i)-\tilde n_h^{m+1}(\mathbf{x}_i)}{\tau}
		=
		\lambda_h^{m+1}(\mathbf{x}_i)+\theta\xi_h^{m+1},
		\qquad \forall \mathbf{x}_i\in\Sigma_h,
		\label{first_order_correction_n_eq}
		\\
		&\lambda_h^{m+1}(\mathbf{x}_i)\ge0,\qquad
		n_h^{m+1}(\mathbf{x}_i)\ge0,\qquad
		\lambda_h^{m+1}(\mathbf{x}_i)n_h^{m+1}(\mathbf{x}_i)=0,
		\qquad \forall \mathbf{x}_i\in\Sigma_h.
		\label{first_order_kkt_n}
	\end{align}
\end{subequations}

\noindent\textbf{Step 5.} Find $b_h^{m+1}$ and $\eta_h^{m+1}$ from
\begin{subequations}\label{eq:first_order_correction_b_system}
	\begin{align}
		&\frac{b_h^{m+1}(\mathbf{x}_i)-\tilde b_h^{m+1}(\mathbf{x}_i)}{\tau}
		=
		\eta_h^{m+1}(\mathbf{x}_i)+(1-\theta)\xi_h^{m+1},
		\qquad \forall \mathbf{x}_i\in\Sigma_h,
		\label{first_order_correction_b_eq}
		\\
		&\eta_h^{m+1}(\mathbf{x}_i)\ge0,\qquad
		b_h^{m+1}(\mathbf{x}_i)\ge0,\qquad
		\eta_h^{m+1}(\mathbf{x}_i)b_h^{m+1}(\mathbf{x}_i)=0,
		\qquad \forall \mathbf{x}_i\in\Sigma_h.
		\label{first_order_kkt_b}
	\end{align}
\end{subequations}

From \textbf{Scheme 1}, we can see that the predicted solutions $\tilde n_h^{m+1}$ and $\tilde b_h^{m+1}$ in \textbf{Steps 1} and \textbf{2} are obtained from standard time-stepping schemes and may not preserve nonnegativity or the total mass. Therefore, a global scalar mass-correction multiplier $\xi_h^{m+1}$ is introduced in \textbf{Step 3} to enforce the conservation of the total discrete mass, and nodal Lagrange multipliers $\lambda_h^{m+1}$ and $\eta_h^{m+1}$ are used to impose the nonnegativity constraints on $n_h^{m+1}$ and $b_h^{m+1}$ in \textbf{Steps 4} and \textbf{5}, respectively.
In fact, after $\xi_h^{m+1}$ is computed in \textbf{Step 3}, \textbf{Steps 4} and \textbf{5} admit the pointwise solution:
\begin{align}\label{eq:first_order_pointwise_n}
	(n_h^{m+1}(\mathbf{x}_i),\lambda_h^{m+1}(\mathbf{x}_i))
	=
	\begin{cases}
		\left(
		\tilde n_h^{m+1}(\mathbf{x}_i)+\tau\theta\xi_h^{m+1}, \,0
		\right),
		&
		\tilde n_h^{m+1}(\mathbf{x}_i)+\tau\theta\xi_h^{m+1}>0,
		\\[6pt]
		\left(
		0,\,
		-\dfrac{\tilde n_h^{m+1}(\mathbf{x}_i)+\tau\theta\xi_h^{m+1}}{\tau}
		\right),
		&
		\tilde n_h^{m+1}(\mathbf{x}_i)+\tau\theta\xi_h^{m+1}\le0,
	\end{cases}
\end{align}
and
\begin{align}\label{first_order_pointwise_b}
	(b_h^{m+1}(\mathbf{x}_i),\eta_h^{m+1}(\mathbf{x}_i))
	=
	\begin{cases}
		\left(
		\tilde b_h^{m+1}(\mathbf{x}_i)+\tau(1-\theta)\xi_h^{m+1},\,0
		\right),
		&
		\tilde b_h^{m+1}(\mathbf{x}_i)+\tau(1-\theta)\xi_h^{m+1}>0,
		\\[6pt]
		\left(
		0,\,
		-\dfrac{\tilde b_h^{m+1}(\mathbf{x}_i)+\tau(1-\theta)\xi_h^{m+1}}{\tau}
		\right),
		&
		\tilde b_h^{m+1}(\mathbf{x}_i)+\tau(1-\theta)\xi_h^{m+1}\le0.
	\end{cases}
\end{align}

In fact, \eqref{eq:first_order_pointwise_n}-\eqref{first_order_pointwise_b} suggest that
$$
n_h^{m+1}(\mathbf{x}_i)
=
\left(
\tilde n_h^{m+1}(\mathbf{x}_i)
+\tau\theta\xi_h^{m+1}
\right)_+,
\qquad
b_h^{m+1}(\mathbf{x}_i)
=
\left(
\tilde b_h^{m+1}(\mathbf{x}_i)
+\tau(1-\theta)\xi_h^{m+1}
\right)_+.
$$
Putting these into \eqref{first_order_xi_equation}, we have
\begin{align}\label{first_order_mass_conservation}
	(n_h^{m+1}+b_h^{m+1},1)_h-(n_h^m+b_h^m,1)_h=0,
\end{align}
which confirms the total mass conservation.
\begin{remark}
	The scalar variable $\xi_h^{m+1}$ is determined by solving the nonlinear algebraic equation \eqref{first_order_xi_equation}. From  Theorem \ref{secant_proof} established in the following, we can see that $\xi_h^{m+1}$ can be calculated exactly by the secant iteration if the initial guesses belongs to a sufficiently small neighborhood of the solution, which ensures the mass conservation  \eqref{first_order_mass_conservation}. So does the second-order fully discrete scheme considered in the next subsection.
\end{remark}
Moreover, for \textbf{Scheme 1}, it is obviously that only linear systems need to be solved in \textbf{Steps 1} and \textbf{2}, and these are assignment operations in \textbf{Steps 4} and \textbf{5}. The remaining issue is the unique solvability of the scalar nonlinear algebraic equation \eqref{first_order_xi_equation} in \textbf{Step 3}, as established in the following theorem.
\begin{theorem}[Unique solvability of \eqref{first_order_xi_equation}]
\label{first_order_solvability}
Under the Assumption \ref{regularity}, if $\theta\in(0,1)$ and $(n_h^0+b_h^0, 1)_h>0$, then the nonlinear algebraic equation \eqref{first_order_xi_equation} generated from  the mass-correction multiplier admits a unique solution $\xi^{m+1}_{h}\in\mathbb{R}$ satisfying
	$$\xi^{m+1}_{h}\leq0.$$
\end{theorem}

\begin{proof}
We will prove the result by induction. It suffices to show that
\eqref{first_order_xi_equation} admits a unique solution, provided that the first step is uniquely solvable. By the assumption on the initial data and the mass conservation constraint \eqref{first_order_mass_conservation}, we obtain that
$$(n_h^m+b_h^m, 1)_h=(n_h^0+b_h^0, 1)_h>0.$$

From \eqref{first_order_xi_equation}, we have that $F$ is continuous on $\mathbb{R}$ since it is a finite sum of positive-part functions. Moreover, noting $\theta\in(0,1)$, we have
$$
\lim_{\xi\to-\infty}F(\xi_{h}^{m+1})=-(n_h^m+b_h^m, 1)_h<0,
\qquad
\lim_{\xi\to+\infty}F(\xi_{h}^{m+1})=+\infty.
$$
Hence, by the intermediate value theorem, there exists at least one $\xi^{m+1}_{h}\in\mathbb{R}$ such that
$$
F(\xi^{m+1}_{h})=0.
$$
	
We next prove the uniqueness. For each node $\mathbf{x}_i$, the functions
$$
\left[\tilde n^{m+1}_{h}(\mathbf{x}_i)+\tau\theta\xi_{h}^{m+1}\right]_{+}
\quad\text{and}\quad
\left[\tilde b^{m+1}_{h}(\mathbf{x}_i)+\tau(1-\theta)\xi_{h}^{m+1}\right]_{+}
$$
are nondecreasing with respect to $\xi_h^{m+1}$. Hence $F$ is nondecreasing on $\mathbb{R}$.

Suppose that $F$ has two zero points $\xi_1$ and $\xi_2$ with $\xi_1<\xi_2$. Since $F(\xi_1)=0$, we have
$$
\sum_{i=1}^{N}\left[\tilde n^{m+1}_{h}(\mathbf{x}_i)+\tau\theta\xi_1\right]_{+}\omega_i
+
\sum_{i=1}^{N}\left[\tilde b^{m+1}_{h}(\mathbf{x}_i)+\tau(1-\theta)\xi_1\right]_{+}\omega_i
=(n_h^m+b_h^m, 1)_h>0.
$$
Therefore, at least one of the above two positive-part term is strictly positive at $\xi_1$. Since $\theta\in(0,1)$, this active term is strictly larger at $\xi_2$ than at $\xi_1$, while all other terms are nondecreasing. Consequently,
$$
F(\xi_2)>F(\xi_1)=0.
$$
This contradicts the assumption that $\xi_2$ is also a zero point of $F$. Hence $F$ cannot have two distinct zeros, and the solution of \eqref{first_order_xi_equation} is unique.
	
It remains to prove the sign of the corresponding multiplier. Let $\xi^{m+1}_{h}$ be the unique solution of \eqref{first_order_xi_equation}, and let $n^{m+1}_{h}$, $b^{m+1}_{h}$, $\lambda^{m+1}_{h}$, and $\eta^{m+1}_{h}$ be determined by the nodal correction relations. Adding \eqref{first_order_predictor_n} and \eqref{first_order_correction_n_eq},
	\begin{align}
		\left(\dfrac{n^{m+1}_{h}-n^{m}_{h}}{\tau},\phi_h\right)_{h}
		+(\nabla \tilde n^{m+1}_{h},\nabla \phi_h)_{h}
		+(b^{m}_{h}\tilde n^{m+1}_{h},\phi_h)_{h}
		=(\lambda^{m+1}_{h}+\theta\xi^{m+1}_{h},\phi_h)_{h}, \qquad \forall \phi_h\in V^{r}_{h}
		\label{eq:n_combined_correction}
	\end{align}
	Similarly, adding \eqref{first_order_predictor_b} and \eqref{first_order_correction_b_eq} yields,
	\begin{align}
		\left(\dfrac{b^{m+1}_{h}-b^{m}_{h}}{\tau},\psi_h\right)_{h}
		+(\sigma n^{m}_{h}b^{m}_{h}\nabla \tilde b^{m+1}_{h},\nabla \psi_h)_{h}
		=&(\tilde n^{m+1}_{h}b^{m}_{h},\psi_h)_{h}
		+\chi_0\left(\frac{\sigma n_h^m(b_h^m)^2}{(1+n_h^m)^2}\nabla n_h^m, \nabla\psi_h\right)_h \notag\\
		&+(\eta^{m+1}_{h}+(1-\theta)\xi^{m+1}_{h},\psi_h)_{h}, \qquad \forall \psi_h\in V^{r}_{h}.
		\label{eq:b_combined_correction}
	\end{align}
	Taking $\phi_h=\psi_h=1$ in \eqref{eq:n_combined_correction} and \eqref{eq:b_combined_correction}, and using $\nabla 1=0$, we obtain
	\begin{align}
		\left(\dfrac{n^{m+1}_{h}-n^{m}_{h}}{\tau},1\right)_{h}
		+(b^{m}_{h}\tilde n^{m+1}_{h},1)_{h}
		=(\lambda^{m+1}_{h}+\theta\xi^{m+1}_{h},1)_{h},
		\label{eq:n_mass_balance_with_multiplier}
	\end{align}
	and
	\begin{align}
		\left(\dfrac{b^{m+1}_{h}-b^{m}_{h}}{\tau},1\right)_{h}
		=(\tilde n^{m+1}_{h}b^{m}_{h},1)_{h}
		+(\eta^{m+1}_{h}+(1-\theta)\xi^{m+1}_{h},1)_{h}.
		\label{eq:b_mass_balance_with_multiplier}
	\end{align}
	Adding \eqref{eq:n_mass_balance_with_multiplier} and \eqref{eq:b_mass_balance_with_multiplier}, and cancelling the common reaction term, gives
	\begin{align}
		\left(\dfrac{n^{m+1}_{h}-n^{m}_{h}}{\tau}
		+\dfrac{b^{m+1}_{h}-b^{m}_{h}}{\tau},1\right)_{h}
		=(\lambda^{m+1}_{h}+\eta^{m+1}_{h}+\xi^{m+1}_{h},1)_{h}.
		\label{eq:total_mass_multiplier_balance}
	\end{align}
	On the other hand, the scalar equation $F(\xi^{m+1}_{h})=0$ is exactly the discrete total mass conservation constraint, namely
	$$
	(n^{m+1}_{h}+b^{m+1}_{h},1)_{h}
	=(n^{m}_{h}+b^{m}_{h},1)_{h}.
	$$
	Hence,
	$$
	\left(\dfrac{n^{m+1}_{h}-n^{m}_{h}}{\tau}
	+\dfrac{b^{m+1}_{h}-b^{m}_{h}}{\tau},1\right)_{h}
	=0.
	$$
	It follows from \eqref{eq:total_mass_multiplier_balance} that
	$$
	(\lambda^{m+1}_{h}+\eta^{m+1}_{h}+\xi^{m+1}_{h},1)_{h}=0.
	$$
	Since $\xi^{m+1}_{h}$ is a scalar multiplier, we have
	$$
	\xi^{m+1}_{h}
	=-\dfrac{(\lambda^{m+1}_{h}+\eta^{m+1}_{h},1)_{h}}{(1,1)_{h}}.
	$$
	Finally, using \eqref{first_order_kkt_b} and \eqref{first_order_kkt_n}, we obtain
	$$
	\lambda^{m+1}_{h}(\mathbf{x}_i)\geq0,
	\qquad
	\eta^{m+1}_{h}(\mathbf{x}_i)\geq0,
	\qquad i=1,\ldots,N.
	$$
	Therefore,
	$$
	(\lambda^{m+1}_{h}+\eta^{m+1}_{h},1)_{h}\geq0,
	$$
	which implies
	$$
	\xi^{m+1}_{h}\leq0.
	$$
	This completes the proof.
\end{proof}
Next, we consider the iteration method of the nonlinear algebraic equation \eqref{first_order_xi_equation}. Since $F'(\xi)$ may be unavailable or difficult to evaluate, we solve \eqref{first_order_xi_equation} by the secant iteration
\begin{align}\label{first_order_secant}
	\xi^{(\ell+1)}
	&=
	\xi^{(\ell)}
	-
	\frac{
		F(\xi^{(\ell)})(\xi^{(\ell)}-\xi^{(\ell-1)})
	}{
		F(\xi^{(\ell)})-F(\xi^{(\ell-1)})
	},
	\qquad \ell= 1, 2, \cdots,
\end{align}
with $\xi^{(0)}=0$ and
$\xi^{(1)}=-c_0\tau$, where $c_0>0$, which belongs to small enough neighborhood of $\xi_h^{m+1}$.

Based on the uniqueness of $\xi_h^{m+1}$ established above, we next show that the secant iteration \eqref{first_order_secant} used for solving \eqref{first_order_xi_equation} is locally convergent.
\begin{theorem}[Local convergence of \eqref{first_order_secant}]
	\label{secant_proof}
	Let $\xi_h^{m+1}$ be the solution of \eqref{first_order_xi_equation}. The secant iteration \eqref{first_order_secant} converges locally to $\xi_h^{m+1}$, and there exists $\ell$, such that the iteration solution $\xi^{(\ell)}=\xi_h^{m+1}$. Moreover, once two distinct consecutive iterations belongs to a sufficiently small neighborhood of $\xi_h^{m+1}$, the secant iteration \eqref{first_order_secant} reaches
	$\xi_h^{m+1}$ after at most three additional updates.
\end{theorem}

\begin{proof}
	Since the positive-part function is continuous, nondecreasing, convex, and piecewise affine, each term
	$$
	\left[
	\tilde n_h^{m+1}(\mathbf{x}_i)+\tau\theta\xi_{h}^{m+1}
	\right]_+
	\quad \text{and} \quad
	\left[
	\tilde b_h^{m+1}(\mathbf{x}_i)+\tau(1-\theta)\xi_{h}^{m+1}
	\right]_+
	$$
	has the same properties as a function of $\xi_{h}^{m+1}$. Therefore, $F$ is
	continuous, nondecreasing, convex, and piecewise affine.
	
    Let $\xi^{(\ell-1)}$ and $\xi^{(\ell)}$ be two distinct consecutive iterates in a sufficiently small neighborhood
    $U_\delta=(\xi_h^{m+1}-\delta,\xi_h^{m+1}+\delta)$
    of $\xi_h^{m+1}$.

	\textit{Case 1: $\xi_h^{m+1}$ is not a breakpoint of $F$.}
	In this case, $\delta$ can be chosen so that $F$ has no breakpoint
	in $U_\delta$. Therefore, $F$ is affine throughout $U_\delta$, and
	$$
	F(\xi)
	=
	F'(\xi_h^{m+1})
	\left(
	\xi-\xi_h^{m+1}
	\right),
	\qquad
	\xi\in U_\delta,
	$$
	where $F'(\xi_h^{m+1})>0$. Since $F$ is continuous, nondecreasing and $\xi_h^{m+1}$ is the unique solution (see Theorem \ref{first_order_solvability}).
	
	The two consecutive iterates lie on the same affine segment of $F$.
	Hence, the secant line through the corresponding points on the graph
	of $F$ coincides with $F$, and its zero is $\xi_h^{m+1}$. Thus,
	$$
	\xi^{(\ell+1)}
	=
	\xi_h^{m+1}.
	$$
	
    \textit{Case 2: $\xi_h^{m+1}$ is a breakpoint of $F$.}
    Then $\xi_{h}^{m+1}$ satisfies
	$$
	\tilde n_h^{m+1}(\mathbf{x}_i)
	+\tau\theta\xi_{h}^{m+1}
	=0
	\quad \text{or} \quad
    \tilde b_h^{m+1}(\mathbf{x}_i)
	+\tau(1-\theta)\xi_{h}^{m+1}
	=0.
	$$
	Since there are only finitely nodal values, $F$ has finitely breakpoints. Hence, there exists $\delta>0$ such that
	$
	U_\delta
	=
	\left(
	\xi_h^{m+1}-\delta,\,
	\xi_h^{m+1}+\delta
	\right)
	$
	contains no breakpoint of $F$ other than possibly
	$\xi_h^{m+1}$.
	
	We next determine the one-sided derivatives of $F$ at
	$\xi_h^{m+1}$. Consider a generic positive-part term
	$$
	g(\xi_h^{m+1})=[a+\gamma\xi_h^{m+1}]_+,
	\qquad \gamma>0.
	$$
	
	If
	$$
	a+\gamma\xi_h^{m+1}>0,
	$$
	then $g$ is affine on both sides of $\xi_h^{m+1}$, and this term
	contributes $\gamma$ to both one-sided derivatives.
	
	If
	$$
	a+\gamma\xi_h^{m+1}<0,
	$$
	then $g$ vanishes in a neighborhood of $\xi_h^{m+1}$ and
	contributes zero to both one-sided derivatives.
	
	If
	$$
	a+\gamma\xi_h^{m+1}=0,
	$$
	then, for every sufficiently small $\varepsilon>0$,
	$$
	a+\gamma(\xi_h^{m+1}-\varepsilon)
	=
	-\gamma\varepsilon<0,
	$$
	whereas
	$$
	a+\gamma(\xi_h^{m+1}+\varepsilon)
	=
	\gamma\varepsilon>0.
	$$
	Consequently,
	$$
	g'_{-}(\xi_h^{m+1})=\lim_{\varepsilon\to 0} \frac{g(\xi_h^{m+1})-g(\xi_h^{m+1}-\varepsilon)}{\varepsilon}=0,
	\qquad
	g'_{+}(\xi_h^{m+1})=\lim_{\varepsilon\to 0} \frac{g(\xi_h^{m+1}+\varepsilon)-g(\xi_h^{m+1})}{\varepsilon}=\gamma.
	$$
	Thus, a term whose argument vanishes at $\xi_h^{m+1}$ does not contribute to the left derivative but does contribute to the right derivative. This explains why the strict inequality is used in the expression for the left derivative, whereas the nonstrict inequality
	is used for the right derivative.
	
	It follows that
	\begin{align}
		F'_{-}(\xi_h^{m+1})
		=\lim_{\varepsilon\to 0} \frac{F(\xi_h^{m+1})-F(\xi_h^{m+1}-\varepsilon)}{\varepsilon}=\tau\theta
		\sum_{
			\tilde n_h^{m+1}(\mathbf{x}_i)
			+\tau\theta\xi_h^{m+1}>0
		}
		\omega_i
		+
		\tau(1-\theta)
		\sum_{
			\tilde b_h^{m+1}(\mathbf{x}_i)
			+\tau(1-\theta)\xi_h^{m+1}>0
		}
		\omega_i,
		\label{eq:left_derivative_F}
	\end{align}
	whereas
	\begin{align}
		F'_{+}(\xi_h^{m+1})
		=\lim_{\varepsilon\to 0} \frac{F(\xi_h^{m+1}+\varepsilon)-F(\xi_h^{m+1})}{\varepsilon}=
		\tau\theta
		\sum_{
			\tilde n_h^{m+1}(\mathbf{x}_i)
			+\tau\theta\xi_h^{m+1}\ge0
		}
		\omega_i
		+
		\tau(1-\theta)
		\sum_{
			\tilde b_h^{m+1}(\mathbf{x}_i)
			+\tau(1-\theta)\xi_h^{m+1}\ge0
		}
		\omega_i.
		\label{eq:right_derivative_F}
	\end{align}
	
	Since $F(\xi_h^{m+1})=0$, equation
	\eqref{first_order_xi_equation} gives
	\begin{align}
		\sum_{i=1}^{N}
		\left[
		\tilde n_h^{m+1}(\mathbf{x}_i)
		+\tau\theta\xi_h^{m+1}
		\right]_+\omega_i
		+
		\sum_{i=1}^{N}
		\left[
		\tilde b_h^{m+1}(\mathbf{x}_i)
		+\tau(1-\theta)\xi_h^{m+1}
		\right]_+\omega_i
		=
		(n_h^m+b_h^m,1)_h.
		\label{eq:positive_mass_at_root}
	\end{align}
    By the mass conservation property
   \eqref{first_order_mass_conservation}
    and the assumption $(n_h^0+b_h^0,1)_h>0$
    in Theorem \ref{first_order_solvability},
	$$
	(n_h^m+b_h^m,1)_h
	=
	(n_h^0+b_h^0,1)_h
	>0.
	$$
	Therefore, at least one argument of the positive-part terms is
	strictly positive at $\xi_h^{m+1}$. Since $\theta\in (0, 1)$ and
	$\omega_i>0$, we have
	$$
	F'_{-}(\xi_h^{m+1})>0.
	$$
	Moreover, every term contributing to the left derivative also
	contributes to the right derivative, while a term whose argument is
	zero contributes only to the right derivative. Hence,
	\begin{align}
		\label{F_P_F_N}
		0<
		F'_{-}(\xi_h^{m+1})
		\le
		F'_{+}(\xi_h^{m+1}).
	\end{align}

	Since $U_\delta$ contains no other breakpoint, $F$ is affine on each side of $\xi_h^{m+1}$:
	$$
	F(\xi)
	=
	F'_{-}(\xi_h^{m+1})
	\left(
	\xi-\xi_h^{m+1}
	\right),
	\qquad
	\xi_h^{m+1}-\delta<\xi<\xi_h^{m+1},
	$$
	and
	$$
	F(\xi)
	=
	F'_{+}(\xi_h^{m+1})
	\left(
	\xi-\xi_h^{m+1}
	\right),
	\qquad
	\xi_h^{m+1}<\xi<\xi_h^{m+1}+\delta.
	$$
	
	Due to the possible locations of the two consecutive
	iterates.

	\textit{Case 2.1: One iteration equals $\xi_h^{m+1}$},
	i.e.,
	$$
	\xi^{(\ell-1)}=\xi_h^{m+1}
	\qquad\text{or}\qquad
	\xi^{(\ell)}=\xi_h^{m+1},
	$$
	then the iteration has already terminated.

	\textit{Case 2.2: Both iterations lie on the left of
		$\xi_h^{m+1}$}, i.e.,
	$$
	\xi^{(\ell-1)}<\xi_h^{m+1},
	\qquad
	\xi^{(\ell)}<\xi_h^{m+1}.
	$$
	Both iterations lie on the same left affine branch of $F$. Therefore,
	the secant line coincides with this branch, and
	$$
	\xi^{(\ell+1)}
	=
	\xi_h^{m+1}.
	$$

	\textit{Case 2.3: Both iterations lie to the right of
		$\xi_h^{m+1}$}, i.e.,
	$$
	\xi^{(\ell-1)}>\xi_h^{m+1},
	\qquad
	\xi^{(\ell)}>\xi_h^{m+1}.
	$$
	Both iterates lie on the same right affine branch. Hence,
	$$
	\xi^{(\ell+1)}
	=
	\xi_h^{m+1}.
	$$

	\textit{Case 2.4: The previous iterate lies to the left and the current iterate lies to the right.}
	Suppose that
	\begin{align}
		\label{xi_l}
		\xi^{(\ell-1)}
		<
		\xi_h^{m+1}
		<
		\xi^{(\ell)}.
	\end{align}
	Then
	$$
	F(\xi^{(\ell-1)})
	=
	F'_{-}(\xi_h^{m+1})
	\left(
	\xi^{(\ell-1)}-\xi_h^{m+1}
	\right),
	$$
	and
	$$
	F(\xi^{(\ell)})
	=
	F'_{+}(\xi_h^{m+1})
	\left(
	\xi^{(\ell)}-\xi_h^{m+1}
	\right).
	$$
	Substituting these expressions into
	\eqref{first_order_secant} gives
	\begin{align}
		\xi^{(\ell+1)}-\xi_h^{m+1}
		=
		-
		\frac{
			\left[
			F'_{+}(\xi_h^{m+1})
			-
			F'_{-}(\xi_h^{m+1})
			\right]
			\left(
			\xi^{(\ell)}-\xi_h^{m+1}
			\right)
			\left(
			\xi_h^{m+1}-\xi^{(\ell-1)}
			\right)
		}{
			F'_{+}(\xi_h^{m+1})
			\left(
			\xi^{(\ell)}-\xi_h^{m+1}
			\right)
			+
			F'_{-}(\xi_h^{m+1})
			\left(
			\xi_h^{m+1}-\xi^{(\ell-1)}
			\right)
		}.
		\label{eq:cross_root_left_right}
	\end{align}
	Equivalently,
	\begin{align*}
		\xi_h^{m+1}-\xi^{(\ell+1)}
		&=
		\left(
		\xi_h^{m+1}-\xi^{(\ell-1)}
		\right)
		\frac{
			\left[
			F'_{+}(\xi_h^{m+1})
			-
			F'_{-}(\xi_h^{m+1})
			\right]
			\left(
			\xi^{(\ell)}-\xi_h^{m+1}
			\right)
		}{
			F'_{+}(\xi_h^{m+1})
			\left(
			\xi^{(\ell)}-\xi_h^{m+1}
			\right)
			+
			F'_{-}(\xi_h^{m+1})
			\left(
			\xi_h^{m+1}-\xi^{(\ell-1)}
			\right)
		}.
	\end{align*}
	By \eqref{F_P_F_N} and \eqref{xi_l}, the denominator in the above fraction is positive.
	
	 Moreover,
	\begin{align*}
		&F'_{+}(\xi_h^{m+1})
		\left(
		\xi^{(\ell)}-\xi_h^{m+1}
		\right)
		+
		F'_{-}(\xi_h^{m+1})
		\left(
		\xi_h^{m+1}-\xi^{(\ell-1)}
		\right)-
		\left[
		F'_{+}(\xi_h^{m+1})
		-
		F'_{-}(\xi_h^{m+1})
		\right]
		\left(
		\xi^{(\ell)}-\xi_h^{m+1}
		\right)
		\\
		&=
		F'_{-}(\xi_h^{m+1})
		\left(
		\xi^{(\ell)}-\xi^{(\ell-1)}
		\right)
		>0.
	\end{align*}
	Therefore,
	$$
	0
	\le
	\frac{
		\left[
		F'_{+}(\xi_h^{m+1})
		-
		F'_{-}(\xi_h^{m+1})
		\right]
		\left(
		\xi^{(\ell)}-\xi_h^{m+1}
		\right)
	}{
		F'_{+}(\xi_h^{m+1})
		\left(
		\xi^{(\ell)}-\xi_h^{m+1}
		\right)
		+
		F'_{-}(\xi_h^{m+1})
		\left(
		\xi_h^{m+1}-\xi^{(\ell-1)}
		\right)
	}
	<1.
	$$
	Since $\xi_h^{m+1}-\xi^{(\ell-1)}>0$, we obtain
	\begin{align}
		0
		\le
		\xi_h^{m+1}-\xi^{(\ell+1)}
		<
		\xi_h^{m+1}-\xi^{(\ell-1)}.
		\label{eq:cross_root_left_right_bound}
	\end{align}
	Thus,
	$$
	\xi^{(\ell-1)}
	<
	\xi^{(\ell+1)}
	\le
	\xi_h^{m+1},
	$$
	and hence $\xi^{(\ell+1)}\in U_\delta$.

	If
	$$
	\xi^{(\ell+1)}
	=
	\xi_h^{m+1},
	$$
	the iteration terminates after one update. Otherwise,
	$$
	\xi^{(\ell+1)}
	<
	\xi_h^{m+1}
	<
	\xi^{(\ell)}.
	$$
	Applying the same calculation to the next pair
	$\xi^{(\ell)}$ and $\xi^{(\ell+1)}$ gives
	$$
	\xi^{(\ell+1)}
	<
	\xi^{(\ell+2)}
	\le
	\xi_h^{m+1}.
	$$
	
	If
	$$
	\xi^{(\ell+2)}
	=
	\xi_h^{m+1},
	$$
	the iteration terminates after two updates. Otherwise,
	$$
	\xi^{(\ell+1)}
	<
	\xi^{(\ell+2)}
	<
	\xi_h^{m+1}.
	$$
	The two consecutive iterates $\xi^{(\ell+1)}$ and
	$\xi^{(\ell+2)}$ now lie on the same left affine branch of $F$.
	Therefore,
	$$
	\xi^{(\ell+3)}
	=
	\xi_h^{m+1}.
	$$
	Hence, at most three additional updates are required in this case.

\textit{Case 2.5: The previous iterate lies to the right and the current iterate lies to the left.}
Suppose that
\begin{align}
	\label{xi_l_l-1}
	\xi^{(\ell)}
	<
	\xi_h^{m+1}
	<
	\xi^{(\ell-1)}.
\end{align}
Then
$$
F(\xi^{(\ell)})
=
F'_{-}(\xi_h^{m+1})
\left(
\xi^{(\ell)}-\xi_h^{m+1}
\right),
$$
and
$$
F(\xi^{(\ell-1)})
=
F'_{+}(\xi_h^{m+1})
\left(
\xi^{(\ell-1)}-\xi_h^{m+1}
\right).
$$
Substituting these expressions into \eqref{first_order_secant} gives
\begin{align}
	\xi^{(\ell+1)}-\xi_h^{m+1}
	=
	-
	\frac{
		\left[
		F'_{+}(\xi_h^{m+1})
		-
		F'_{-}(\xi_h^{m+1})
		\right]
		\left(
		\xi^{(\ell-1)}-\xi_h^{m+1}
		\right)
		\left(
		\xi_h^{m+1}-\xi^{(\ell)}
		\right)
	}{
		F'_{+}(\xi_h^{m+1})
		\left(
		\xi^{(\ell-1)}-\xi_h^{m+1}
		\right)
		+
		F'_{-}(\xi_h^{m+1})
		\left(
		\xi_h^{m+1}-\xi^{(\ell)}
		\right)
	}.
	\label{eq:cross_root_right_left}
\end{align}
Equivalently,
\begin{align*}
	\xi_h^{m+1}-\xi^{(\ell+1)}
	&=
	\left(
	\xi_h^{m+1}-\xi^{(\ell)}
	\right)
	\frac{
		\left[
		F'_{+}(\xi_h^{m+1})
		-
		F'_{-}(\xi_h^{m+1})
		\right]
		\left(
		\xi^{(\ell-1)}-\xi_h^{m+1}
		\right)
	}{
		F'_{+}(\xi_h^{m+1})
		\left(
		\xi^{(\ell-1)}-\xi_h^{m+1}
		\right)
		+
		F'_{-}(\xi_h^{m+1})
		\left(
		\xi_h^{m+1}-\xi^{(\ell)}
		\right)
	}.
\end{align*}
By \eqref{F_P_F_N} and \eqref{xi_l_l-1}, the denominator in the above fraction is positive.

Moreover,
\begin{align*}
	&F'_{+}(\xi_h^{m+1})
	\left(
	\xi^{(\ell-1)}-\xi_h^{m+1}
	\right)
	+
	F'_{-}(\xi_h^{m+1})
	\left(
	\xi_h^{m+1}-\xi^{(\ell)}
	\right)
	-
	\left[
	F'_{+}(\xi_h^{m+1})
	-
	F'_{-}(\xi_h^{m+1})
	\right]
	\left(
	\xi^{(\ell-1)}-\xi_h^{m+1}
	\right)
	\\
	&=
	F'_{-}(\xi_h^{m+1})
	\left(
	\xi^{(\ell-1)}-\xi^{(\ell)}
	\right)
	>0.
\end{align*}
Therefore,
$$
0
\le
\frac{
	\left[
	F'_{+}(\xi_h^{m+1})
	-
	F'_{-}(\xi_h^{m+1})
	\right]
	\left(
	\xi^{(\ell-1)}-\xi_h^{m+1}
	\right)
}{
	F'_{+}(\xi_h^{m+1})
	\left(
	\xi^{(\ell-1)}-\xi_h^{m+1}
	\right)
	+
	F'_{-}(\xi_h^{m+1})
	\left(
	\xi_h^{m+1}-\xi^{(\ell)}
	\right)
}
<1.
$$
Since $\xi_h^{m+1}-\xi^{(\ell)}>0$, we obtain
$$
0
\le
\xi_h^{m+1}-\xi^{(\ell+1)}
<
\xi_h^{m+1}-\xi^{(\ell)}.
$$
Therefore,
$$
\xi^{(\ell)}
<
\xi^{(\ell+1)}
\le
\xi_h^{m+1},
$$
and hence $\xi^{(\ell+1)}\in U_\delta$.

If
	$$
	\xi^{(\ell+1)}
	=
	\xi_h^{m+1},
	$$
	the iteration terminates after one update. Otherwise,
	$$
	\xi^{(\ell)}
	<
	\xi^{(\ell+1)}
	<
	\xi_h^{m+1}.
	$$
	Thus, $\xi^{(\ell)}$ and $\xi^{(\ell+1)}$ lie on the same left affine branch, which belongs to Case 2.2 and the next update gives
	$$
	\xi^{(\ell+2)}
	=
	\xi_h^{m+1}.
	$$
	Hence, at most two additional updates are required in this case.
	
	Combining the above cases, once two distinct consecutive iterations belong to $U_\delta$, the secant iteration remains in $U_\delta$ and reaches $\xi_h^{m+1}$ after at most three additional updates.
\end{proof}
\begin{remark}
	In the Lagrange multiplier correction, the mass conservation constraint leads to a scalar nonlinear algebraic equation for the global Lagrange multiplier. The secant method is commonly employed to solve such equations in related works, whereas its local convergence has not been theoretically established. To the best of our knowledge, Theorem \ref{secant_proof} provides the first rigorous proof of the local convergence of the corresponding secant iteration.
\end{remark}

Next, we consider the total mass conservation and nodal nonnegativity of the first-order scheme \eqref{first_order_predictor_n}-\eqref{eq:first_order_correction_b_system}. It follows from \eqref{first_order_mass_conservation}, \eqref{first_order_kkt_n}, and \eqref{first_order_kkt_b} that the resulting first-order fully discrete scheme preserves these two desired structural properties.
\begin{theorem}[Total mass conservation
]\label{first_mass}
	The first-order fully discrete scheme \eqref{first_order_predictor_n}-\eqref{eq:first_order_correction_b_system} is unconditionally total mass conservative, namely,
\begin{align}
	(n_h^m+b_h^m,1)_h=(n_h^0+b_h^0,1)_h,
	\qquad m=1,\ldots,M.
\end{align}
\end{theorem}
%\begin{proof}
%	The conclusion follows directly from \eqref{first_order_mass_conservation}. The proof is complete.
%\end{proof}
\begin{theorem}[Nonnegativity-preserving]
	\label{first_positivity}
	Let $n_h^0$, $b_h^0 \in V_h^r$ be nonnegative,
	Then the numerical solutions of \eqref{first_order_predictor_n}-\eqref{first_order_mass_conservation} satisfy
	$$
	n_h^m(\mathbf{x}_i) \ge 0, \quad b_h^m(\mathbf{x}_i) \ge 0,
	\qquad m=1,\ldots,M, \quad\forall \mathbf{x}_i\in\Sigma_h .
	$$
\end{theorem}
\begin{remark}
	A notable feature of the present scheme is that both $n_h^m$ and $b_h^m$ remain nonnegative for arbitrary time step sizes.
\end{remark}
\subsection{Second-order fully discrete scheme}
In this subsection, we consider the second-order fully discrete scheme which is constructed using the BDF2 formula together with the KKT shift and the Lagrange multiplier method. Moreover, different from the standard Adams-Bashforth extrapolation \cite{HJQ2023,JCB2016} does not necessarily preserve positivity, the extrapolated coefficients are modified following the positivity-preserving extrapolation developed for degenerate diffusion problems \cite{SX20}.

For $m\ge 1$, define
\begin{align}
	n_h^{m+1,*}(\mathbf{x}_i)
	=
	\begin{cases}
		2n_h^m(\mathbf{x}_i)-n_h^{m-1}(\mathbf{x}_i),
		& n_h^m(\mathbf{x}_i)\ge n_h^{m-1}(\mathbf{x}_i),\\[0.6em]
		\dfrac{n_h^m(\mathbf{x}_i) n_h^{m-1}(\mathbf{x}_i)}{2n_h^{m-1}(\mathbf{x}_i)-n_h^m(\mathbf{x}_i)},
		& n_h^m(\mathbf{x}_i)<n_h^{m-1}(\mathbf{x}_i),
	\end{cases}
	\label{positive_extrapolation_n}
\end{align}
and
\begin{align}
	b_h^{m+1,*}(\mathbf{x}_i)
	=
	\begin{cases}
		2b_h^m(\mathbf{x}_i)-b_h^{m-1}(\mathbf{x}_i),
		& b_h^m(\mathbf{x}_i)\ge b_h^{m-1}(\mathbf{x}_i),\\[0.6em]
		\dfrac{b_h^m(\mathbf{x}_i) b_h^{m-1}(\mathbf{x}_i)}{2b_h^{m-1}(\mathbf{x}_i)-b_h^m(\mathbf{x}_i)},
		& b_h^m(\mathbf{x}_i)<b_h^{m-1}(\mathbf{x}_i).
	\end{cases}
	\label{positive_extrapolation_b}
\end{align}

Consequently, if $n_h^m(\mathbf{x}_i)\ge 0$,
$n_h^{m-1}(\mathbf{x}_i)\ge 0$,
$b_h^m(\mathbf{x}_i)\ge 0$, and
$b_h^{m-1}(\mathbf{x}_i)\ge 0$ for all $i$, then
$$
n_h^{m+1,*}(\mathbf{x}_i)\ge 0,
\qquad
b_h^{m+1,*}(\mathbf{x}_i)\ge 0.
$$

With the positivity-preserving extrapolated values defined above, we are now ready to formulate the second-order fully discrete scheme.

\noindent\textbf{Scheme 2.}
Given $\lambda_h^0=\eta_h^0=\xi_h^0=0$, $n_h^0=\Pi_h n_0$, and $b_h^0=\Pi_h b_0$, let $(n_h^1,b_h^1,\lambda_h^1,\eta_h^1,\xi_h^1)$ be generated by \textbf{Scheme 1}. Then, we compute $n_h^{m+1}$ and $b_h^{m+1}$ for $m=1,2,\ldots,M-1$ as follows.

\noindent\textbf{Step 1.}
Find $\tilde n_h^{m+1}\in V_h^r$ such that,
\begin{align}
	\left(
	\frac{3\tilde n_h^{m+1}-4n_h^m+n_h^{m-1}}{2\tau},
	\phi_h
	\right)_h
	+
	(\nabla\tilde n_h^{m+1},\nabla\phi_h)_h
	+
	(n_h^{m+1,*}b_h^{m+1,*}, \phi_h)_h
	=
	(\lambda_h^m+\theta\xi_h^m, \phi_h)_h, \quad \forall \phi_h\in V_h^r.
	\label{second_order_predictor_n}
\end{align}

\noindent\textbf{Step 2.}
Find $\tilde b_h^{m+1}\in V_h^r$ such that,
\begin{align}
	\left(
	\frac{3\tilde b_h^{m+1}-4b_h^m+b_h^{m-1}}{2\tau},
	\psi_h
	\right)_h
	+
	\left(
	\sigma n_h^{m+1,*}b_h^{m+1,*}
	\nabla\tilde b_h^{m+1},
	\nabla\psi_h
	\right)_h
	&=
	(n_h^{m+1,*}b_h^{m+1,*}, \psi_h)_h
	\nonumber\\
	&\quad
	+
	\chi_0
	\left(
	\frac{\sigma n_h^{m+1,*}(b_h^{m+1,*})^2}{(1+n_h^{m+1,*})^2}
	\nabla n_h^{m+1,*},
	\nabla\psi_h
	\right)_h
	\nonumber\\
	&\quad
	+
	(\eta_h^m+(1-\theta)\xi_h^m, \psi_h)_h, \quad \forall\psi_h\in V_h^r.
	\label{second_order_predictor_b}
\end{align}

\noindent\textbf{Step 3.}
Find $\xi_h^{m+1}\in\mathbb{R}$ from the scalar nonlinear algebraic equation:
\begin{align}
	F(\xi_h^{m+1})
	&:=
	\sum_{i=1}^{N}
	\left(
	\tilde n_h^{m+1}(\mathbf{x}_i)
	+
	\frac{2\tau}{3}
	\left(
	\theta(\xi_h^{m+1}-\xi_h^m)
	-
	\lambda_h^m(\mathbf{x}_i)
	\right)
	\right)_+
	\omega_i
	\nonumber\\
	&\quad
	+
	\sum_{i=1}^{N}
	\left(
	\tilde b_h^{m+1}(\mathbf{x}_i)
	+
	\frac{2\tau}{3}
	\left(
	(1-\theta)(\xi_h^{m+1}-\xi_h^m)
	-
	\eta_h^m(\mathbf{x}_i)
	\right)
	\right)_+
	\omega_i	
	-
	(n_h^m+b_h^m, 1)_h
	=0.
	\label{second_order_xi_equation}
\end{align}
Here, $(s)_+=\max\{s,0\}$ and $\theta\in(0,1)$.

\noindent\textbf{Step 4.}
Find $n_h^{m+1}$ and $\lambda_h^{m+1}$ from
\begin{subequations}
	\label{second_order_correction_n_system}
	\begin{align}
		&\frac{3}{2\tau}
		\left(
		n_h^{m+1}(\mathbf{x}_i)
		-
		\tilde n_h^{m+1}(\mathbf{x}_i)
		\right)
		=
		\lambda_h^{m+1}(\mathbf{x}_i)
		-
		\lambda_h^m(\mathbf{x}_i)
		+
		\theta(\xi_h^{m+1}-\xi_h^m),
		\qquad
		\forall\mathbf{x}_i\in\Sigma_h,
		\label{second_order_correction_n_eq}
		\\
		&\lambda_h^{m+1}(\mathbf{x}_i)
		\geq 0,\qquad
		n_h^{m+1}(\mathbf{x}_i)\geq 0,\qquad
		\lambda_h^{m+1}(\mathbf{x}_i)
		n_h^{m+1}(\mathbf{x}_i)=0,
		\qquad
		\forall\mathbf{x}_i\in\Sigma_h.
		\label{second_order_kkt_n}
	\end{align}
\end{subequations}

\noindent\textbf{Step 5.}
Find $b_h^{m+1}$ and $\eta_h^{m+1}$ from
\begin{subequations}
	\label{second_order_correction_b_system}
	\begin{align}
		&\frac{3}{2\tau}
		\left(
		b_h^{m+1}(\mathbf{x}_i)
		-
		\tilde b_h^{m+1}(\mathbf{x}_i)
		\right)
		=
		\eta_h^{m+1}(\mathbf{x}_i)
		-
		\eta_h^m(\mathbf{x}_i)
		+
		(1-\theta)(\xi_h^{m+1}-\xi_h^m),
		\qquad
		\forall\mathbf{x}_i\in\Sigma_h,
		\label{second_order_correction_b_eq}
		\\
		&\eta_h^{m+1}(\mathbf{x}_i)
		\geq 0,\qquad
		b_h^{m+1}(\mathbf{x}_i)\geq 0,\qquad
		\eta_h^{m+1}(\mathbf{x}_i)
		b_h^{m+1}(\mathbf{x}_i)=0,
		\qquad
		\forall\mathbf{x}_i\in\Sigma_h.
		\label{second_order_kkt_b}
	\end{align}
\end{subequations}

Similar to that in \textbf{Scheme 1}, we can see that the predicted solutions $\tilde n_h^{m+1}$ and $\tilde b_h^{m+1}$ in \textbf{Steps 1} and \textbf{2} of \textbf{Scheme 2} are obtained from standard BDF2 time-stepping schemes and may not preserve nonnegativity or the total mass. Therefore, we introduce a global scalar mass-correction multiplier $\xi_h^{m+1}$ in \textbf{Step 3} to enforce the conservation of the total discrete mass, together with the nodal Lagrange multipliers
$\lambda_h^{m+1}$ and $\eta_h^{m+1}$ to impose the nonnegativity constraints on $n_h^{m+1}$ and $b_h^{m+1}$ in \textbf{Steps 4} and \textbf{5}, respectively.

In fact, after $\xi_h^{m+1}$ is computed in \textbf{Step 3}, \textbf{Steps 4} and \textbf{5} admit explicit pointwise solutions.
For any $\mathbf{x}_i\in\Sigma_h$, $n_h^{m+1}$ and $\lambda_h^{m+1}$ satisfy
\begin{align}
	\left(
	n_h^{m+1}(\mathbf{x}_i),
	\lambda_h^{m+1}(\mathbf{x}_i)
	\right)
	&=
	\begin{cases}
		\left(
		\displaystyle
		\tilde n_h^{m+1}(\mathbf{x}_i)
		+
		\frac{2\tau}{3}
		\left(
		\theta(\xi_h^{m+1}-\xi_h^m)
		-
		\lambda_h^m(\mathbf{x}_i)
		\right),
		0
		\right),
		&
		\begin{aligned}
			&\tilde n_h^{m+1}(\mathbf{x}_i)
			+
			\frac{2\tau}{3}
			\left(
			\theta(\xi_h^{m+1}-\xi_h^m)
			-
			\lambda_h^m(\mathbf{x}_i)
			\right)
			>0,
		\end{aligned}
		\\
		\left(
		0,
		\displaystyle
		-\frac{3}{2\tau}
		\left(
		\tilde n_h^{m+1}(\mathbf{x}_i)
		+
		\frac{2\tau}{3}
		\left(
		\theta(\xi_h^{m+1}-\xi_h^m)
		-
		\lambda_h^m(\mathbf{x}_i)
		\right)
		\right)
		\right),
		&
		\begin{aligned}
			&\tilde n_h^{m+1}(\mathbf{x}_i)
			+
			\frac{2\tau}{3}
			\left(
			\theta(\xi_h^{m+1}-\xi_h^m)
			-
			\lambda_h^m(\mathbf{x}_i)
			\right)
			\leq 0.
		\end{aligned}
	\end{cases}
	\label{second_order_pointwise_n}
\end{align}

Similarly, for any $\mathbf{x}_i\in\Sigma_h$, $b_h^{m+1}$ and $\eta_h^{m+1}$ satisfy
\begin{align}
	\left(
	b_h^{m+1}(\mathbf{x}_i),
	\eta_h^{m+1}(\mathbf{x}_i)
	\right)
	&=
	\begin{cases}
		\left(
		\displaystyle
		\tilde b_h^{m+1}(\mathbf{x}_i)
		+
		\frac{2\tau}{3}
		\left(
		(1-\theta)(\xi_h^{m+1}-\xi_h^m)
		-
		\eta_h^m(\mathbf{x}_i)
		\right),
		0
		\right),
		&
		\begin{aligned}
			&\tilde b_h^{m+1}(\mathbf{x}_i)
			+
			\frac{2\tau}{3}
			\left(
			(1-\theta)(\xi_h^{m+1}-\xi_h^m)
			-
			\eta_h^m(\mathbf{x}_i)
			\right)
			>0,
		\end{aligned}
		\\
		\left(
		0,
		\displaystyle
		-\frac{3}{2\tau}
		\left(
		\tilde b_h^{m+1}(\mathbf{x}_i)
		+
		\frac{2\tau}{3}
		\left(
		(1-\theta)(\xi_h^{m+1}-\xi_h^m)
		-
		\eta_h^m(\mathbf{x}_i)
		\right)
		\right)
		\right),
		&
		\begin{aligned}
			&\tilde b_h^{m+1}(\mathbf{x}_i)
			+
			\frac{2\tau}{3}
			\left(
			(1-\theta)(\xi_h^{m+1}-\xi_h^m)
			-
			\eta_h^m(\mathbf{x}_i)
			\right)
			\leq 0.
		\end{aligned}
	\end{cases}
	\label{second_order_pointwise_b}
\end{align}

Equations \eqref{second_order_pointwise_n} and
\eqref{second_order_pointwise_b} imply that
\begin{align}
	n_h^{m+1}(\mathbf{x}_i)
	&=
	\left(
	\tilde n_h^{m+1}(\mathbf{x}_i)
	+
	\frac{2\tau}{3}
	\left(
	\theta(\xi_h^{m+1}-\xi_h^m)
	-
	\lambda_h^m(\mathbf{x}_i)
	\right)
	\right)_+,
	\label{second_order_positive_part_n}
	\\
	b_h^{m+1}(\mathbf{x}_i)
	&=
	\left(
	\tilde b_h^{m+1}(\mathbf{x}_i)
	+
	\frac{2\tau}{3}
	\left(
	(1-\theta)(\xi_h^{m+1}-\xi_h^m)
	-
	\eta_h^m(\mathbf{x}_i)
	\right)
	\right)_+.
	\label{second_order_positive_part_b}
\end{align}
Therefore, the nodal nonnegativity conditions
$n_h^{m+1}(\mathbf{x}_i)\geq 0$ and
$b_h^{m+1}(\mathbf{x}_i)\geq 0$ for all
$\mathbf{x}_i\in\Sigma_h$ are satisfied automatically.

Substituting \eqref{second_order_positive_part_n} and
\eqref{second_order_positive_part_b} into
\eqref{second_order_xi_equation}, we obtain
\begin{align}
	(n_h^{m+1}+b_h^{m+1},1)_h
	-
	(n_h^m+b_h^m,1)_h
	=0,
	\label{second_order_mass_conservation}
\end{align}
which confirms the conservation of the total discrete mass.

Moreover, for \textbf{Scheme 2}, it is clear that only linear systems need to be solved in \textbf{Steps 1} and \textbf{2}, while \textbf{Steps 4} and \textbf{5} consist only of nodal assignment operations. Therefore, the remaining issue is to establish the existence and uniqueness of the solution to the scalar nonlinear
algebraic equation \eqref{second_order_xi_equation} in
\textbf{Step 3}. This result is stated in the following theorem.
\begin{theorem}[Unique solvability of \eqref{second_order_xi_equation}]
	Let $\theta\in (0, 1)$ and $(n^{0}_{h}+b^{0}_{h},1)_{h}>0$. Then the nonlinear algebraic equation
	\eqref{second_order_xi_equation} admits a unique solution $\xi^{m+1}_{h}\in\mathbb{R}$. Moreover, the corresponding mass-correction multiplier satisfies
	$$
	\xi^{m+1}_{h}\leq0.
	$$
\end{theorem}

\begin{proof}
	The proof of existence and uniqueness is analogous to that of the Theorem \ref{first_order_solvability}. Indeed, $F$ is continuous on $\mathbb{R}$, since it is a finite sum of positive-part functions. Moreover, since $\theta\in (0, 1)$,
	$$
	\lim_{\xi\to-\infty}F(\xi)=-(n^{m}_{h}+b^{m}_{h},1)_{h}<0,
	\qquad
	\lim_{\xi\to+\infty}F(\xi)=+\infty.
	$$
	Hence, the intermediate value theorem gives the existence of at least one zero point.
	
	Furthermore, each term $	\biggl[
	\tilde n_h^{m+1}(\mathbf x_i)
	+\frac{2\tau}{3}
	\Bigl(
	\theta(\xi_h^{m+1}-\xi_h^m)
	-\lambda_h^m(\mathbf x_i)
	\Bigr)
	\biggr]_+$ and $\biggl[
	\tilde b_h^{m+1}(\mathbf x_i)
	+\frac{2\tau}{3}
	\Bigl(
	(1-\theta)(\xi_h^{m+1}-\xi_h^m)
	-\eta_h^m(\mathbf x_i)
	\Bigr)
	\biggr]_+$ is nondecreasing with respect to $\xi_h^{m+1}$. Thus $F$ is nondecreasing. If $F$ had two distinct zero points $\xi_1<\xi_2$, then, using the same argument as in the first-order case, at least one positive-part term would be strictly positive at $\xi_1$ and hence strictly larger at $\xi_2$, while all other terms are nondecreasing. Therefore,
	$$
	F(\xi_2)>F(\xi_1)=0,
	$$
	which contradicts $F(\xi_2)=0$. Hence the zero point is unique.
	
	It remains to prove the sign of the corresponding multiplier. Adding \eqref{second_order_predictor_n} and \eqref{second_order_correction_n_eq}, for any $\phi_h\in V^{r}_{h}$,
	\begin{align}
		\left(
		\frac{3n^{m+1}_{h}-4n^{m}_{h}+n^{m-1}_{h}}{2\tau},
		\phi_h
		\right)_{h}
		+(\nabla\tilde n^{m+1}_{h},\nabla\phi_h)_{h}
		+(b^{m+1,*}_{h}n^{m+1,*}_{h},\phi_h)_{h}
		=
		(\lambda^{m+1}_{h}+\theta\xi^{m+1}_{h},\phi_h)_{h}.
	\end{align}
	Similarly, adding \eqref{second_order_predictor_b} and \eqref{second_order_correction_b_eq}, for any $\psi_h\in V^{r}_{h}$,
	\begin{align}
		\left(
		\frac{3b^{m+1}_{h}-4b^{m}_{h}+b^{m-1}_{h}}{2\tau},
		\psi_h
		\right)_{h}
		&+
		(\sigma n^{m+1,*}_{h}b^{m+1,*}_{h}\nabla\tilde b^{m+1}_{h},\nabla\psi_h)_{h}
		\notag\\
		&=
		(n^{m+1,*}_{h}b^{m+1,*}_{h},\psi_h)_{h}
		+\chi_0	\left(
		\frac{\sigma n_h^{m+1,*}(b_h^{m+1,*})^2}{(1+n_h^{m+1,*})^2}
		\nabla n_h^{m+1,*},
		\nabla\psi_h
		\right)_h
		\notag\\
		&\quad
		+(\eta^{m+1}_{h}+(1-\theta)\xi^{m+1}_{h},\psi_h)_{h}.
	\end{align}
	Taking $\phi_h=\psi_h=1$, we obtain
	\begin{align}
		\left(
		\frac{3n^{m+1}_{h}-4n^{m}_{h}+n^{m-1}_{h}}{2\tau}
		+
		\frac{3b^{m+1}_{h}-4b^{m}_{h}+b^{m-1}_{h}}{2\tau},
		1
		\right)_{h}
		=
		(\lambda^{m+1}_{h}+\eta^{m+1}_{h}+\xi^{m+1}_{h},1)_{h}.
	\end{align}
	By the mass conservation constraint of the correction step and the induction hypothesis,
	\[
	(n^{m+1}_{h}+b^{m+1}_{h},1)_{h}
	=
	(n^{m}_{h}+b^{m}_{h},1)_{h}
	=
	(n^{m-1}_{h}+b^{m-1}_{h},1)_{h}.
	\]
	Therefore,
	\[
	\left(
	\frac{3n^{m+1}_{h}-4n^{m}_{h}+n^{m-1}_{h}}{2\tau}
	+
	\frac{3b^{m+1}_{h}-4b^{m}_{h}+b^{m-1}_{h}}{2\tau},
	1
	\right)_{h}
	=0.
	\]
	It follows that
	\[
	(\lambda^{m+1}_{h}+\eta^{m+1}_{h}+\xi^{m+1}_{h},1)_{h}=0.
	\]
	Since $\xi^{m+1}_{h}$ is a scalar multiplier, we have
	\[
	\xi^{m+1}_{h}
	=
	-\frac{(\lambda^{m+1}_{h}+\eta^{m+1}_{h},1)_{h}}{(1,1)_{h}}.
	\]
	Finally, by the complementarity conditions,
	\[
	\lambda^{m+1}_{h}(\mathbf{x}_i)\geq0,
	\qquad
	\eta^{m+1}_{h}(\mathbf{x}_i)\geq0,
	\qquad i=1,\ldots,N.
	\]
	Hence
	\[
	(\lambda^{m+1}_{h}+\eta^{m+1}_{h},1)_{h}\geq0,
	\]
	which implies
	\[
	\xi^{m+1}_{h}\leq0.
	\]
	This completes the proof.
\end{proof}
Similar to Theorem \ref{secant_proof}, we can prove the following theorem.
\begin{theorem}[Local convergence of \eqref{first_order_secant}]
	Let $\xi_h^{m+1}$ denote the solution of \eqref{second_order_xi_equation}.
	The secant iteration \eqref{first_order_secant} enjoys the same local convergence property as in the first-order case. In particular, once two distinct consecutive iterations belongs to a sufficiently small neighborhood of $\xi_h^{m+1}$, the secant iteration   reaches $\xi_h^{m+1}$ after at most three additional updates.
\end{theorem}
We next verify that the resulting second-order fully discrete scheme preserves the desired structural properties. Together with \eqref{second_order_mass_conservation}, the KKT conditions \eqref{second_order_kkt_n} and \eqref{second_order_kkt_b} ensure that the proposed second-order fully discrete scheme satisfies total mass conservation and preserves nodal nonnegativity.

\begin{theorem}[Total mass conservation]
	\label{second_mass}
	The second-order fully discrete scheme
	\eqref{second_order_predictor_n}-\eqref{second_order_correction_b_system} unconditionally preserves the total mass, i.e.,
	\begin{align}
	(n_h^m + b_h^m,\, 1)_h
		= (n_h^0 + b_h^0,\, 1)_h,
		\qquad m=1,\ldots, M.
	\end{align}
\end{theorem}
\begin{theorem}[Nonnegativity-preserving]
	\label{second_positivity}
	Let $n_h^0, b_h^0 \in V_h^r$ be nonnegative. Then the numerical solutions of the second-order fully discrete scheme \eqref{second_order_predictor_n}-\eqref{second_order_correction_b_system} satisfy
	$$
   n_h^m(\mathbf{x}_i) \ge 0, \quad b_h^m(\mathbf{x}_i) \ge 0,
   \qquad m=1,\ldots,M, \quad\forall \mathbf{x}_i\in\Sigma_h .
   $$
\end{theorem}
\begin{remark}
Compared with the structure-preserving scheme for a simplified model in \cite{BR2026}, which consists of a first-order scheme in time together with the $\mathbb{P}^1$ element in space, the construction in this paper provides a more general framework to achieve high-order simulations. The same idea can be directly extended to higher-order time discretization schemes in the BDF$k$ type (see, e.g., \cite{CS2022I}) and higher-order  finite element methods with mass-lumping (see, e.g., \cite{GPJ2001}).

\end{remark}

\section{Numerical experiments}
\label{4_test}
In this section, we first verify the convergence order for the proposed schemes. Then we investigate the structure-preserving properties of the schemes. Finally, we will study the growth of Bacillus subtilis colonies on thin agar plates under varying nutrient levels, agar stiffnesses, and chemotactic sensitivities.
\subsection{Convergence order}
\begin{figure}[!htbp]
	\centering
	\subfigure[First-order scheme]{\includegraphics[width=0.45\textwidth]{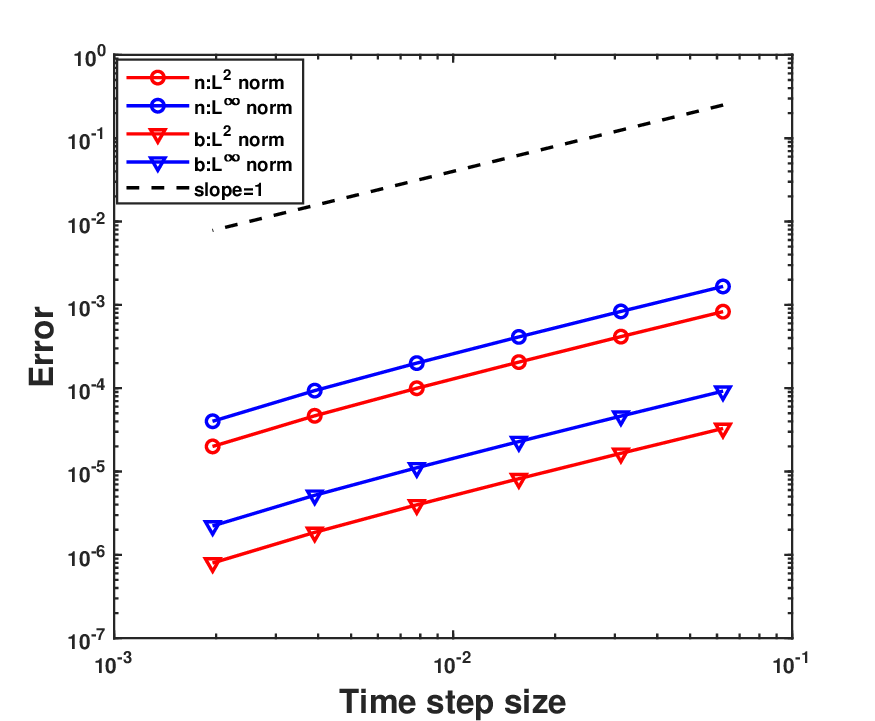}}
	\subfigure[Second-order scheme]{\includegraphics[width=0.45\textwidth]{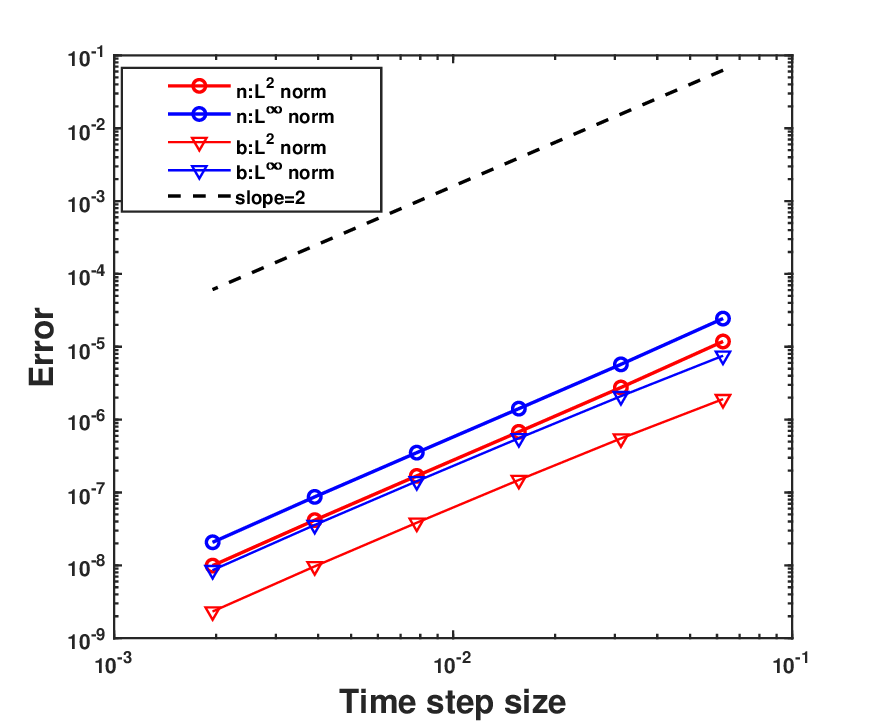}}
	\caption{Temporal convergence orders of the first- and second-order schemes}
	\label{time_convergence}
\end{figure}

\begin{figure}[!htbp]
	\centering
	\subfigure[$\widetilde{\mathbb{P}}^1$ element]{\includegraphics[width=0.45\textwidth]{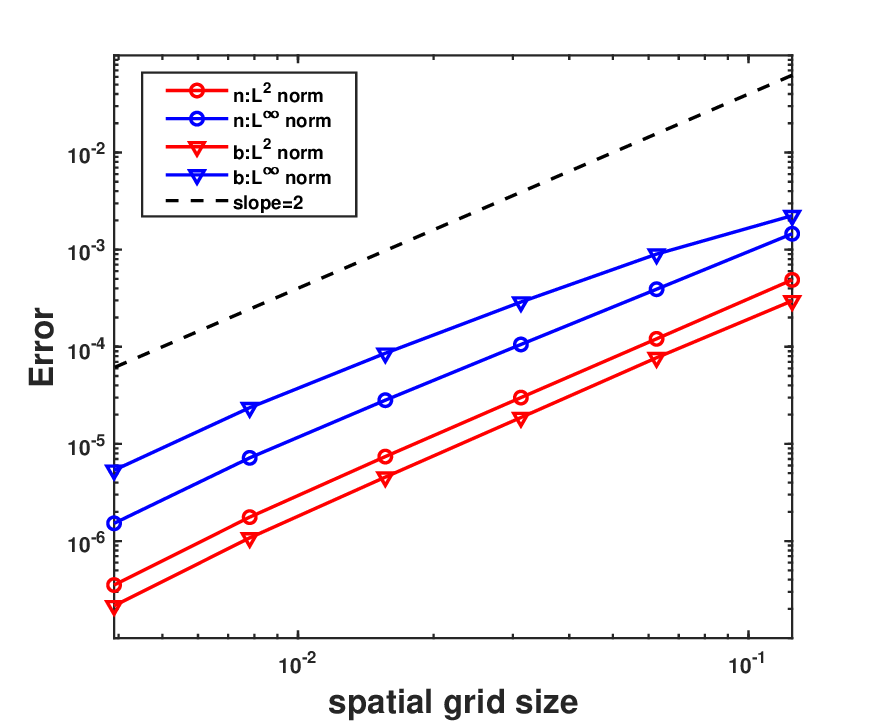}}
	\subfigure[$\widetilde{\mathbb{P}}^2$ element]{\includegraphics[width=0.45\textwidth]{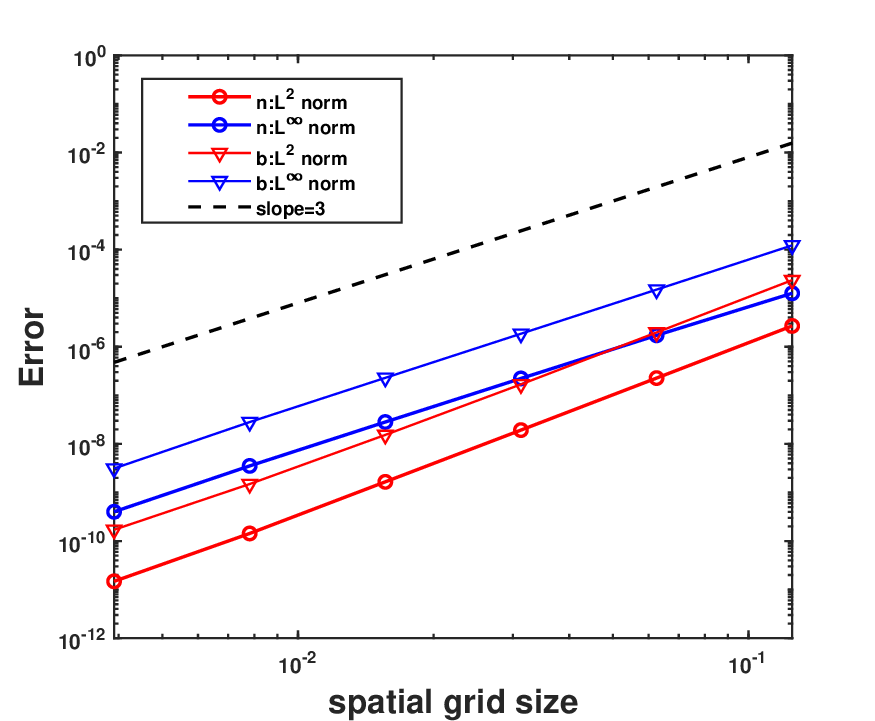}}
	\caption{Spatial convergence orders of the $\widetilde{\mathbb{P}}^1$ and $\widetilde{\mathbb{P}}^2$ elements.}
	\label{space_convergence}
\end{figure}

The computational domain is $\Omega=[0,1]^2$, and the initial data are given by $n(x,y,0)=b(x,y,0)=0.1\cos(\pi x)\cos(\pi y)+0.2$. The final time is set to $T=0.1$, and the weighting parameter in the mass correction is chosen as $\theta=0.5$. First, fixing the spatial mesh size at $h=\frac{1}{128}$ and setting $\tau=\frac{T}{2^m}$, $m=4,5,\ldots,9$, we examine the temporal convergence of the schemes. Since no exact solution is available, the numerical solution computed with $\tau_e=\frac{T}{2048}$ is used as the reference solution. Figure \ref{time_convergence} displays the errors in $L^2$ and $L^\infty$ norms with respect to the time step. It is observed that the first- and second-order temporal convergence rates are obtained in Figures \ref{time_convergence}(a) and (b), respectively. These results are in good agreement with the expected temporal accuracy in both schemes. Next, fixing the time step as $\tau=\frac{T}{1024}$ and setting $h=\frac{1}{2^m}$, $m=3, 4, \ldots, 8$, we investigate the spatial convergence order. Similarly, numerical solutions computed on finer meshes are used as reference solutions, i.e., the reference solutions are obtained with $h_e=\frac{1}{512}$ for the $\widetilde{\mathbb{P}}^1$ and $\widetilde{\mathbb{P}}^2$ elements, respectively. Figure \ref{space_convergence}   presents the errors in $L^2$ and $L^\infty$ norms. The observed rates indicate that the approximation with the $\widetilde{\mathbb{P}}^1$ element yields second-order convergence, while the approximation with $\widetilde{\mathbb{P}}^2$ element achieves third-order convergence. These results demonstrate the validity of the proposed scheme.

	\subsection{Structure-preserving properties}
\begin{figure}[!htbp]
	\centering
	\subfigure[Total mass]{\includegraphics[width=0.32\textwidth]{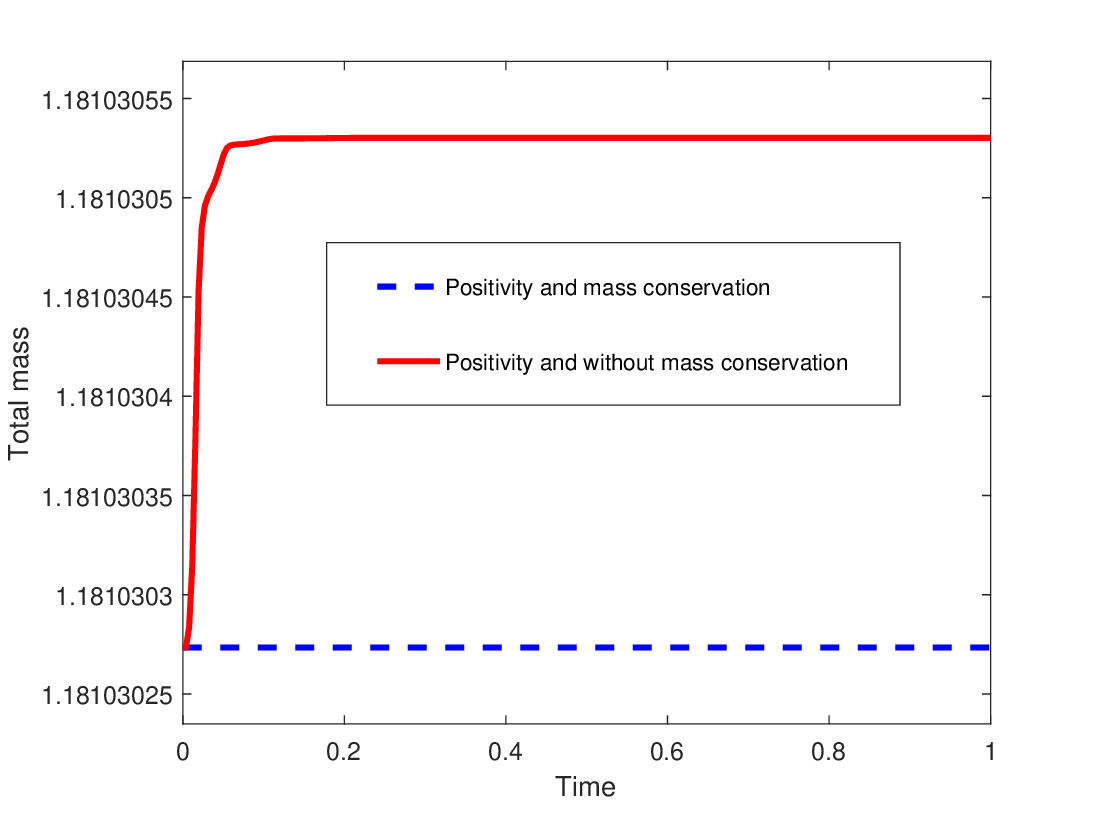}}
	\subfigure[min ($n_{h}^{m}$)]{\includegraphics[width=0.32\textwidth]{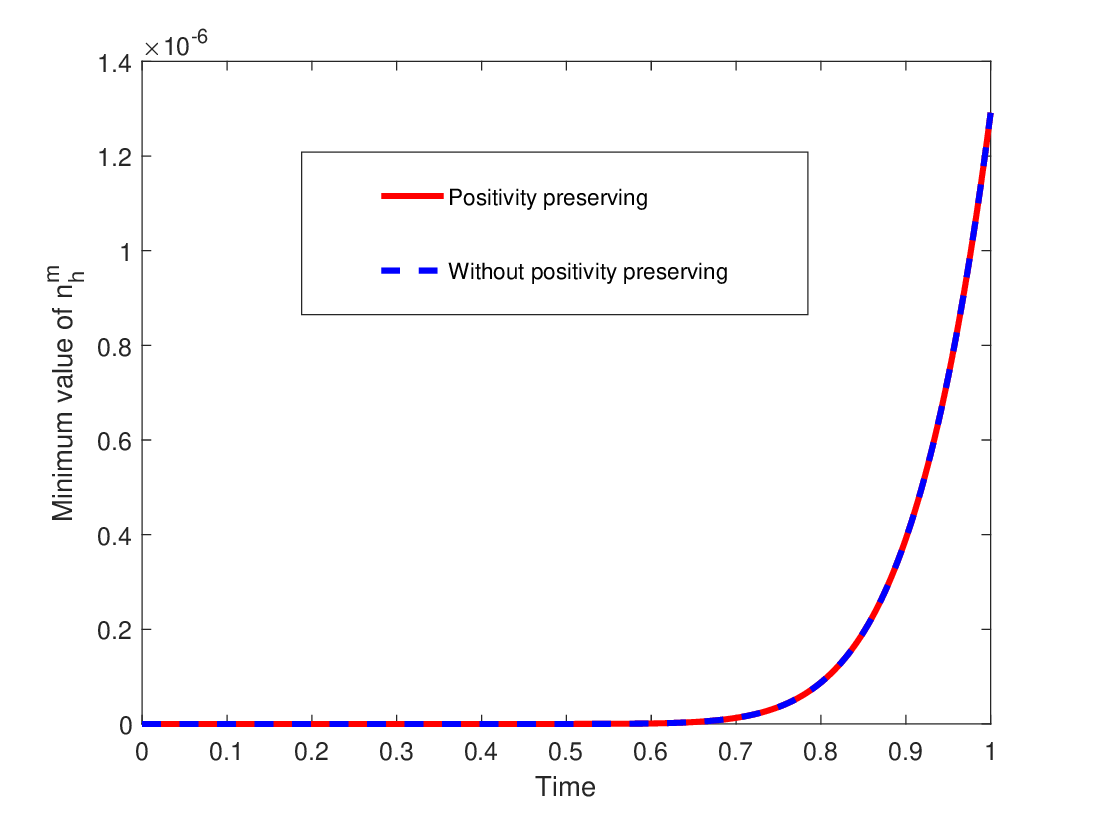}}
	\subfigure[min ($b_{h}^{m}$)]{\includegraphics[width=0.32\textwidth]{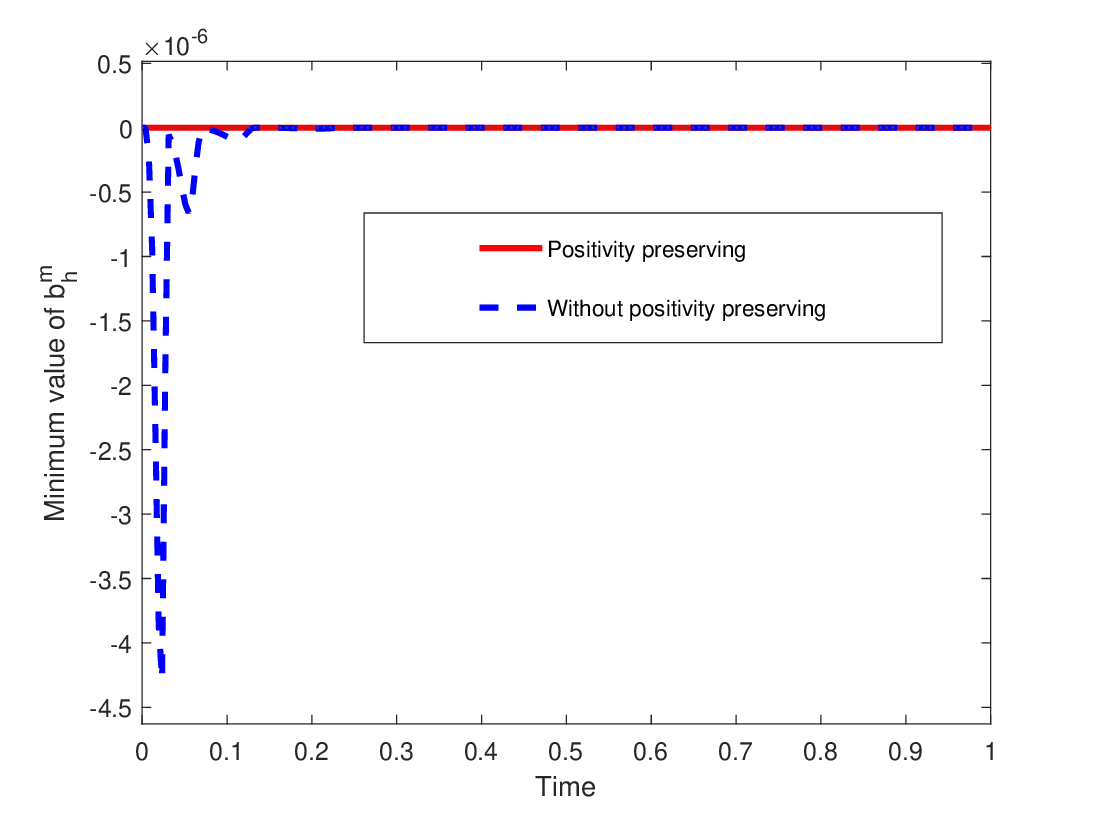}}
	\caption{Nonnegativity and total mass conservation ($\widetilde{\mathbb{P}}^1$ element).}
	\label{preserving}
\end{figure}

\begin{figure}[!htbp]
	\centering
	\subfigure[$\lambda_h^{m}$]{\includegraphics[width=0.40\textwidth]{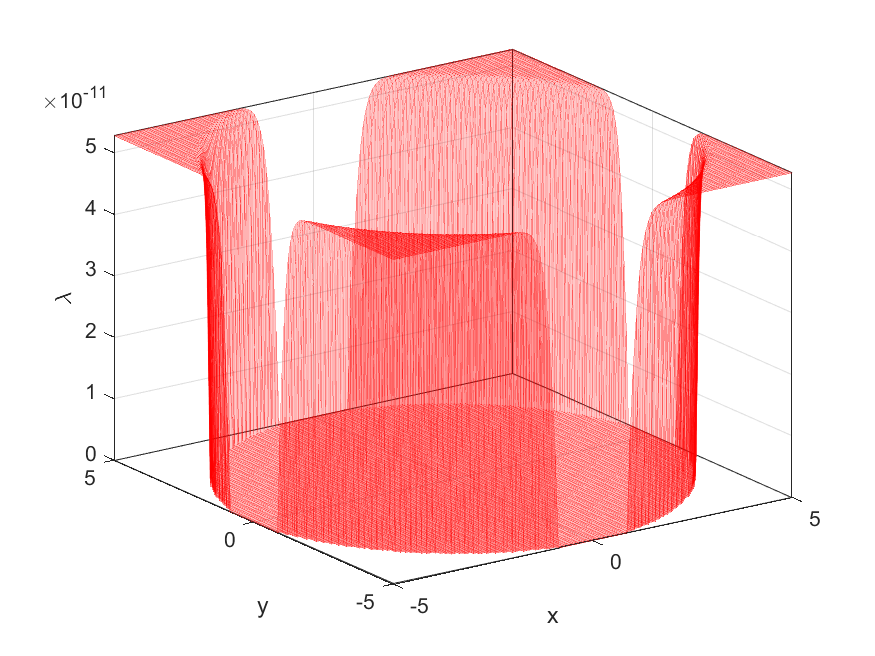}}
	\subfigure[$\eta_h^{m}$]{\includegraphics[width=0.40\textwidth]{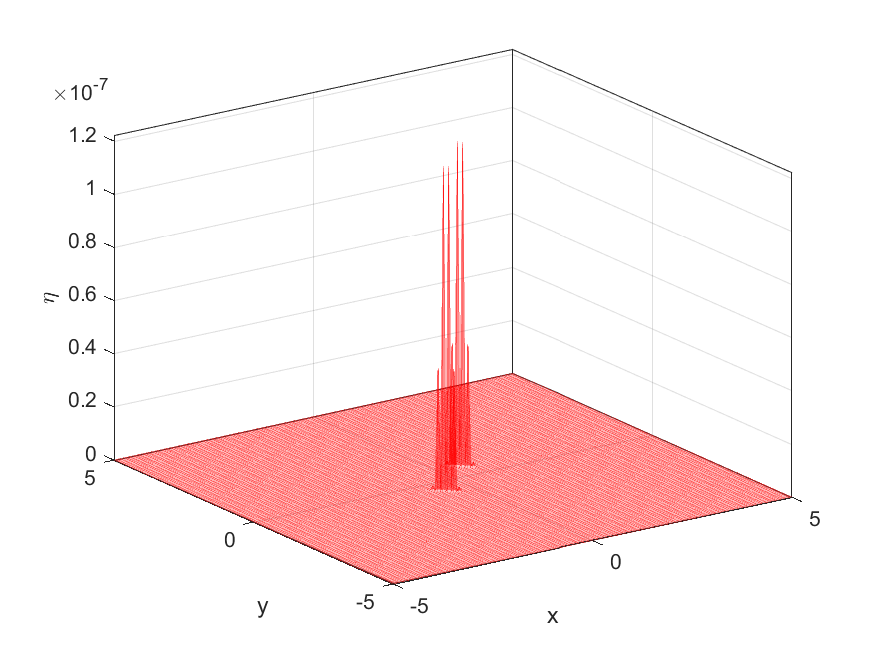}}
	\caption{Lagrange multipliers $\lambda_h^{m}$ and $\eta_h^{m}$ ($t=0.2$) ($\widetilde{\mathbb{P}}^1$ element).}
	\label{Lagrange_multiplier}
\end{figure}

\begin{figure}[!htbp]
		\centering
	\subfigure[$\xi_{h}^{m}$]{\includegraphics[width=0.40\textwidth]{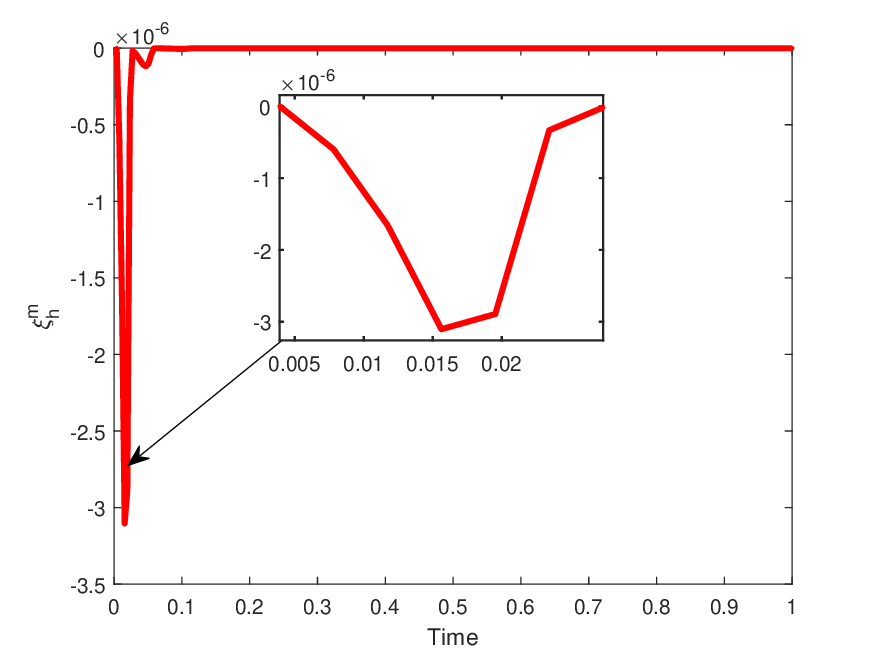}}
	\subfigure[Iteration numbers]{\includegraphics[width=0.40\textwidth]{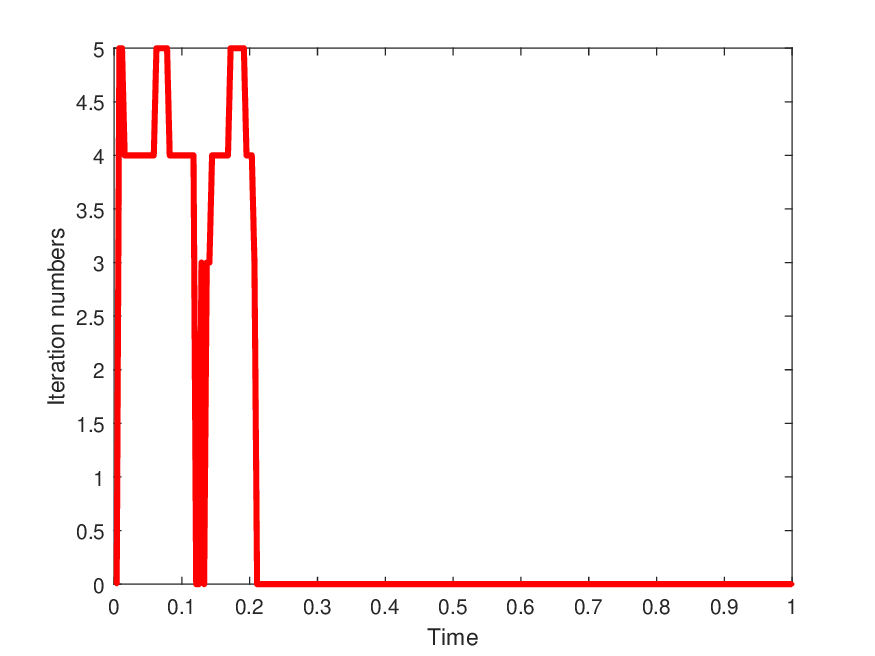}}
	\caption{Mass conservation multiplier $\xi_h^{m}$ and iteration numbers for the nonlinear equation ($\widetilde{\mathbb{P}}^1$ element).}
	\label{xi_simulation_p1}
\end{figure}

\begin{figure}[!htbp]
	\centering
	\subfigure[Total mass]{\includegraphics[width=0.32\textwidth]{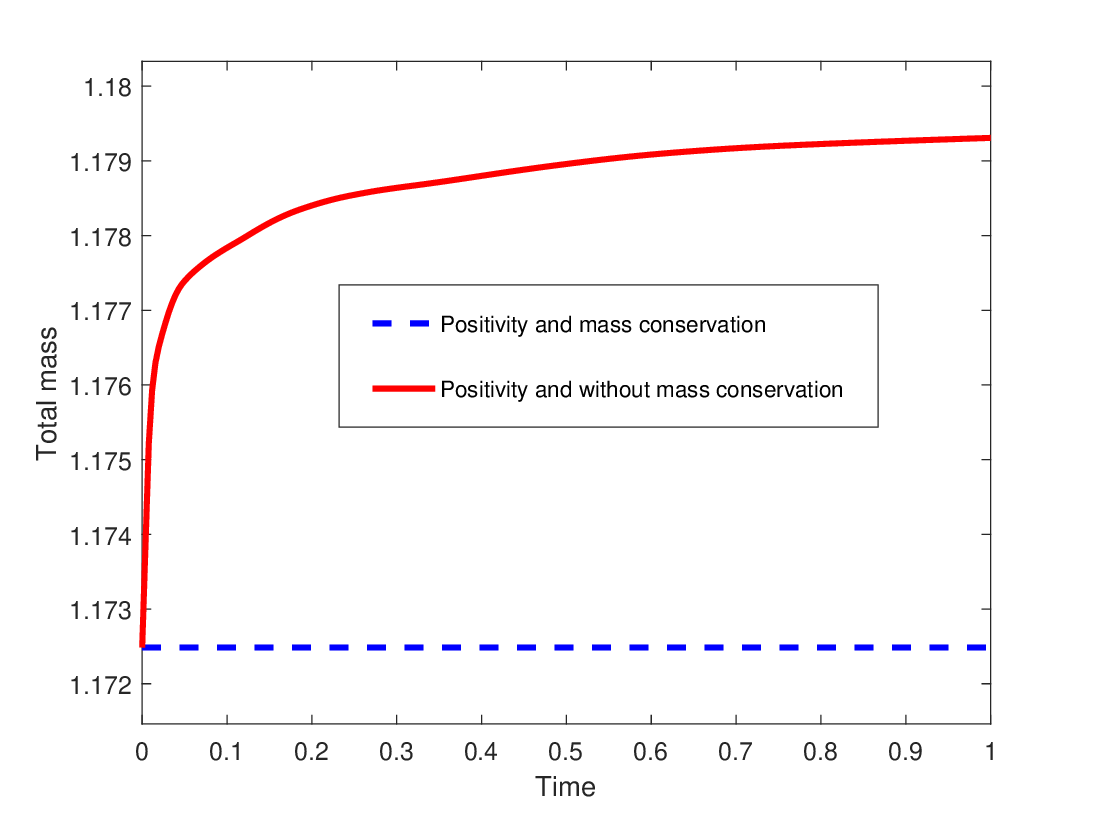}}
	\subfigure[min ($n_{h}^{m}$)]{\includegraphics[width=0.32\textwidth]{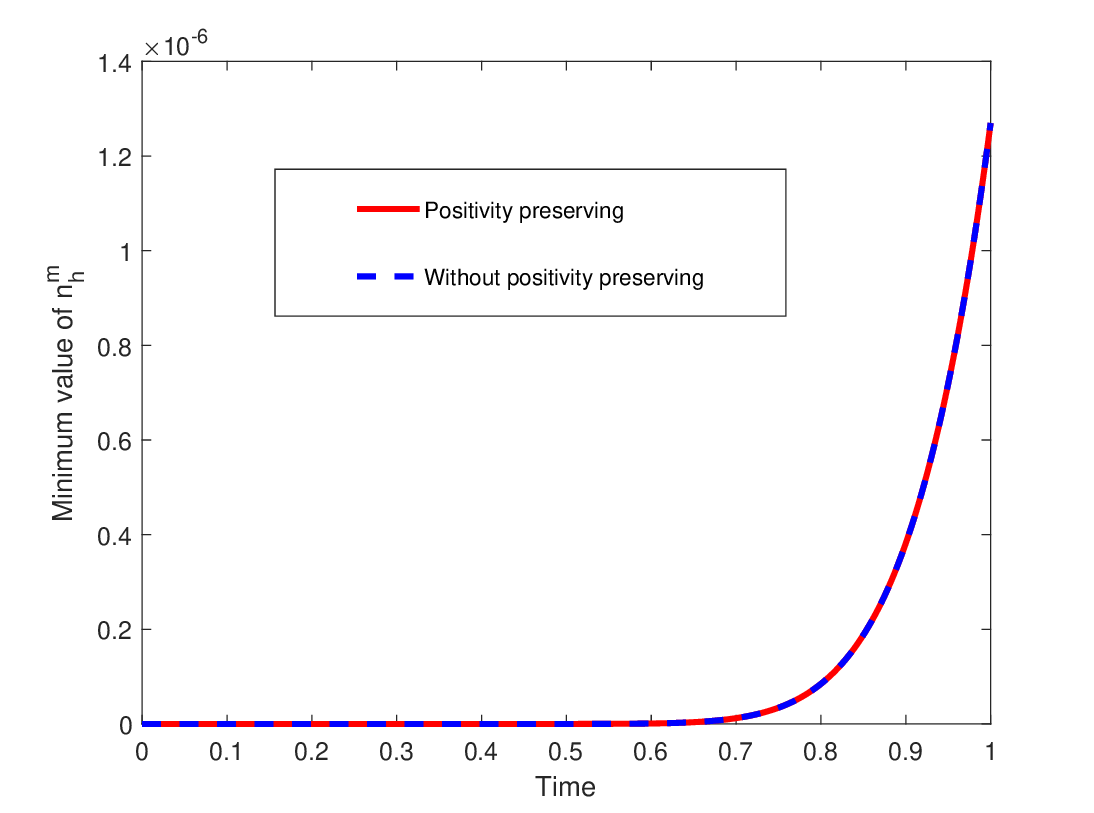}}
	\subfigure[min ($b_{h}^{m}$)]{\includegraphics[width=0.32\textwidth]{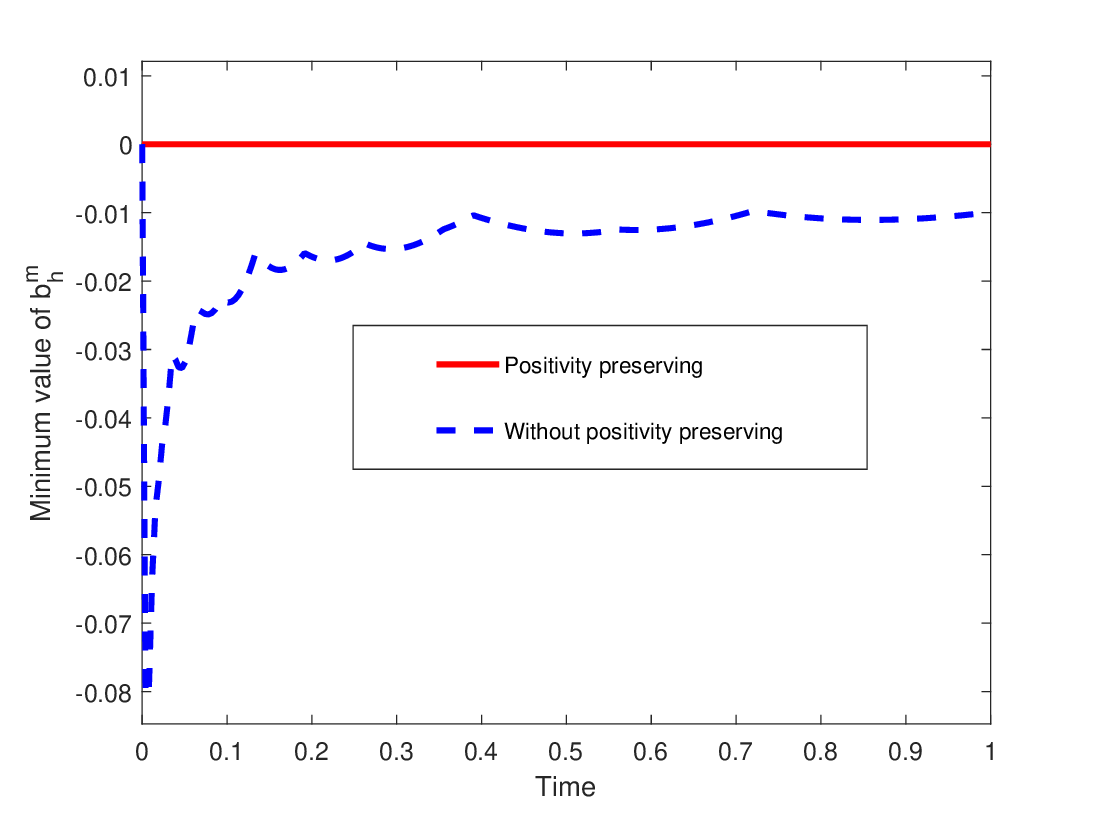}}
	\caption{Nonnegativity and total mass conservation ($\widetilde{\mathbb{P}}^2$ element).}
	\label{P2preserving}
\end{figure}

\begin{figure}[!htbp]
	\centering
	\subfigure[$\lambda_h^{m}$]{\includegraphics[width=0.40\textwidth]{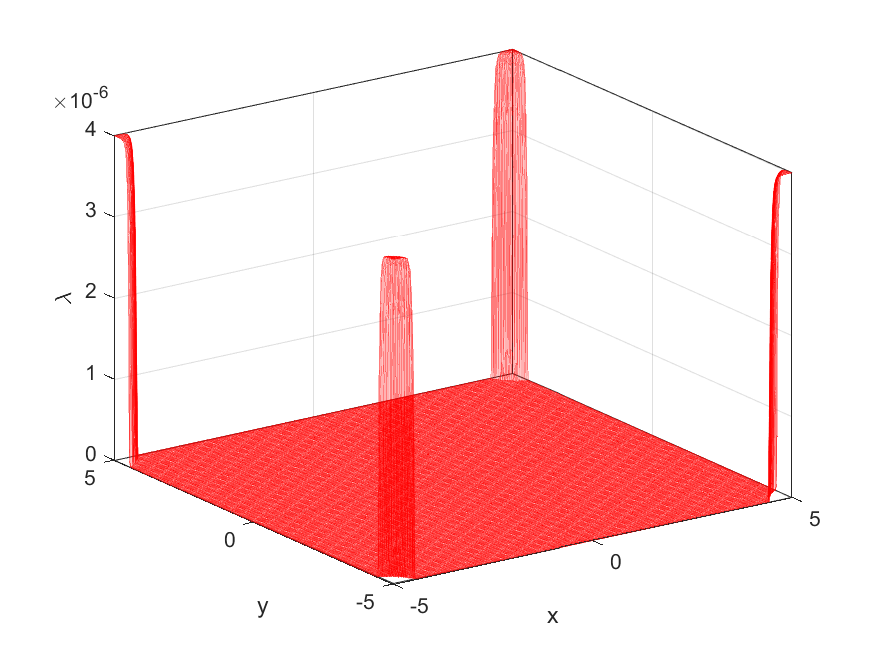}}
	\subfigure[$\eta_h^{m}$]{\includegraphics[width=0.40\textwidth]{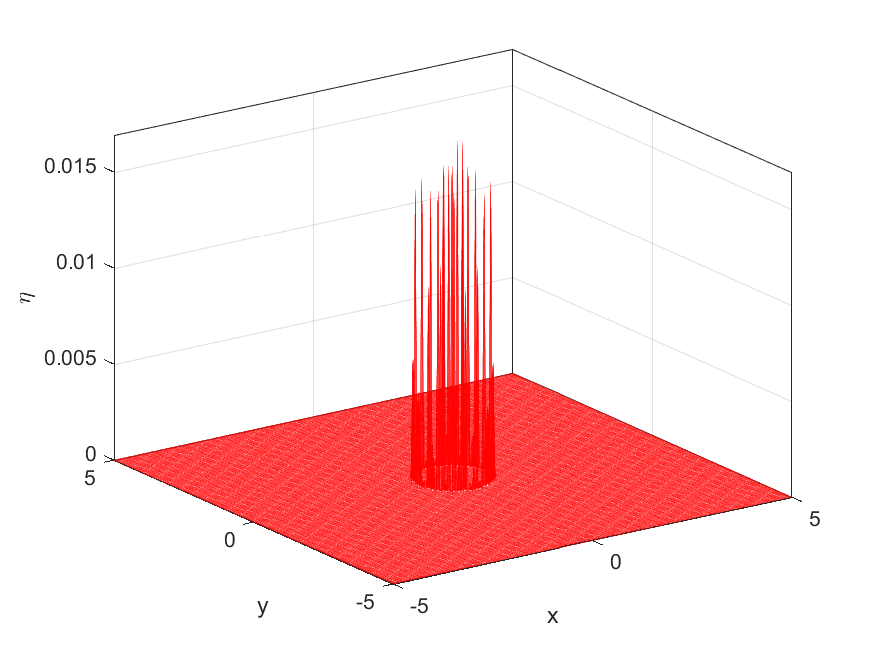}}
	\caption{Lagrange multipliers $\lambda_h^{m}$ and $\eta_h^{m}$ ($t=0.89$) ($\widetilde{\mathbb{P}}^2$ element).}
	\label{Lagrange_multiplier_p2}
\end{figure}

\begin{figure}[!htbp]
	\centering
	\subfigure[$\xi_h^{m}$]{\includegraphics[width=0.4\textwidth]{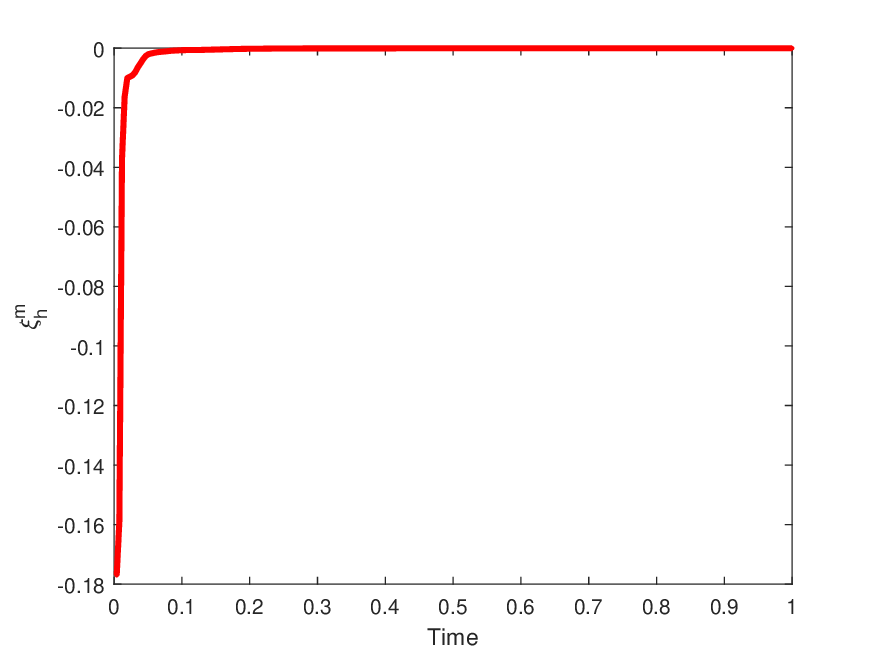}}
	\subfigure[Iteration numbers]{\includegraphics[width=0.4\textwidth]{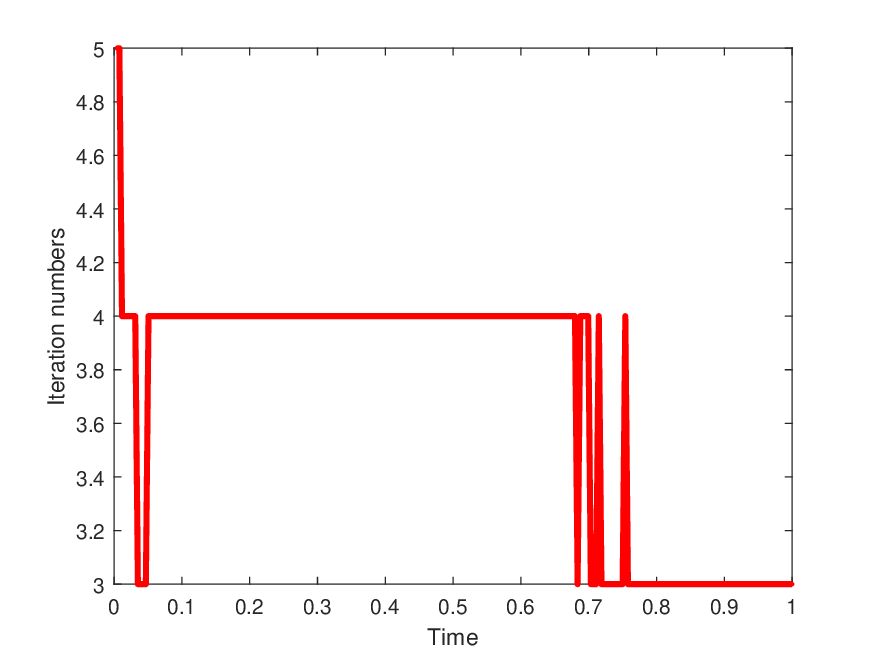}}
	\caption{Mass conservation multiplier $\xi_h^{m}$ and iteration numbers for the nonlinear equation ($\widetilde{\mathbb{P}}^2$ element).}
	\label{xi_simulation_p2}
\end{figure}

To illustrate the structural properties derived above, we perform numerical simulations on $\Omega=[-5,5]^2$. The initial data are chosen as
\begin{align}
	n(x,y,0) &=
	\begin{cases}
		1, & x^2+y^2\le 0.25,\\
		0, & \text{otherwise},
	\end{cases}
	&
	b(x,y,0) &=
	\begin{cases}
		0.5, & x^2+y^2\le 0.25,\\
		0, & \text{otherwise}.
	\end{cases}
\end{align}
Setting the final time $T=1$, $h=\frac{10}{128}$ and $\tau=\frac{1}{256}$, $\theta=0.5$. We collect the corresponding numerical results got by the $\widetilde{\mathbb{P}}^1$ element and the second-order scheme  in Figures \ref{preserving}-\ref{xi_simulation_p1}. Figure \ref{preserving} reports the evolution of the discrete total mass and the minimum nodal values of $n_h^{m}$ and $b_h^{m}$. The numerical results confirm that the minimum nodal values remain nonnegative at all time levels and that the discrete total mass is preserved. Whereas a noticeable drift occurs in the absence of the mass correction and positivity preserving. Figures \ref{Lagrange_multiplier}(a) and (b) display the spatial distributions of the Lagrange multipliers $\lambda_{h}^{m}$ and $\eta_{h}^{m}$ at $t=0.2$. Both multipliers exhibit localized peak structures and very close to $0$, indicating that the positivity constraints are activated mainly in local regions. Figure \ref{xi_simulation_p1}(a) shows the time evolution of the multiplier $\xi_{h}^{m}$ ($\xi_{h}^{m}\leq 0$), which adjusts rapidly at the initial stage and then approaches a steady state. This behavior demonstrates that $\xi_{h}^{m}$ can  continuously enforce the total mass conservation constraint. Figure \ref{xi_simulation_p1}(b) reports the number of secant iterations required to solve the nonlinear equation associated with $\xi_{h}^{m}$, showing that the iteration converges within only a few steps. These results are consistent with the above theoretical analysis and confirm that the scheme preserves nonnegativity and total mass. Using the same parameters as above, the corresponding results for the $\widetilde{\mathbb{P}}^2$ element approximation are presented in Figures \ref{P2preserving}-\ref{xi_simulation_p2}. Similar results are got. These results verify the structure-preserving of the proposed scheme.
\subsection{Colony pattern formation}
In this subsection, we study the growth of the bacterial colony. As in \cite{KAMT1997}, the initial bacterial density is prescribed by
$$
b(x,y,0)=b_M \exp\!\left(-\frac{x^2+y^2}{6.25}\right),
$$
where $b_M=0.71$ denotes the maximum density at the center of the domain. All computations are performed on the two-dimensional square domain $[-170,170]^2$ with $h=0.5$. The time step is chosen as $\tau=0.2$ for $\sigma_0=1$ and $\tau=0.05$ for $\sigma_0=4$. The initial nutrient concentration is taken to be spatially uniform, i.e., $n_0$ is a constant. The parameter $\chi_0$ in the chemotactic term measures the strength of chemotaxis. In the simulations, several values of $\chi_0$ are considered, including $\chi_0=0$, which corresponds to the absence of chemotaxis. The agar hardness is characterized by the parameter $\sigma_0$. The complete set of free parameters used in the numerical experiments is summarized in Table \ref{tab:parameters}.
\begin{table}[htbp]
	\centering
	\caption{Free parameters employed in the numerical simulations.}
	\begin{tabular}{l l l}
		\midrule
		Description & Symbol & Values \\
		\midrule
		Initial nutrient concentration & \(n_0\) & 0.71, 1.07 \\
	    Agar softness & \(\sigma_0\) & 1.0, 4.0 \\
		Chemotactic sensitivity strength & \(\chi_0\) & 0, 5, 7.5\\
		\midrule
	\end{tabular}
	\label{tab:parameters}
\end{table}

\begin{figure}[!htbp]
	\centering
	\subfigure[$t=100$]{\includegraphics[width=0.245\textwidth,trim=18 8 18 8,clip]{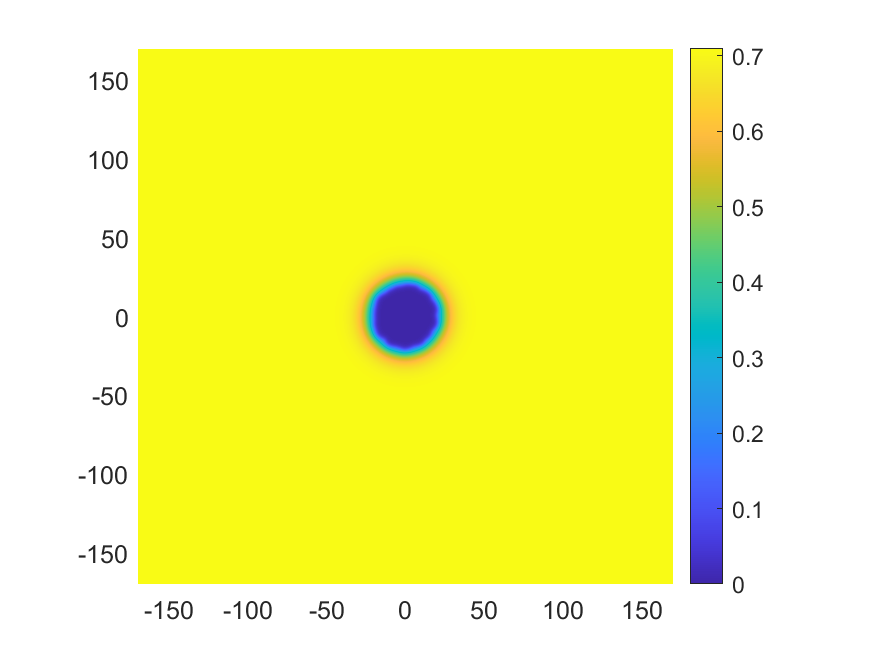}}%
	\subfigure[$t=300$]{\includegraphics[width=0.245\textwidth,trim=18 8 18 8,clip]{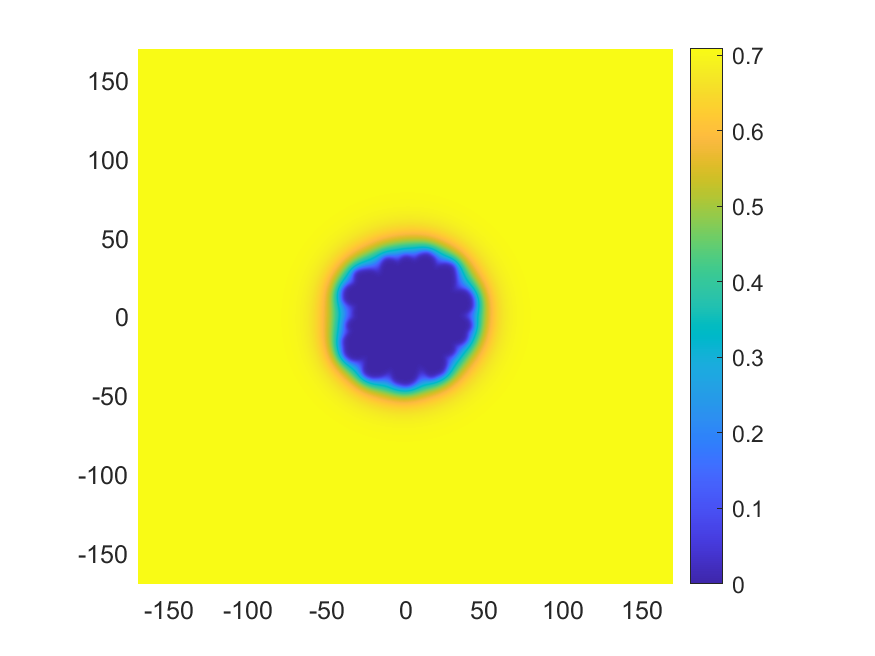}}%
	\subfigure[$t=500$]{\includegraphics[width=0.245\textwidth,trim=18 8 18 8,clip]{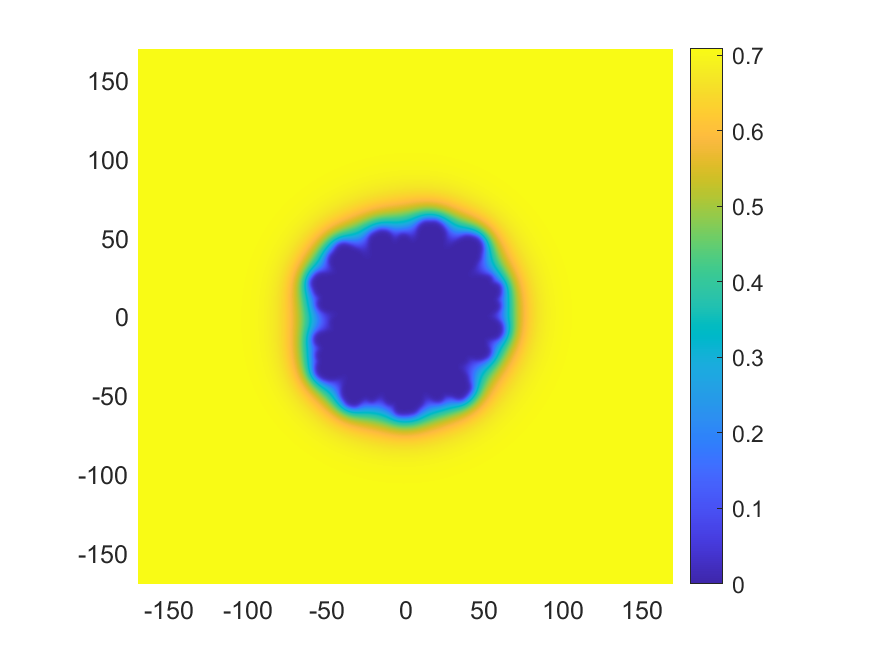}}%
	\subfigure[$t=1000$]{\includegraphics[width=0.245\textwidth,trim=18 8 18 8,clip]{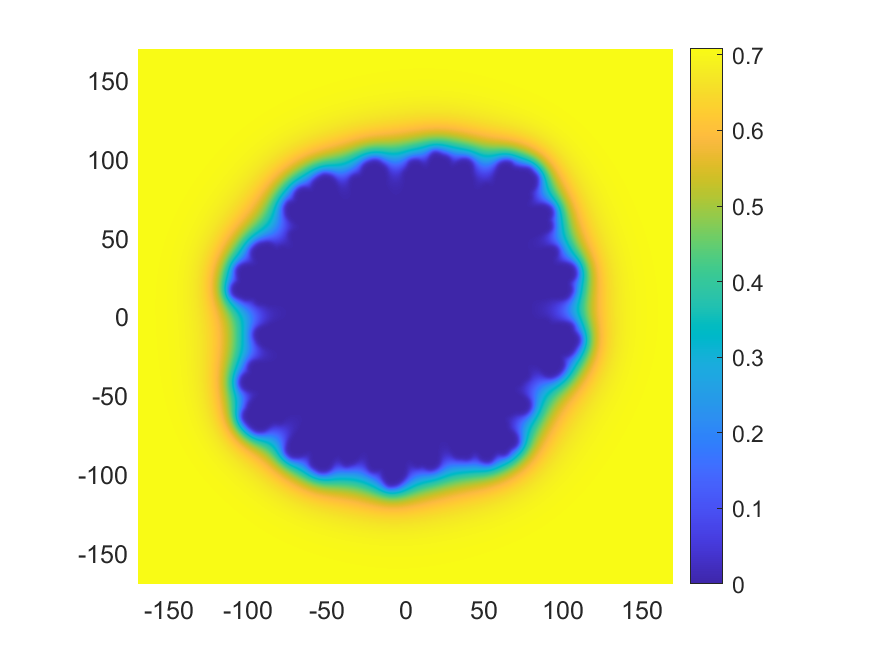}}\\[-1mm]
	\subfigure[$t=100$]{\includegraphics[width=0.245\textwidth,trim=18 8 18 8,clip]{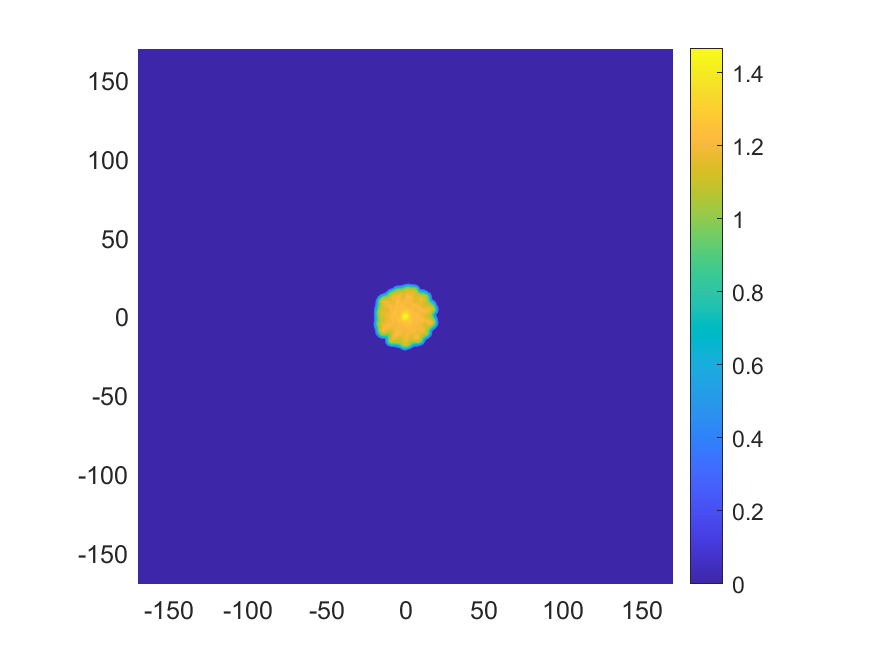}}%
	\subfigure[$t=300$]{\includegraphics[width=0.245\textwidth,trim=18 8 18 8,clip]{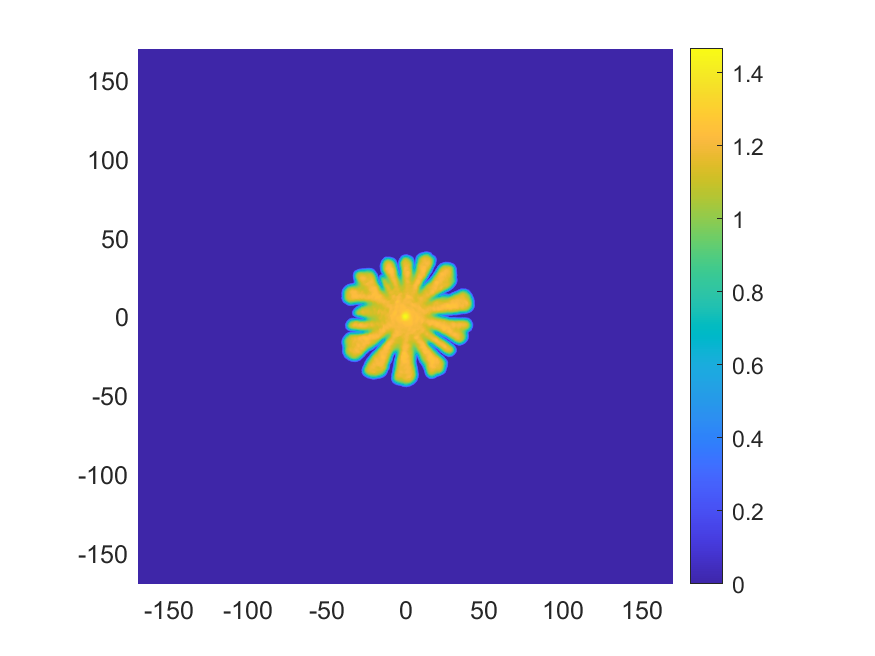}}%
	\subfigure[$t=500$]{\includegraphics[width=0.245\textwidth,trim=18 8 18 8,clip]{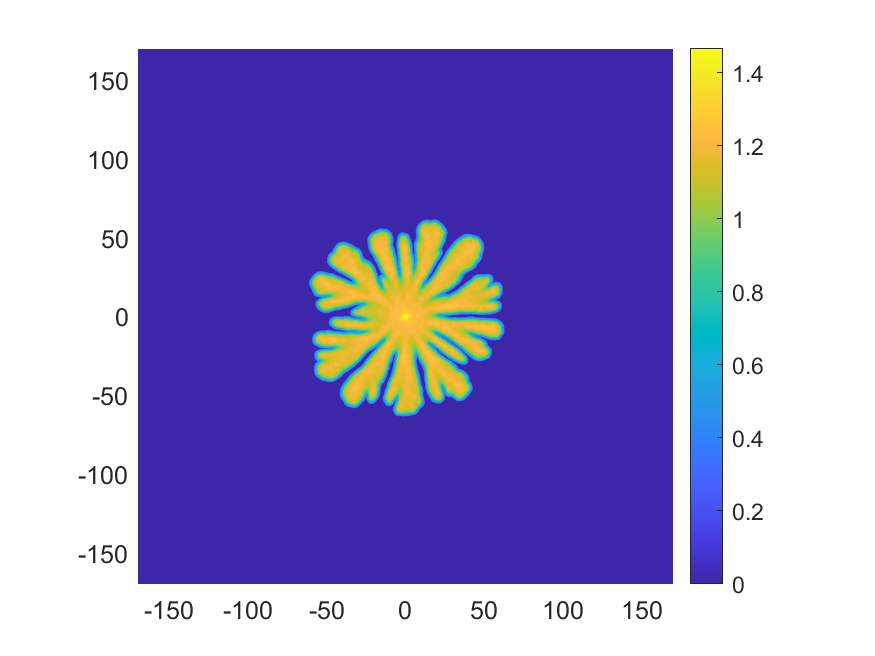}}%
	\subfigure[$t=1000$]{\includegraphics[width=0.245\textwidth,trim=18 8 18 8,clip]{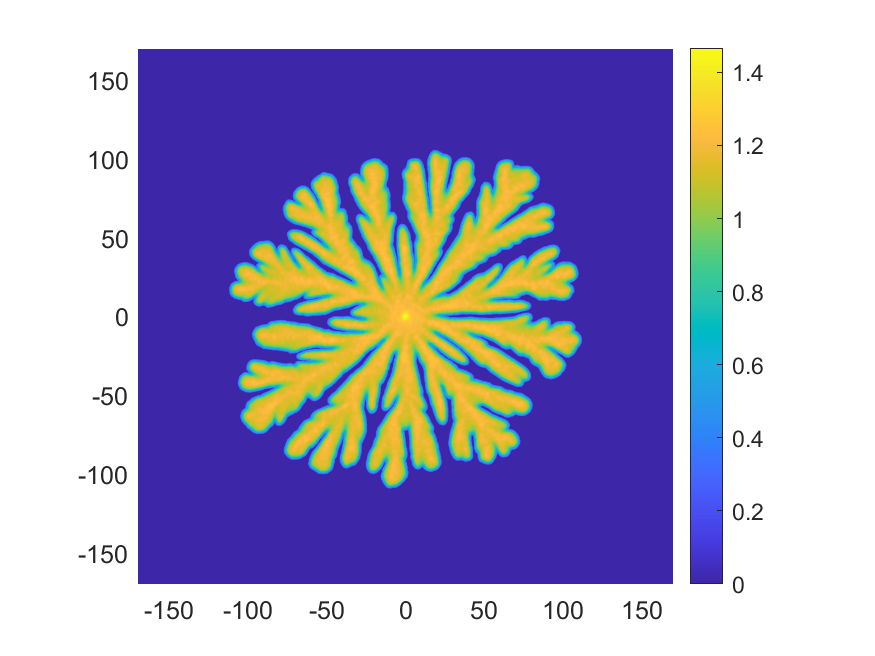}}
	\caption{The evolution of $n_h^m$ (top) and $b_h^m$ (bottom), $n_0=0.71$,  $\sigma_0=1$,  $\chi_0=0$.}
	\label{n_071_sigma_1_chi_0}
\end{figure}

\begin{figure}[!htbp]
	\centering
	\subfigure[$t=100$]{\includegraphics[width=0.245\textwidth,trim=18 8 18 8,clip]{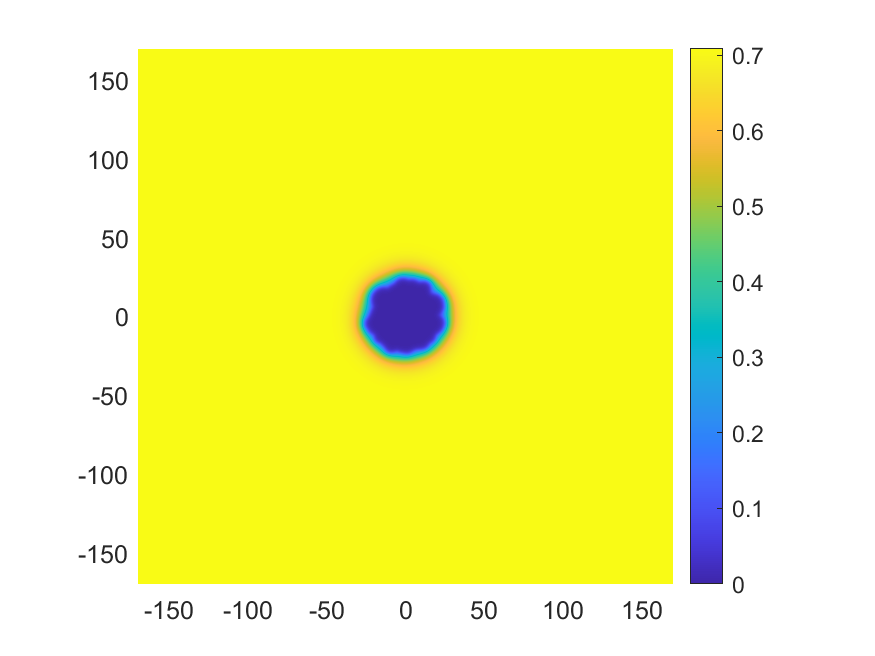}}%
	\subfigure[$t=300$]{\includegraphics[width=0.245\textwidth,trim=18 8 18 8,clip]{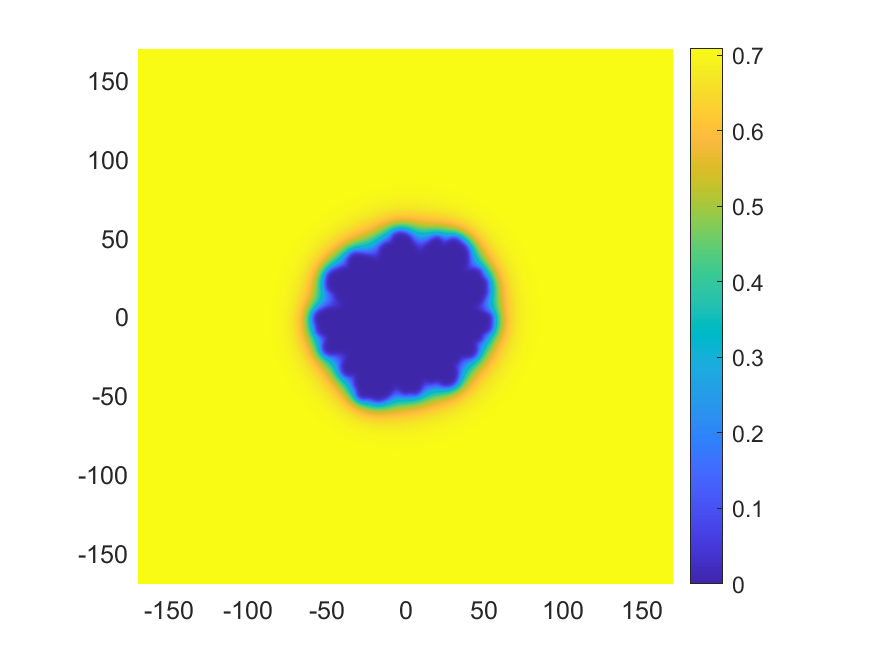}}%
	\subfigure[$t=500$]{\includegraphics[width=0.245\textwidth,trim=18 8 18 8,clip]{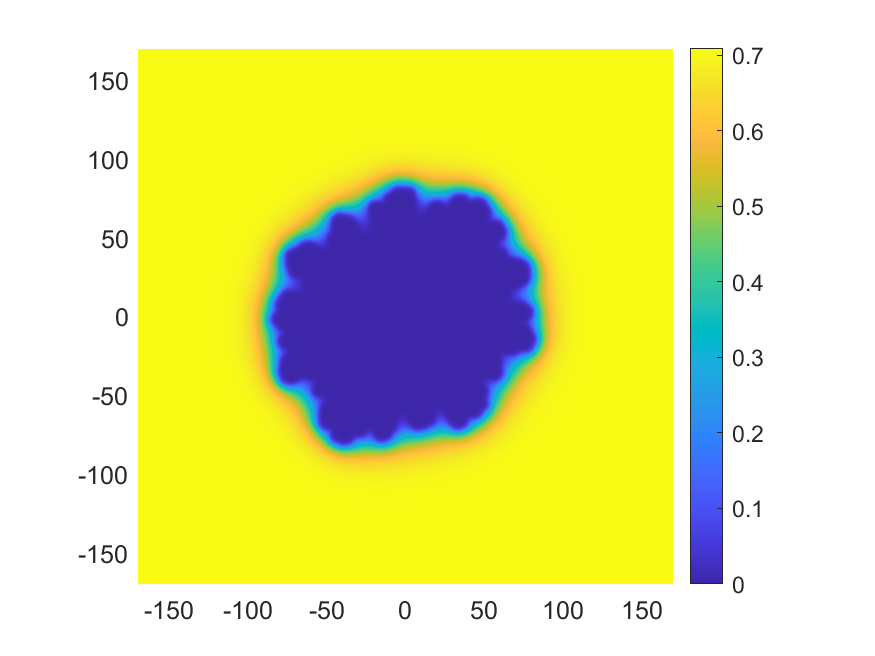}}%
	\subfigure[$t=1000$]{\includegraphics[width=0.245\textwidth,trim=18 8 18 8,clip]{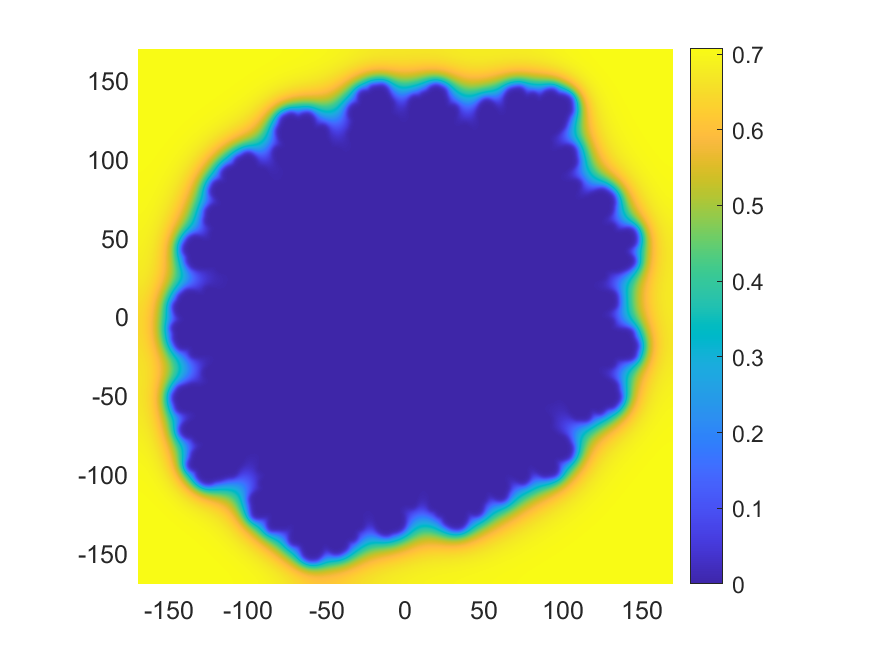}}\\[-1mm]
	\subfigure[$t=100$]{\includegraphics[width=0.245\textwidth,trim=18 8 18 8,clip]{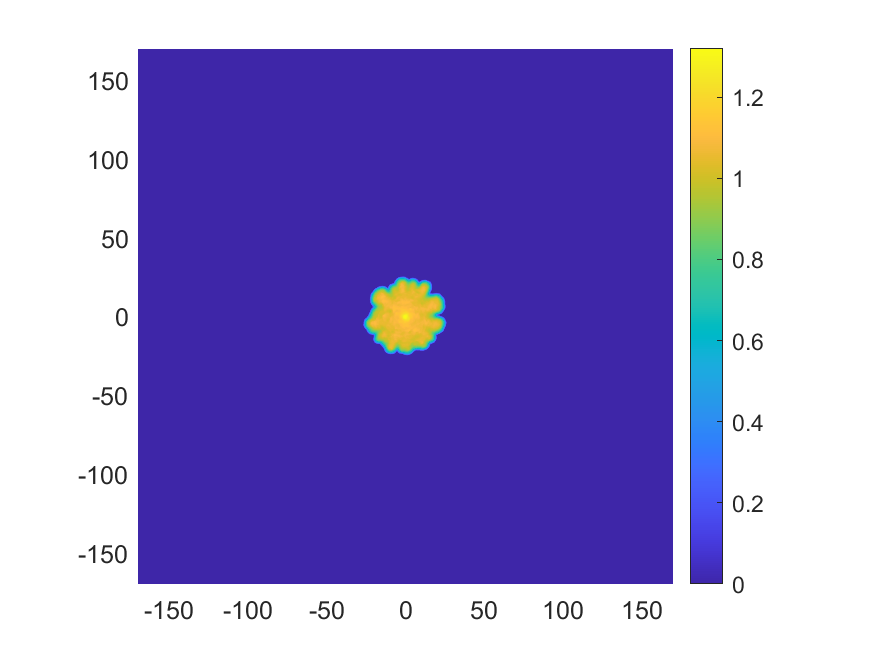}}%
	\subfigure[$t=300$]{\includegraphics[width=0.245\textwidth,trim=18 8 18 8,clip]{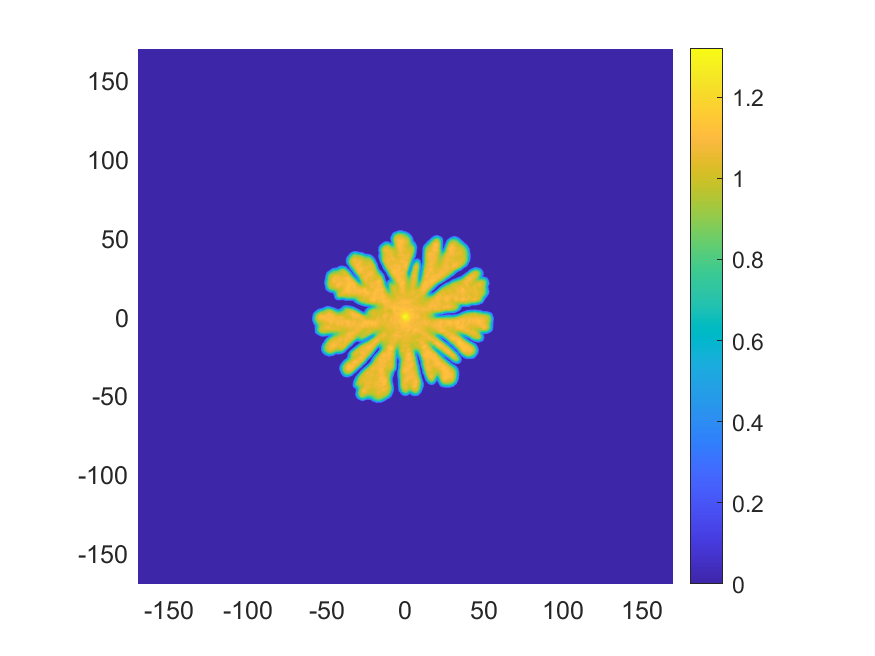}}%
	\subfigure[$t=500$]{\includegraphics[width=0.245\textwidth,trim=18 8 18 8,clip]{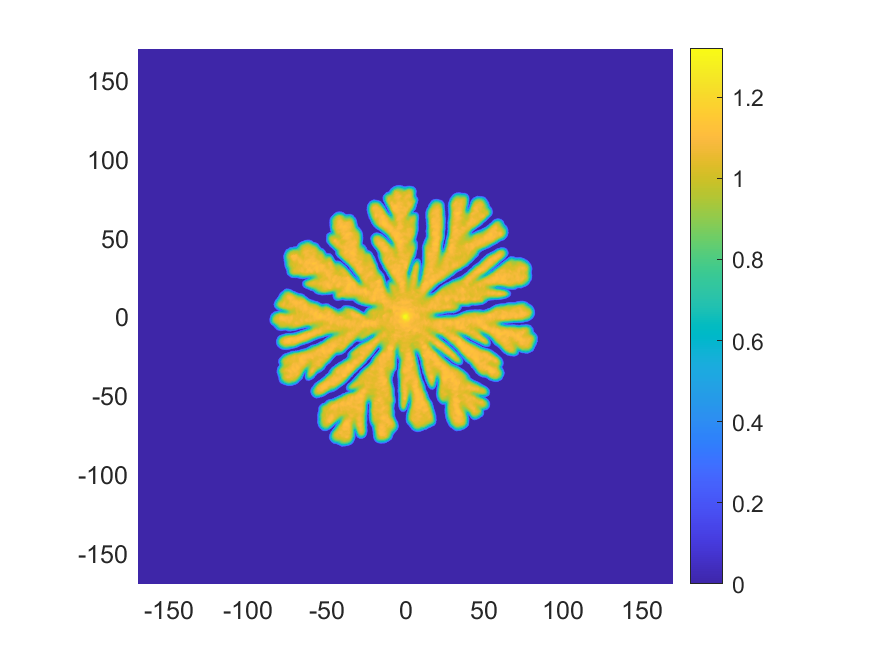}}%
	\subfigure[$t=1000$]{\includegraphics[width=0.245\textwidth,trim=18 8 18 8,clip]{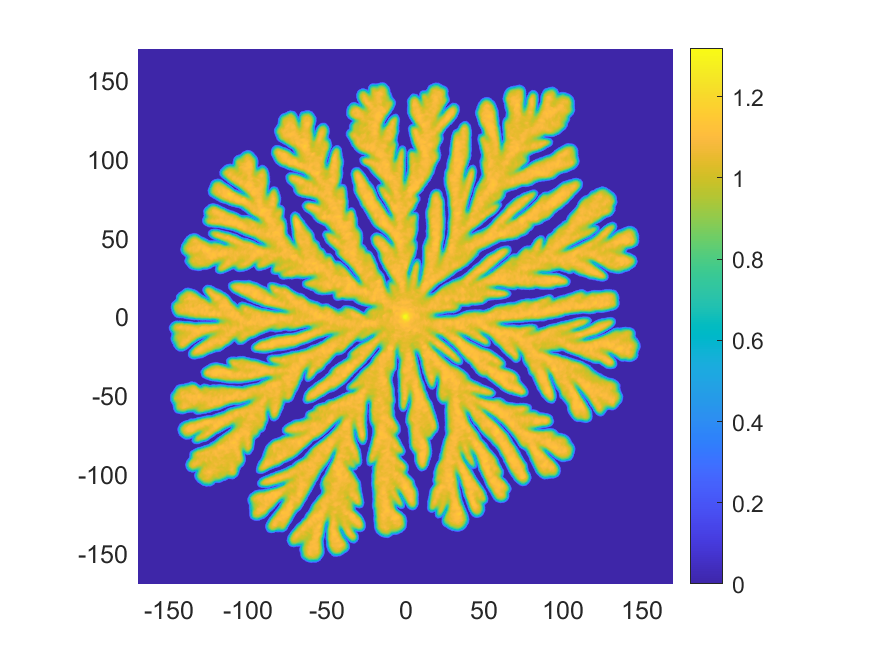}}
	\caption{The evolution of $n_h^m$ (top) and $b_h^m$ (bottom), $n_0=0.71$,  $\sigma_0=1$,  $\chi_0=5$.}
	\label{n_071_sigma_1_chi_5}
\end{figure}

\begin{figure}[!htbp]
	\centering
	\subfigure[$t=50$]{\includegraphics[width=0.245\textwidth,trim=18 8 18 8,clip]{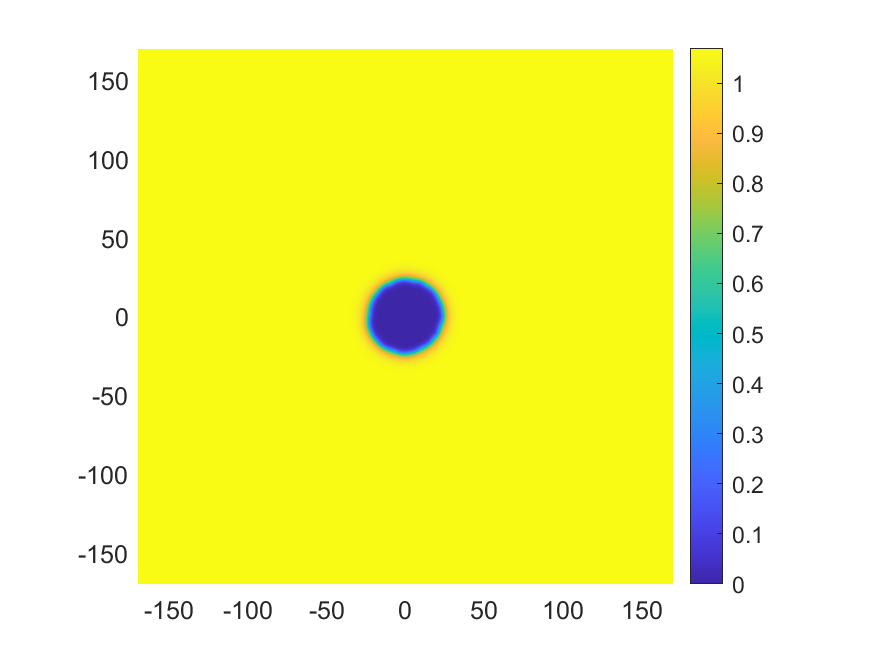}}%
	\subfigure[$t=100$]{\includegraphics[width=0.245\textwidth,trim=18 8 18 8,clip]{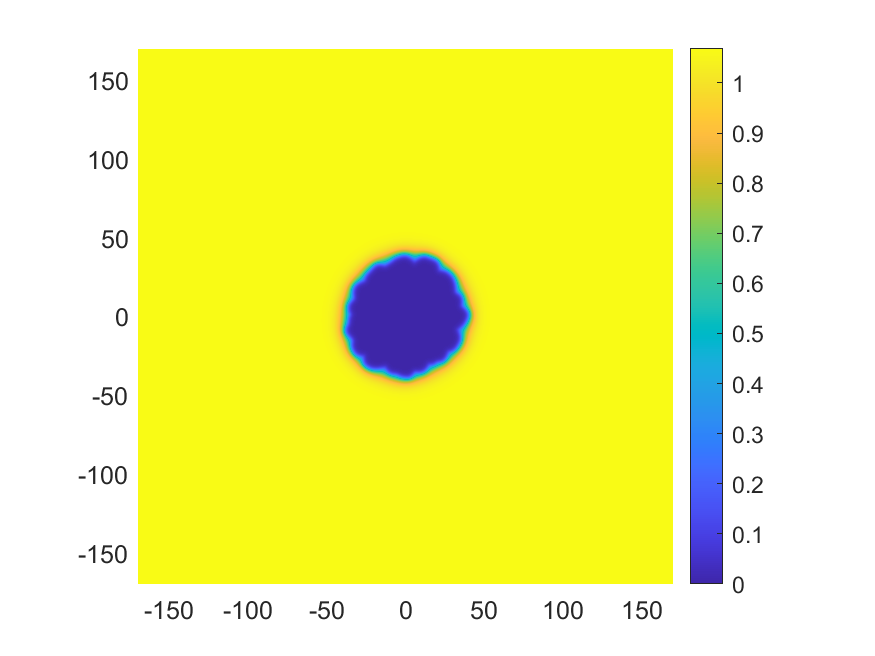}}%
	\subfigure[$t=250$]{\includegraphics[width=0.245\textwidth,trim=18 8 18 8,clip]{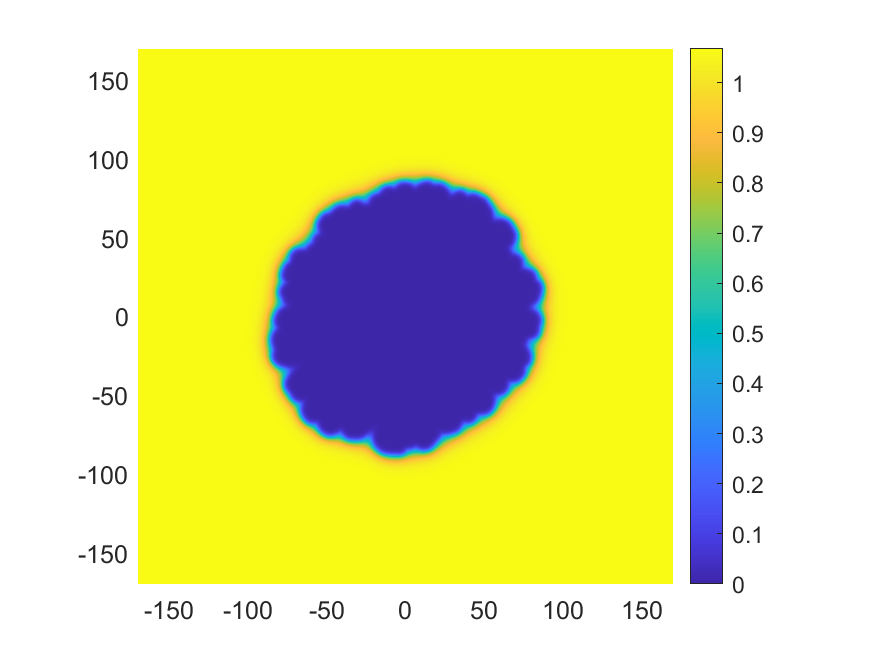}}%
	\subfigure[$t=400$]{\includegraphics[width=0.245\textwidth,trim=18 8 18 8,clip]{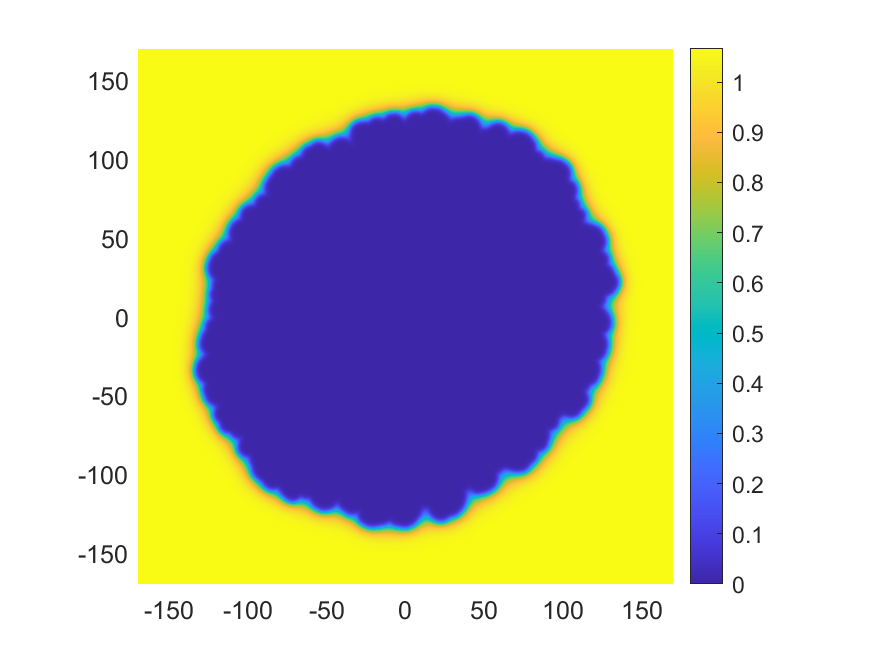}}\\[-1mm]
	\subfigure[$t=50$]{\includegraphics[width=0.245\textwidth,trim=18 8 18 8,clip]{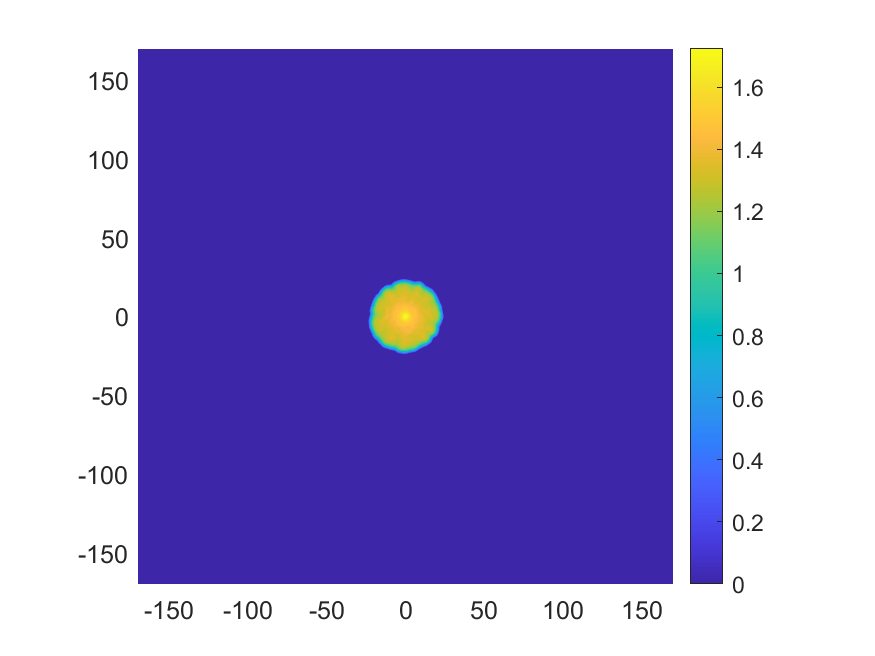}}%
	\subfigure[$t=100$]{\includegraphics[width=0.245\textwidth,trim=18 8 18 8,clip]{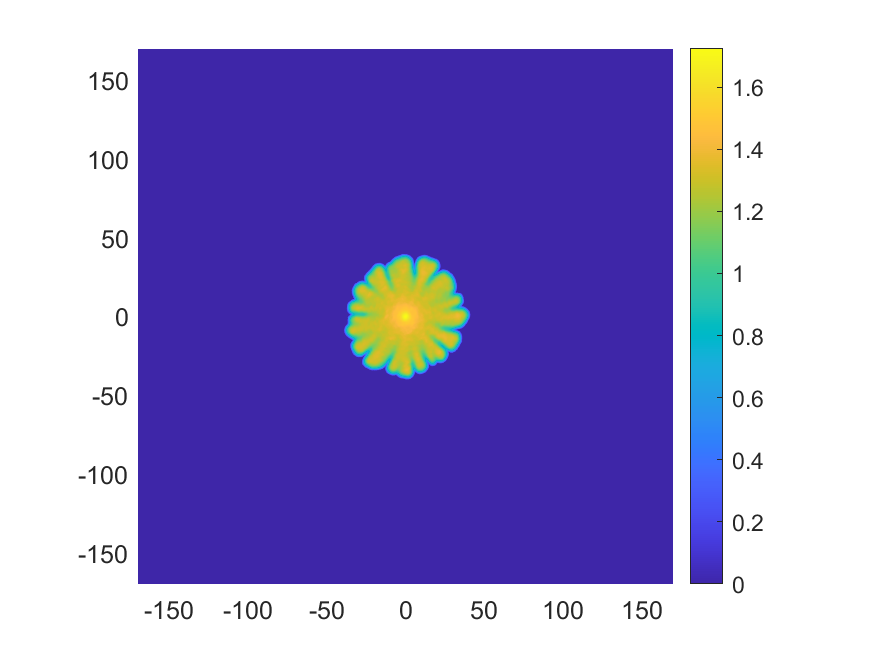}}%
	\subfigure[$t=250$]{\includegraphics[width=0.245\textwidth,trim=18 8 18 8,clip]{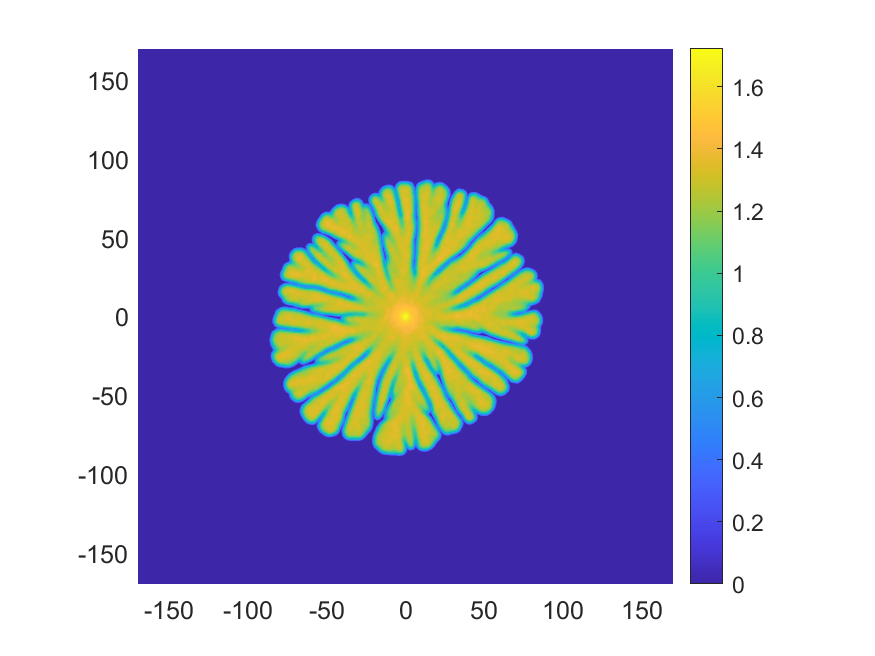}}%
	\subfigure[$t=400$]{\includegraphics[width=0.245\textwidth,trim=18 8 18 8,clip]{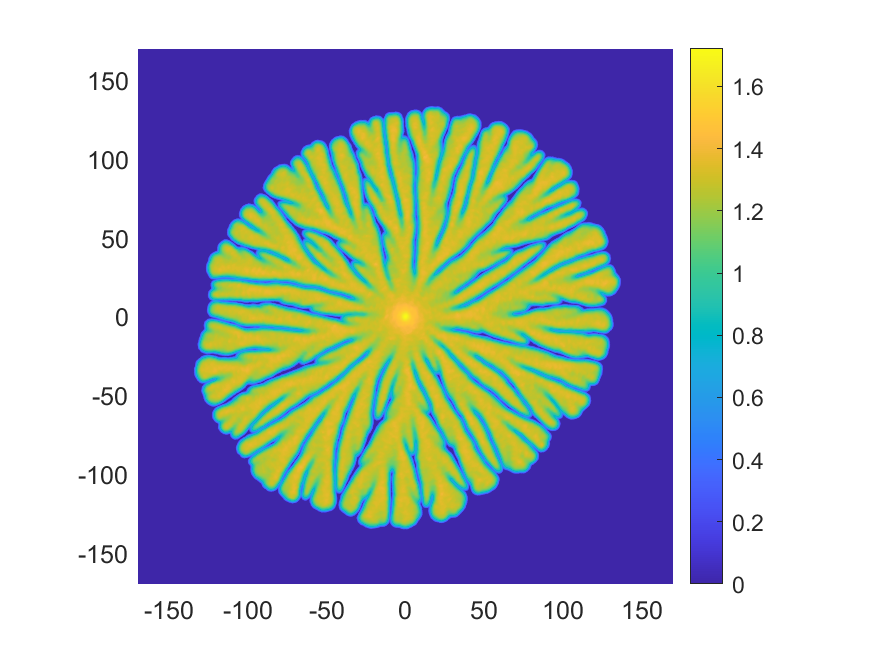}}
	\caption{The evolution of $n_h^m$ (top) and $b_h^m$ (bottom), $n_0=1.07$,  $\sigma_0=1$,  $\chi_0=0$.}
	\label{n_107_sigma_1_chi_0}
\end{figure}

\begin{figure}[!htbp]
	\centering
	\subfigure[$t=50$]{\includegraphics[width=0.245\textwidth,trim=18 8 18 8,clip]{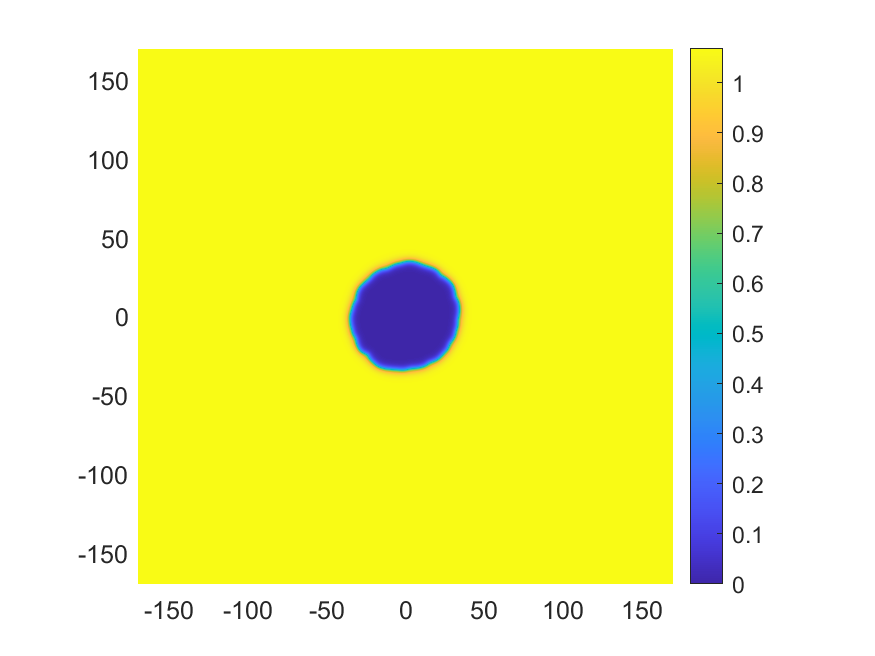}}%
	\subfigure[$t=100$]{\includegraphics[width=0.245\textwidth,trim=18 8 18 8,clip]{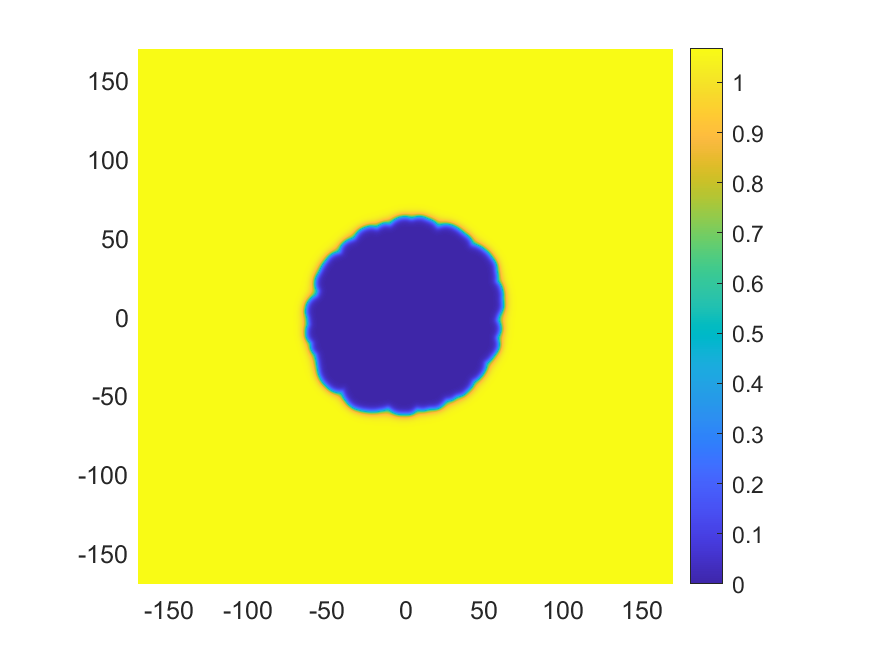}}%
	\subfigure[$t=200$]{\includegraphics[width=0.245\textwidth,trim=18 8 18 8,clip]{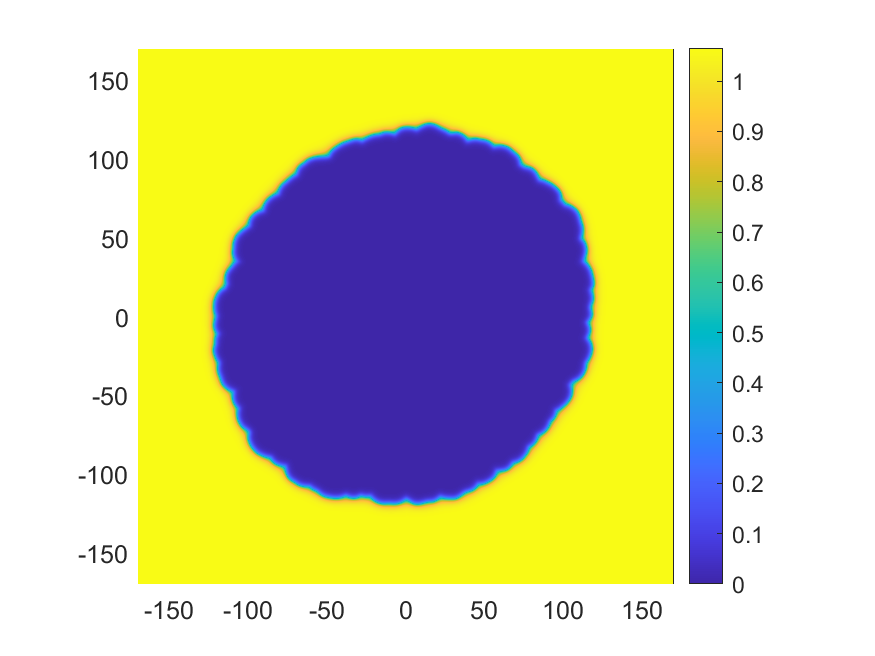}}%
	\subfigure[$t=250$]{\includegraphics[width=0.245\textwidth,trim=18 8 18 8,clip]{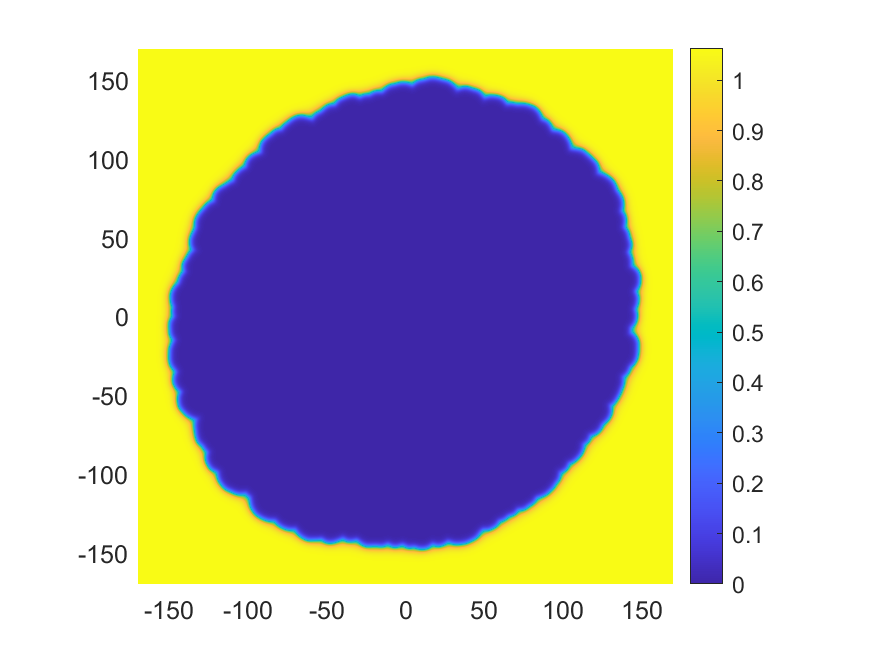}}\\[-1mm]
	\subfigure[$t=50$]{\includegraphics[width=0.245\textwidth,trim=18 8 18 8,clip]{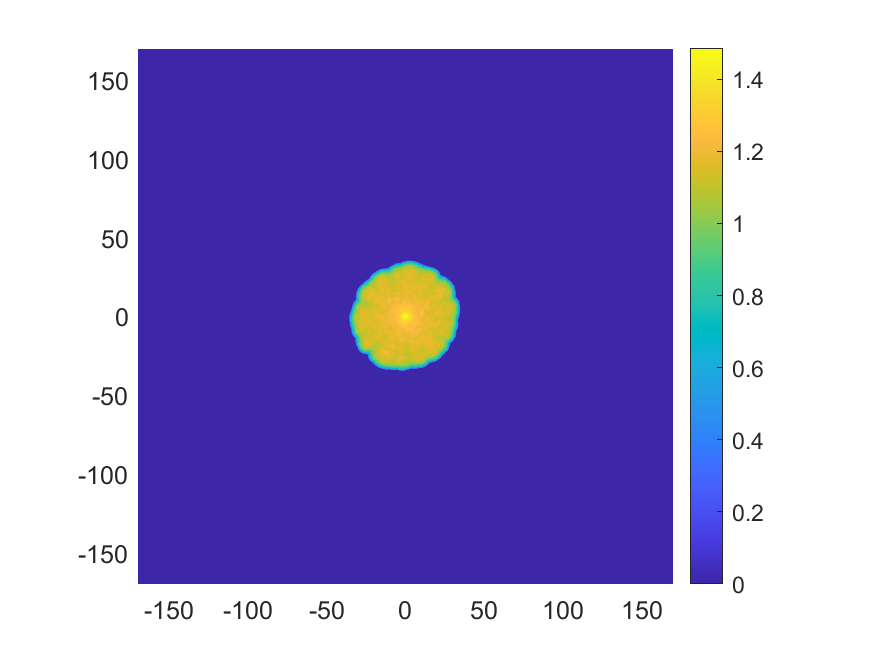}}%
	\subfigure[$t=100$]{\includegraphics[width=0.245\textwidth,trim=18 8 18 8,clip]{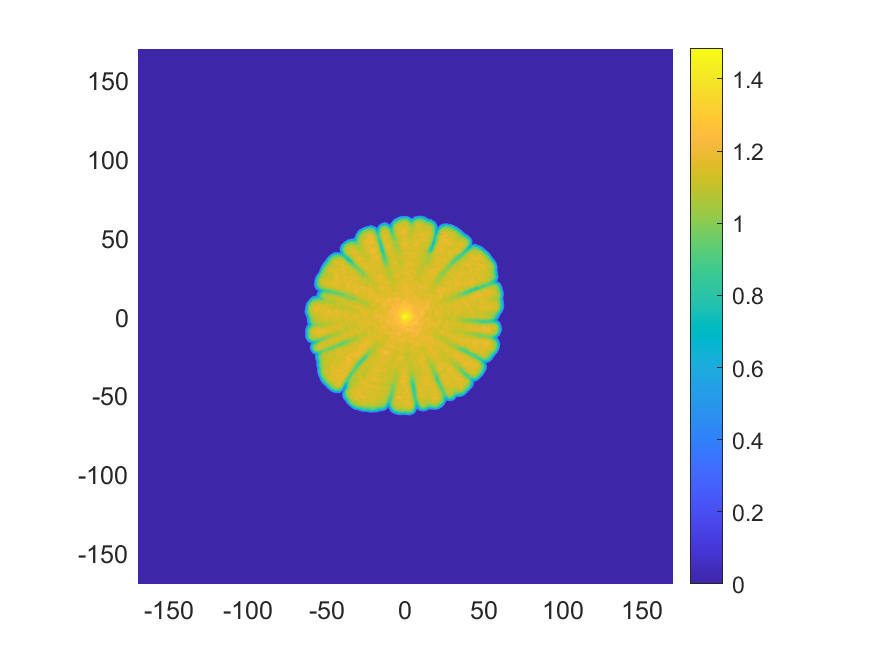}}%
	\subfigure[$t=200$]{\includegraphics[width=0.245\textwidth,trim=18 8 18 8,clip]{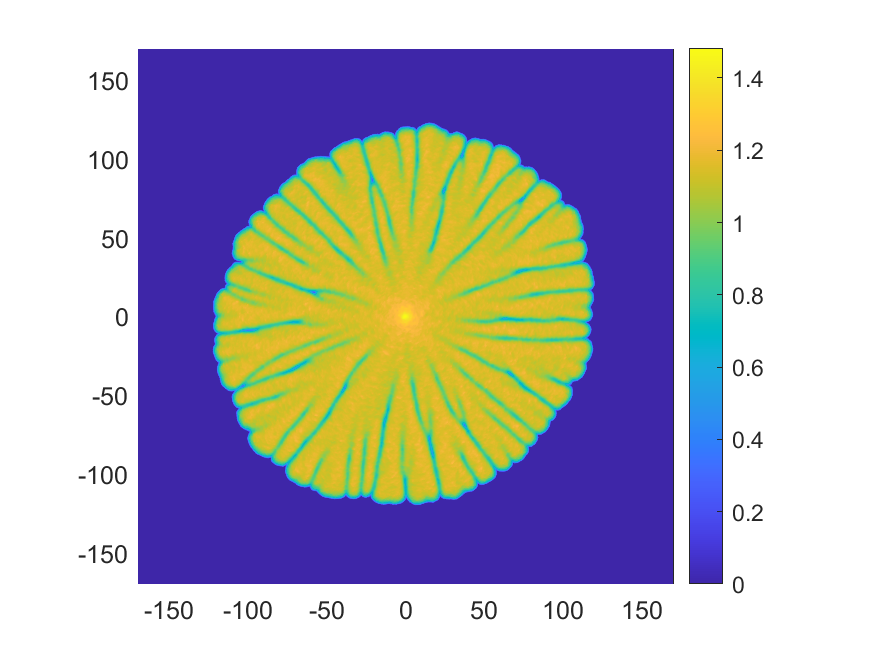}}%
	\subfigure[$t=250$]{\includegraphics[width=0.245\textwidth,trim=18 8 18 8,clip]{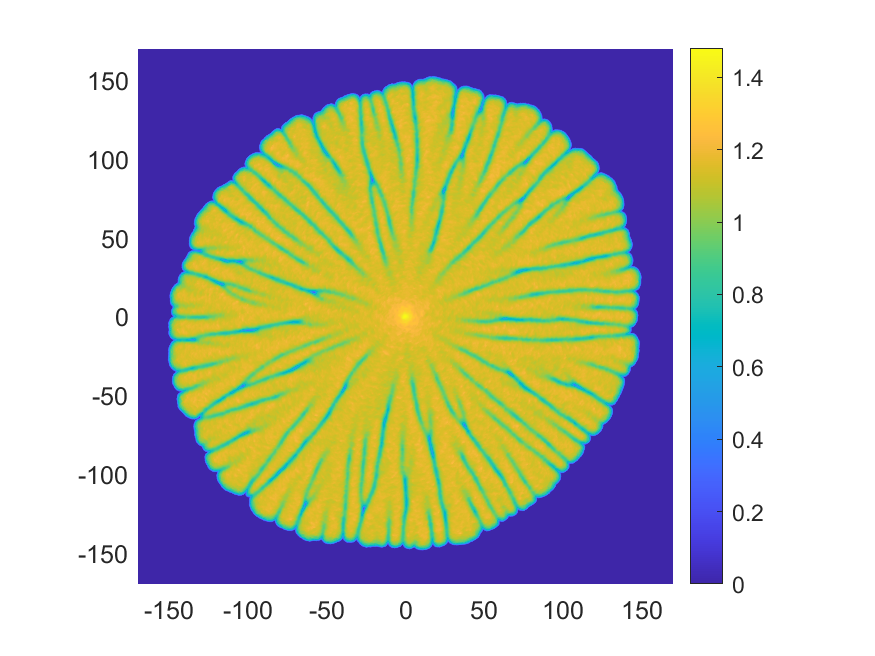}}
	\caption{The evolution of $n_h^m$ (top) and $b_h^m$ (bottom), $n_0=1.07$,  $\sigma_0=1$,  $\chi_0=5$.}
	\label{n_107_sigma_1_chi_5}
\end{figure}

We begin with a semi-solid agar medium characterized by $\sigma_0=1$. In this regime, the moderate agar concentration limits bacterial diffusion. Under this condition, the computed colony morphology lies in region E of the morphological diagram and displays a DBM-like branching structure (see \cite{KAMT1997}). Figures \ref{n_071_sigma_1_chi_0} and \ref{n_071_sigma_1_chi_5} show the numerical solutions for $n_0=0.71$ without chemotaxis ($\chi_0=0$) and with chemotactic sensitivity $\chi_0=5$, respectively.  In both nutrient settings, the overall morphology remains qualitatively similar in the presence and absence of chemotaxis. A closer inspection indicates that the number of branches is essentially unchanged. Nevertheless, once chemotactic signaling is introduced, the outward propagation of the branch envelope becomes noticeably faster. For the higher nutrient level $n_0=1.07$, the corresponding numerical results are shown in Figures \ref{n_107_sigma_1_chi_0} and \ref{n_107_sigma_1_chi_5}. In the non-chemotactic case $\chi_0=0$, the colony develops a dense radially branching pattern enclosed by an approximately circular envelope, while the bacterial density attains its maximum near the colony center. In fact, this peak mainly corresponds to the remnant of the initial inoculum. At the early stage, bacteria consume the nutrient around the inoculation site, and cells in the interior gradually lose their motility. By contrast, cells located near the colony front continue to migrate toward regions with higher nutrient availability, thereby driving branch formation. This indicates that the bacteria at the advancing front are the main source of interface instability, which is one of the essential features of the nonlinear cross-diffusion model.


\begin{figure}[!htbp]
	\centering
	\subfigure[$t=50$]{\includegraphics[width=0.245\textwidth,trim=18 8 18 8,clip]{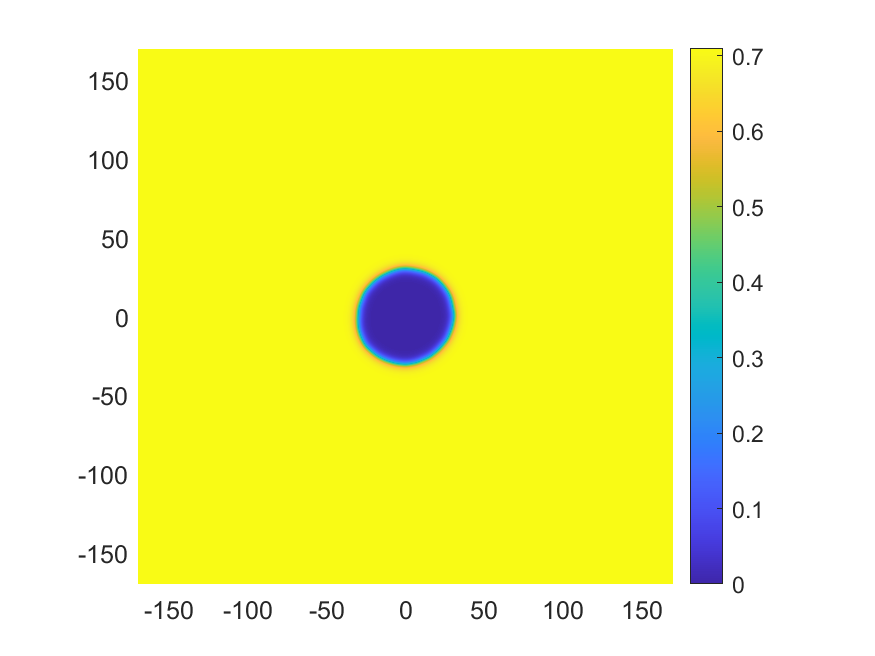}}%
	\subfigure[$t=100$]{\includegraphics[width=0.245\textwidth,trim=18 8 18 8,clip]{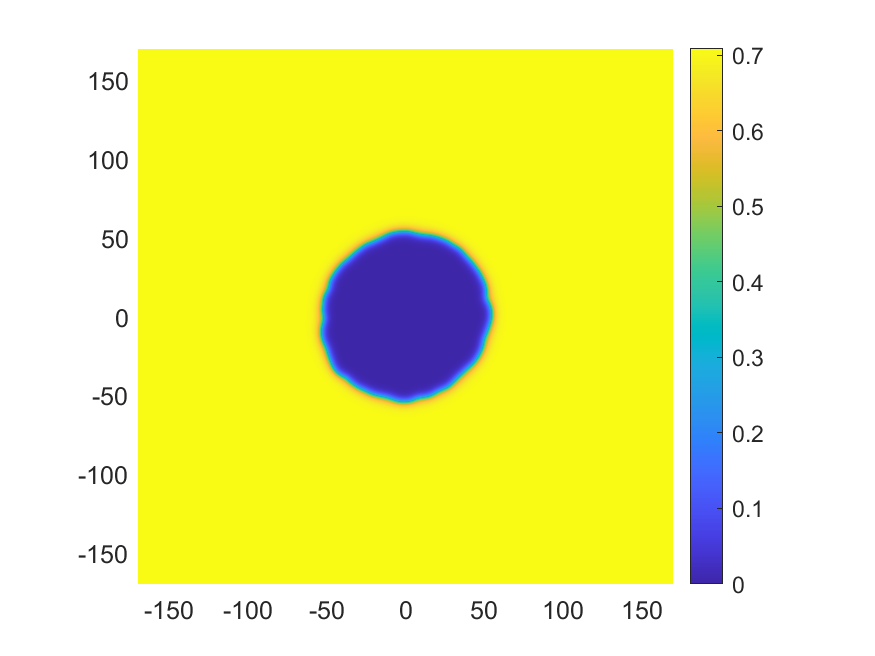}}%
	\subfigure[$t=150$]{\includegraphics[width=0.245\textwidth,trim=18 8 18 8,clip]{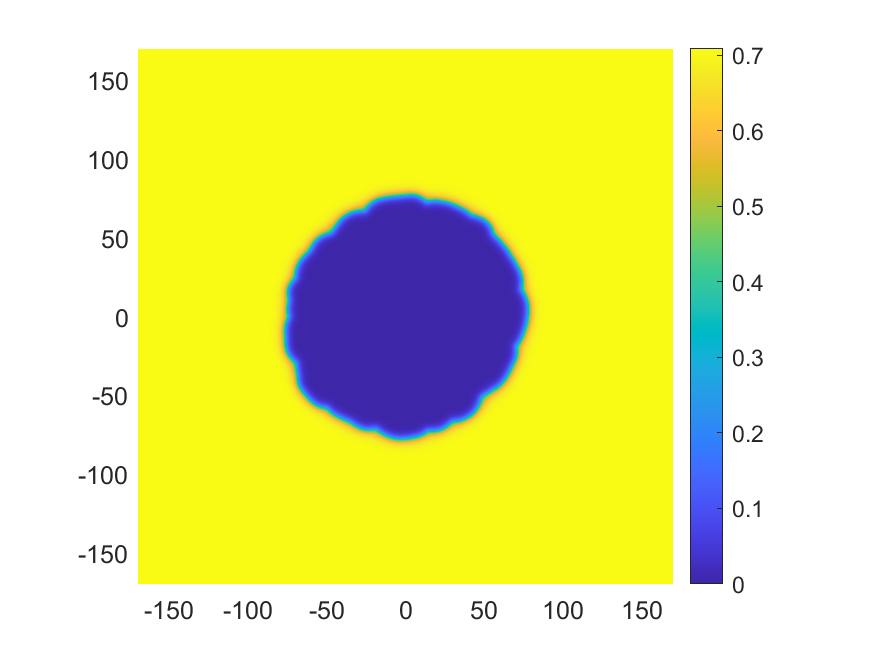}}%
	\subfigure[$t=300$]{\includegraphics[width=0.245\textwidth,trim=18 8 18 8,clip]{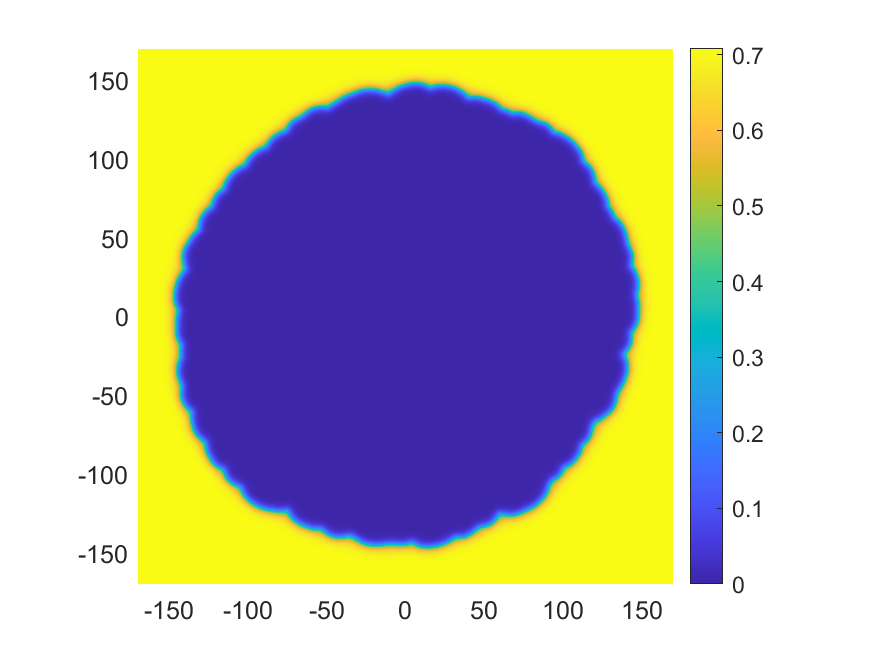}}\\[-1mm]
	\subfigure[$t=50$]{\includegraphics[width=0.245\textwidth,trim=18 8 18 8,clip]{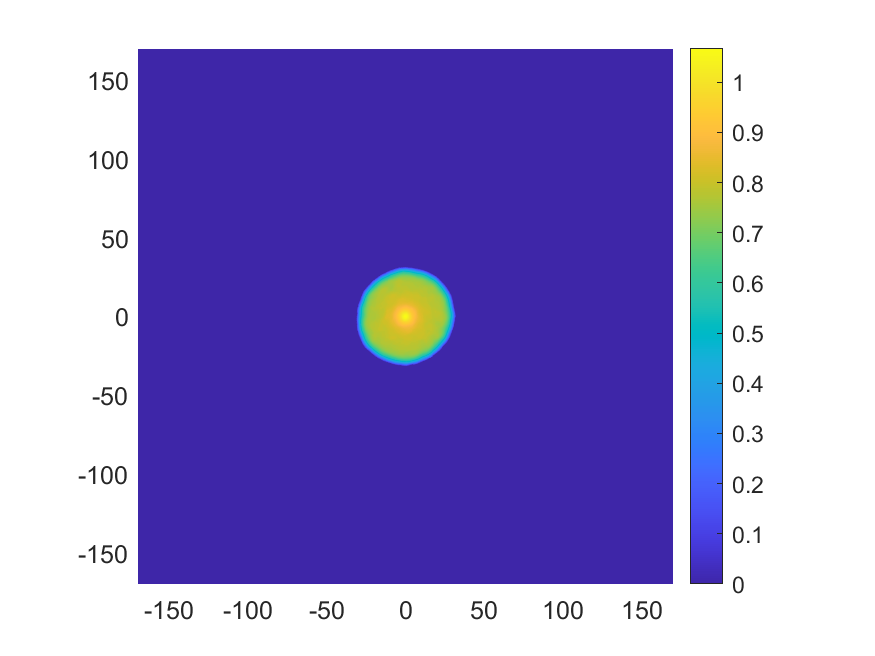}}%
	\subfigure[$t=100$]{\includegraphics[width=0.245\textwidth,trim=18 8 18 8,clip]{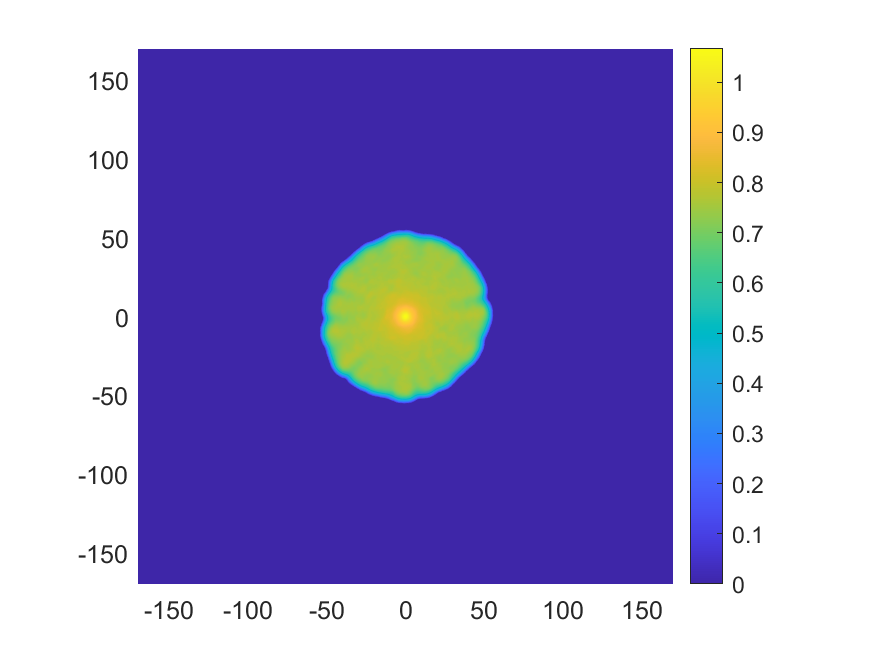}}%
	\subfigure[$t=150$]{\includegraphics[width=0.245\textwidth,trim=18 8 18 8,clip]{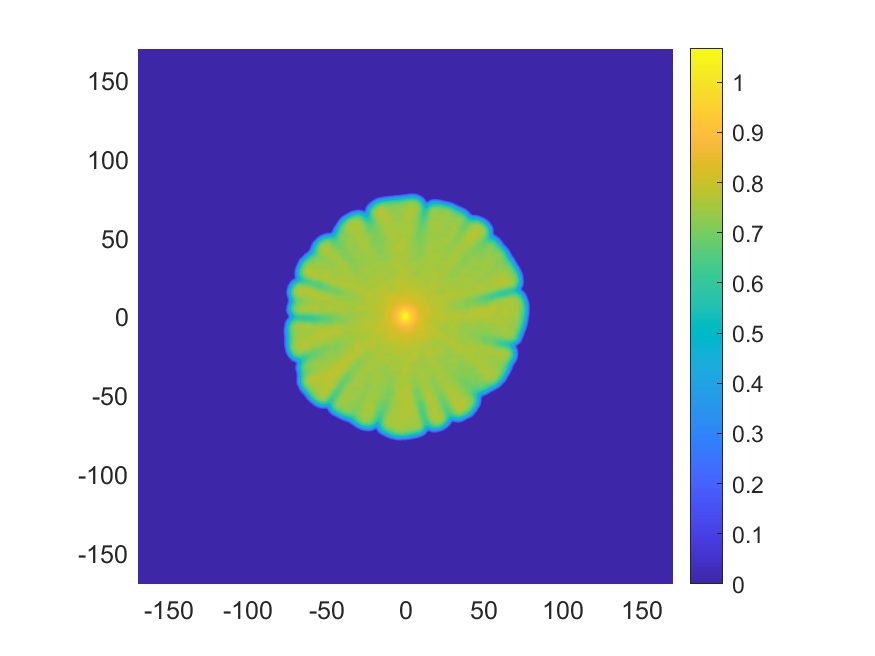}}%
	\subfigure[$t=300$]{\includegraphics[width=0.245\textwidth,trim=18 8 18 8,clip]{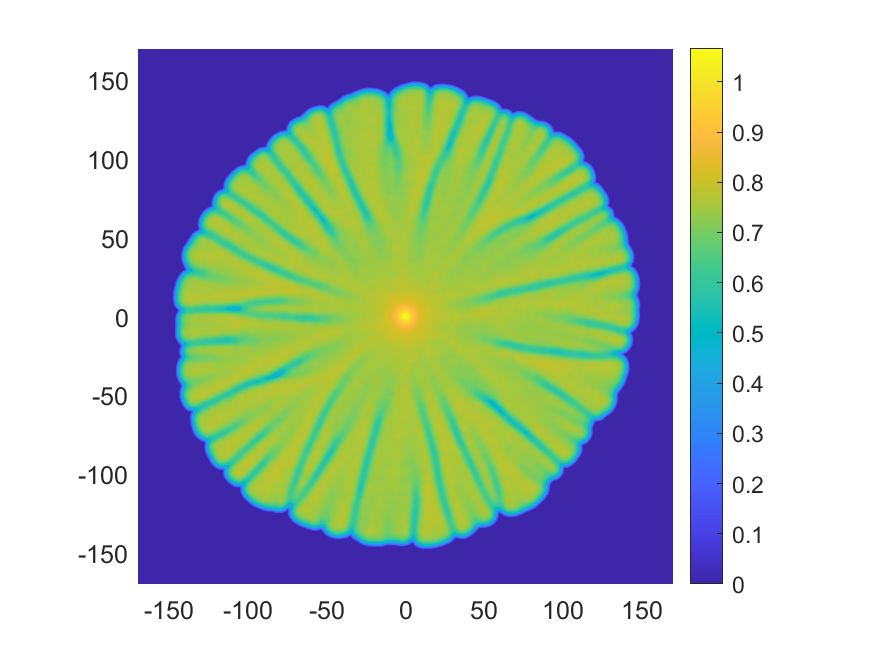}}
	\caption{The evolution of $n_h^m$ (top) and $b_h^m$ (bottom), $n_0=0.71$,  $\sigma_0=4$,  $\chi_0=0$.}
	\label{n_071_sigma_4_chi_0}
\end{figure}

\begin{figure}[!htbp]
	\centering
	\subfigure[$t=30$]{\includegraphics[width=0.245\textwidth,trim=18 8 18 8,clip]{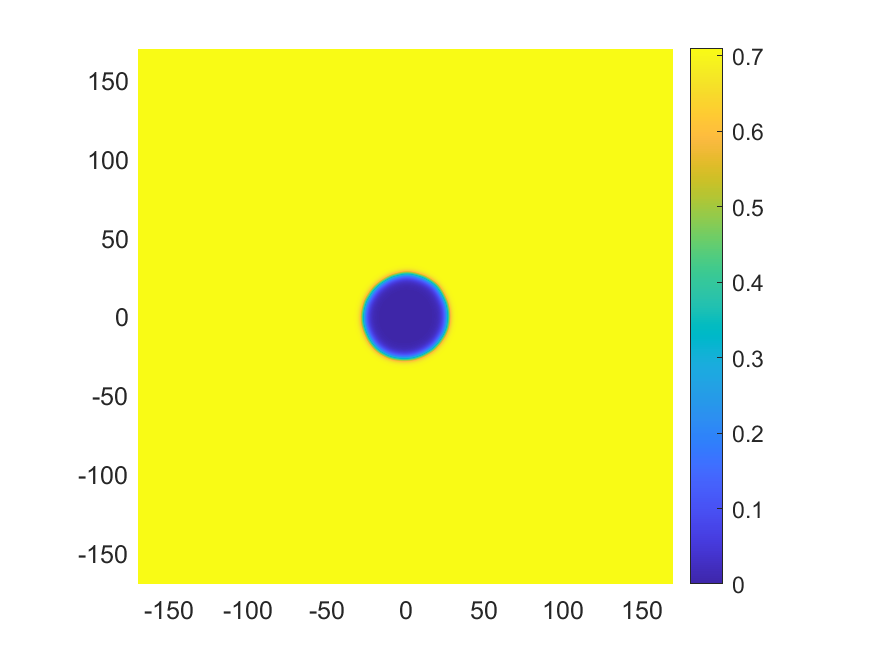}}%
	\subfigure[$t=50$]{\includegraphics[width=0.245\textwidth,trim=18 8 18 8,clip]{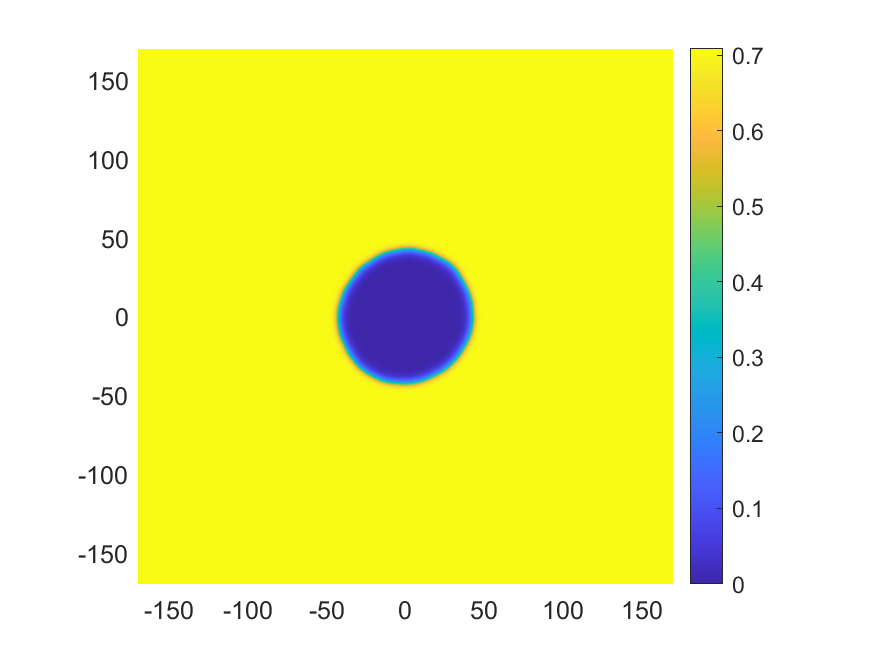}}%
	\subfigure[$t=100$]{\includegraphics[width=0.245\textwidth,trim=18 8 18 8,clip]{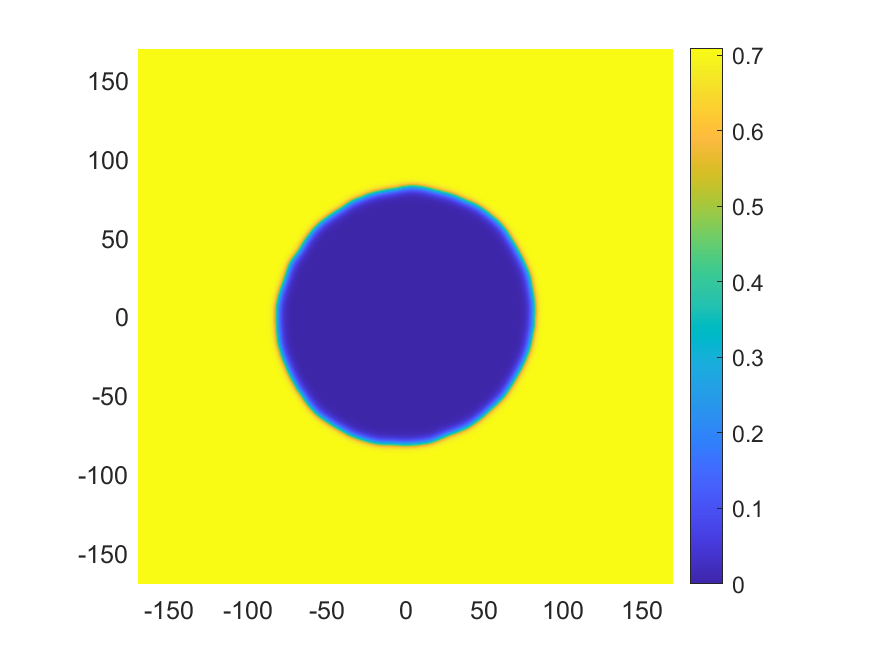}}%
	\subfigure[$t=150$]{\includegraphics[width=0.245\textwidth,trim=18 8 18 8,clip]{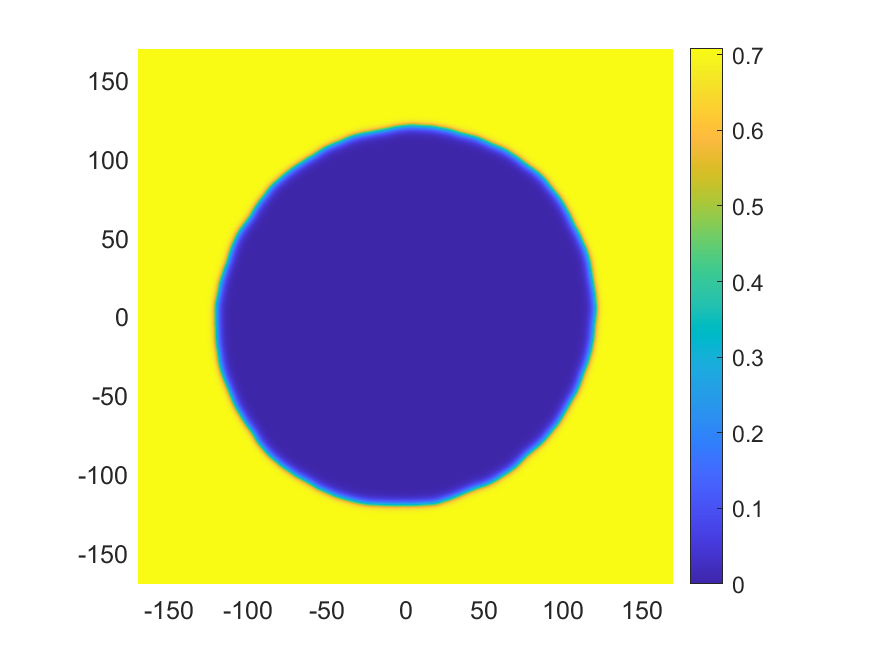}}\\[-1mm]
	\subfigure[$t=30$]{\includegraphics[width=0.245\textwidth,trim=18 8 18 8,clip]{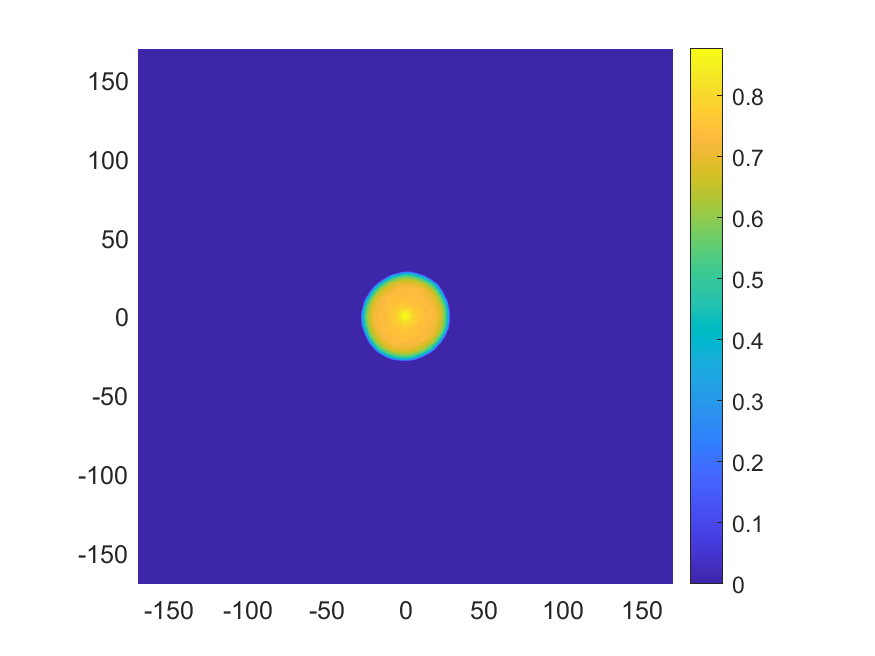}}%
	\subfigure[$t=50$]{\includegraphics[width=0.245\textwidth,trim=18 8 18 8,clip]{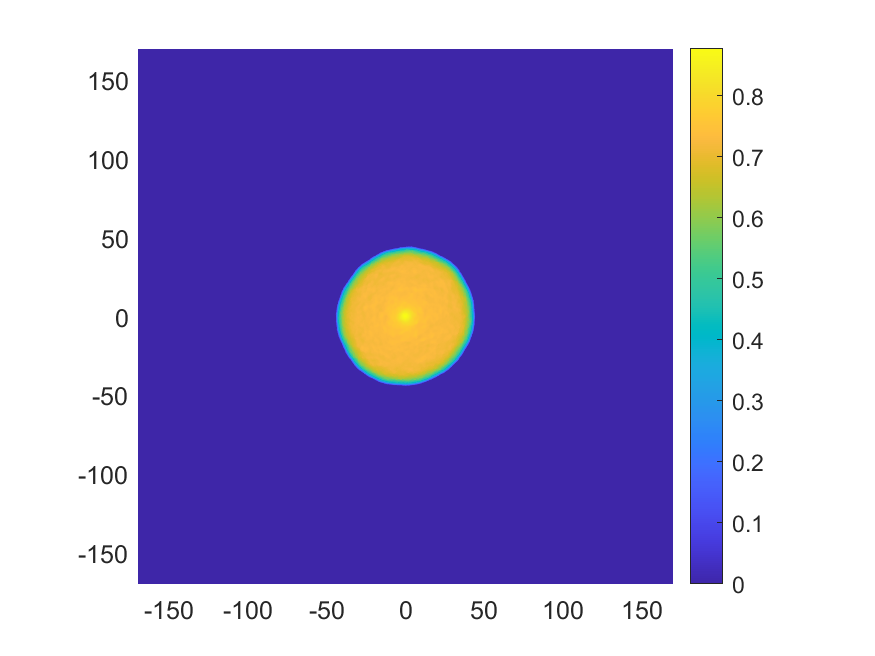}}%
	\subfigure[$t=100$]{\includegraphics[width=0.245\textwidth,trim=18 8 18 8,clip]{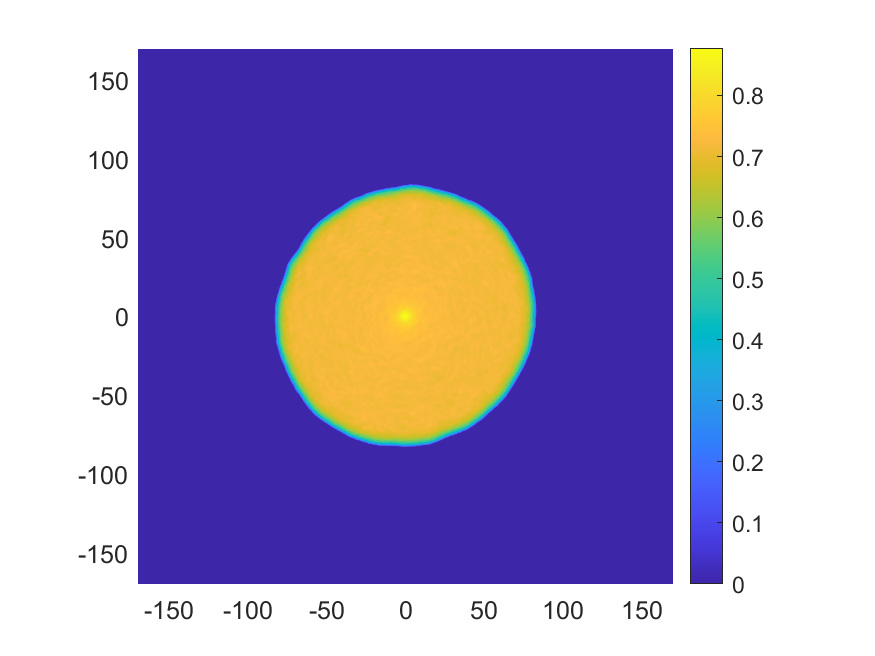}}%
	\subfigure[$t=150$]{\includegraphics[width=0.245\textwidth,trim=18 8 18 8,clip]{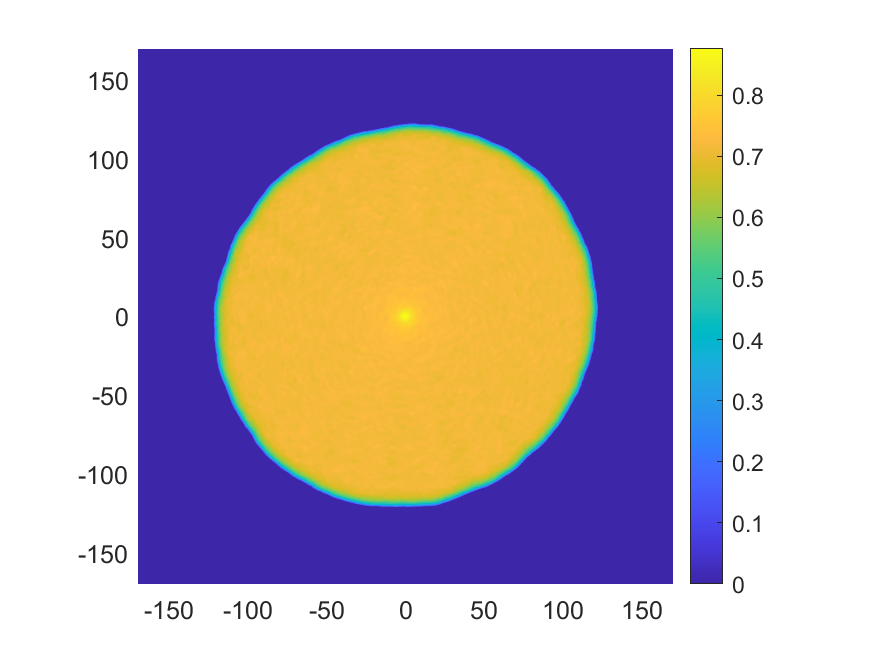}}
	\caption{The evolution of $n_h^m$ (top) and $b_h^m$ (bottom), $n_0=0.71$,  $\sigma_0=4$,  $\chi_0=5$.}
	\label{n_071_sigma_4_chi_5}
\end{figure}

\begin{figure}[!htbp]
	\centering
	\subfigure[$t=30$]{\includegraphics[width=0.245\textwidth,trim=18 8 18 8,clip]{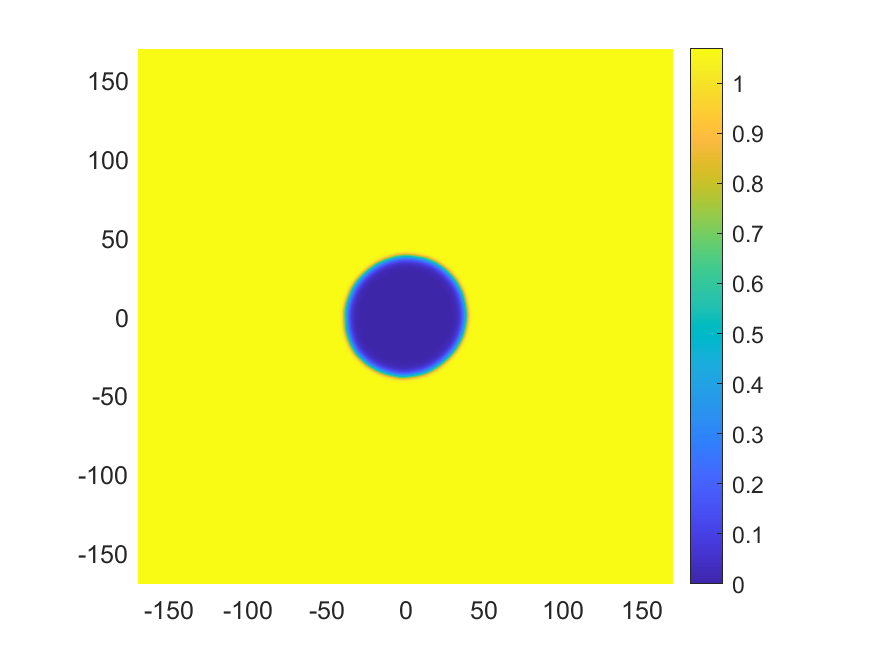}}%
	\subfigure[$t=50$]{\includegraphics[width=0.245\textwidth,trim=18 8 18 8,clip]{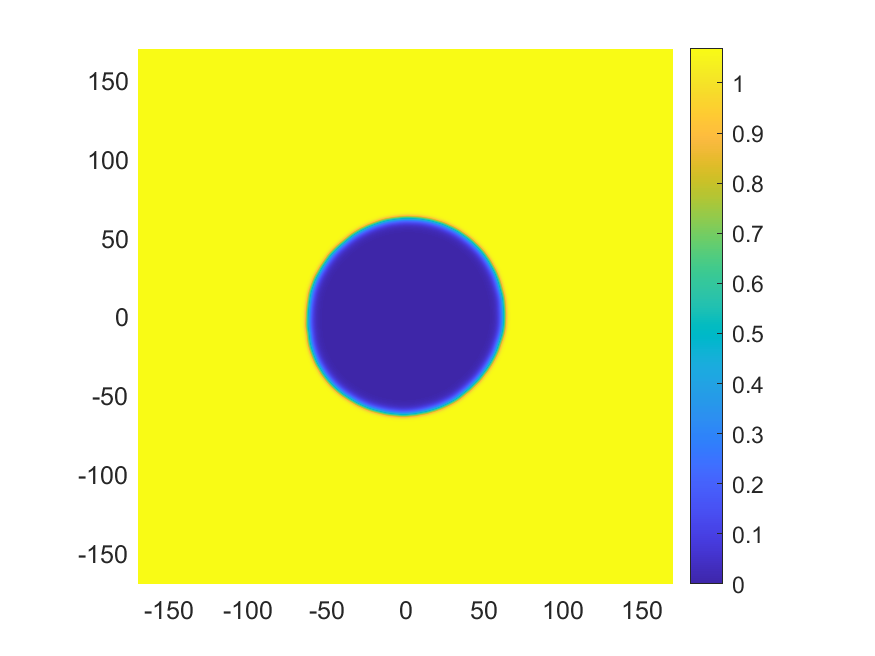}}%
	\subfigure[$t=100$]{\includegraphics[width=0.245\textwidth,trim=18 8 18 8,clip]{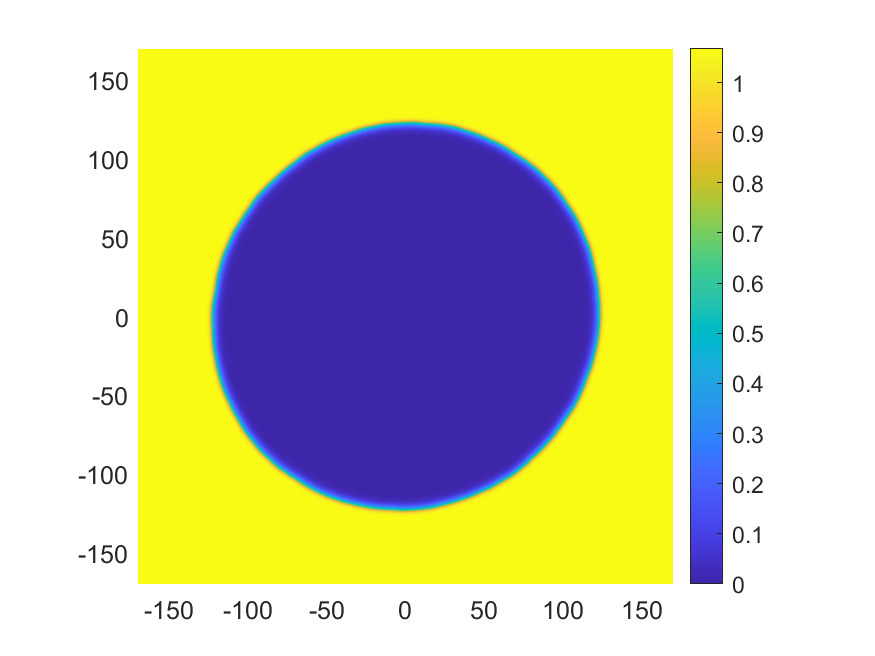}}%
	\subfigure[$t=130$]{\includegraphics[width=0.245\textwidth,trim=18 8 18 8,clip]{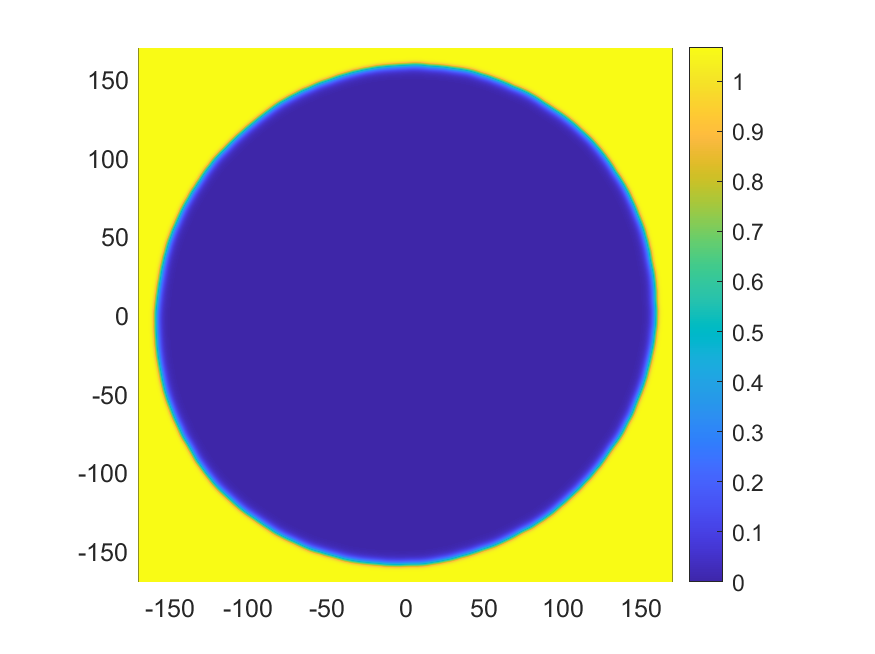}}\\[-1mm]
	\subfigure[$t=30$]{\includegraphics[width=0.245\textwidth,trim=18 8 18 8,clip]{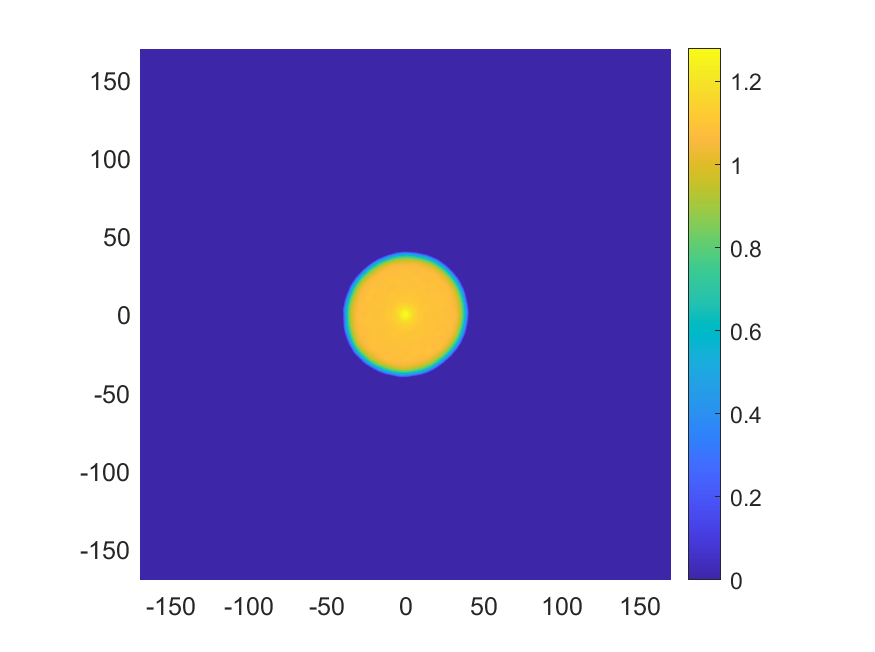}}%
	\subfigure[$t=50$]{\includegraphics[width=0.245\textwidth,trim=18 8 18 8,clip]{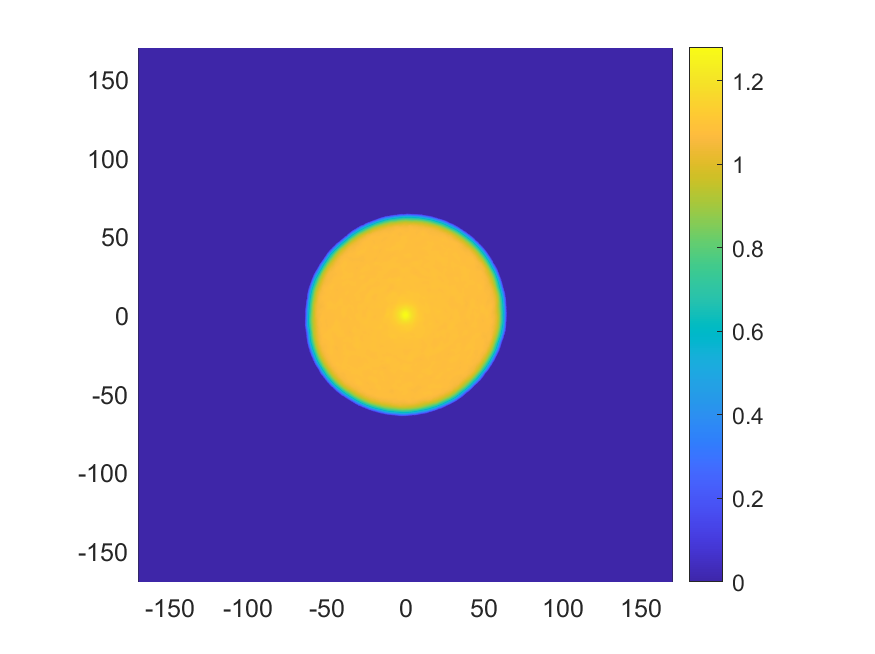}}%
	\subfigure[$t=100$]{\includegraphics[width=0.245\textwidth,trim=18 8 18 8,clip]{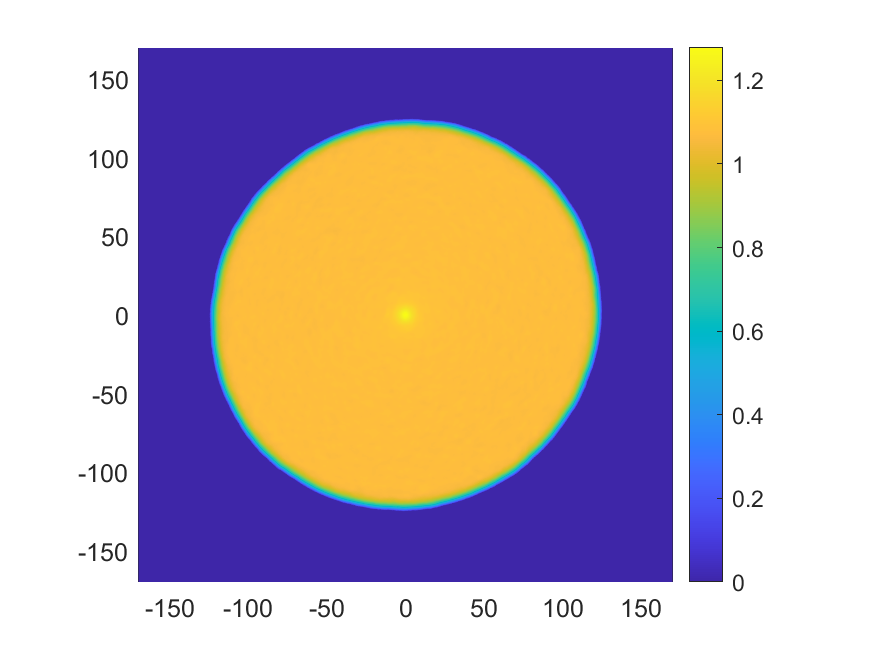}}%
	\subfigure[$t=130$]{\includegraphics[width=0.245\textwidth,trim=18 8 18 8,clip]{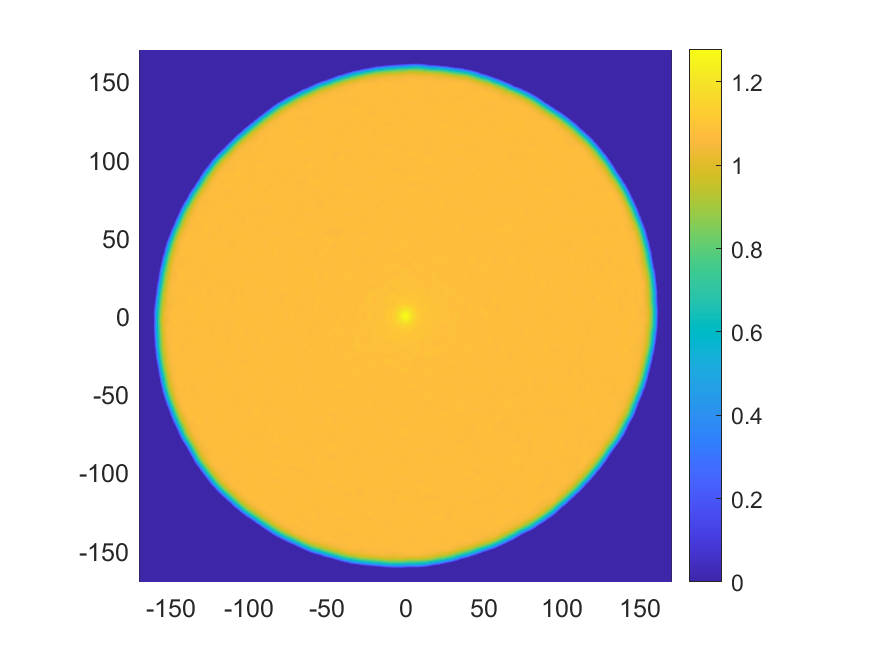}}
	\caption{The evolution of $n_h^m$ (top) and $b_h^m$ (bottom), $n_0=1.07$,  $\sigma_0=4$,  $\chi_0=0$.}
	\label{n_107_sigma_4_chi_0}
\end{figure}

\begin{figure}[!htbp]
	\centering
	\subfigure[$t=10$]{\includegraphics[width=0.245\textwidth,trim=18 8 18 8,clip]{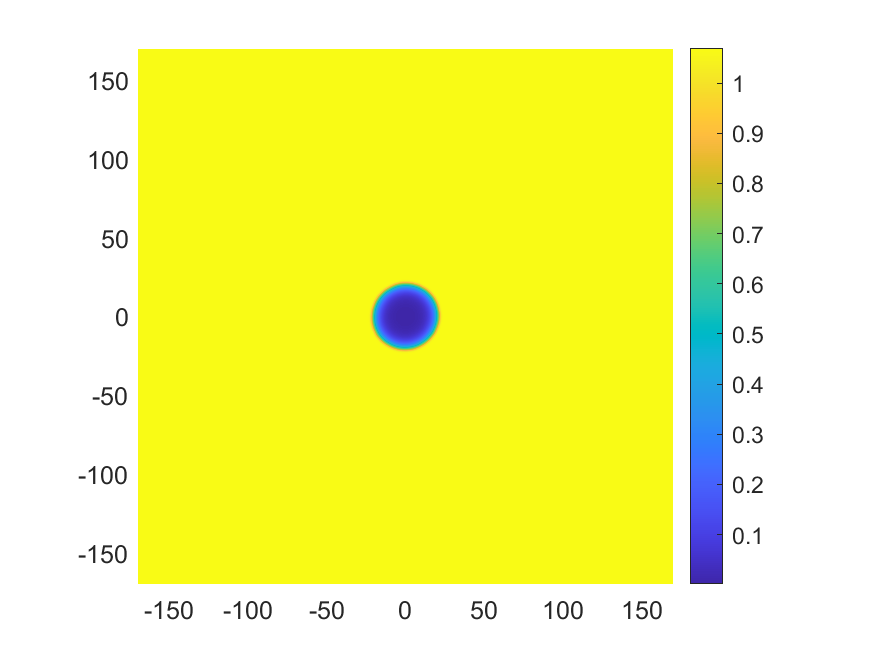}}%
	\subfigure[$t=30$]{\includegraphics[width=0.245\textwidth,trim=18 8 18 8,clip]{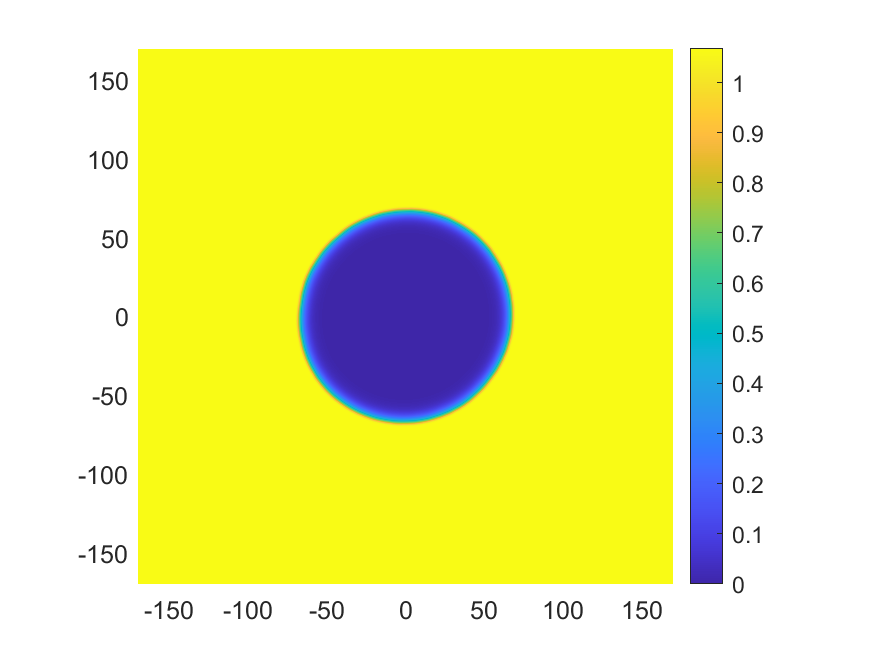}}%
	\subfigure[$t=50$]{\includegraphics[width=0.245\textwidth,trim=18 8 18 8,clip]{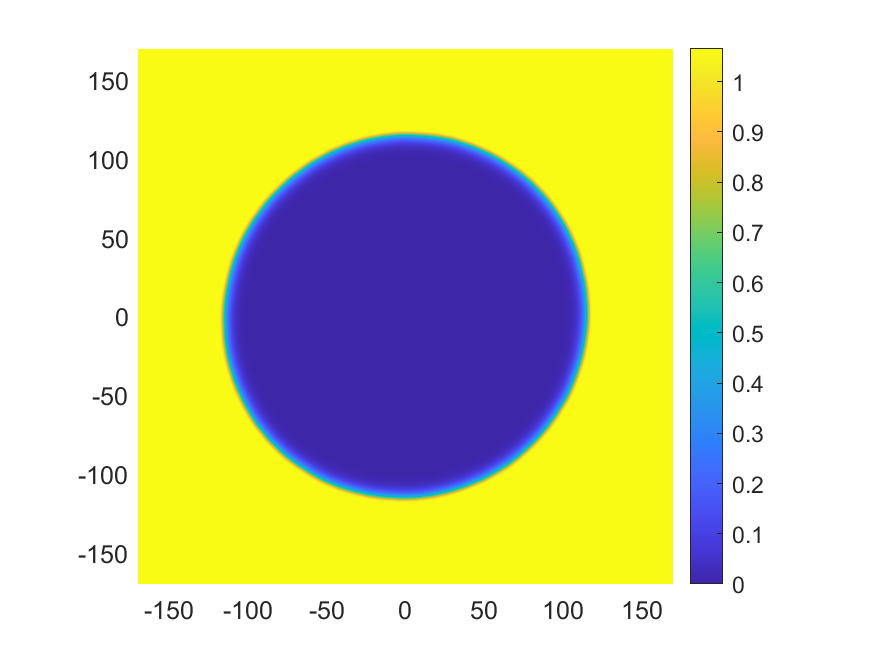}}%
	\subfigure[$t=60$]{\includegraphics[width=0.245\textwidth,trim=18 8 18 8,clip]{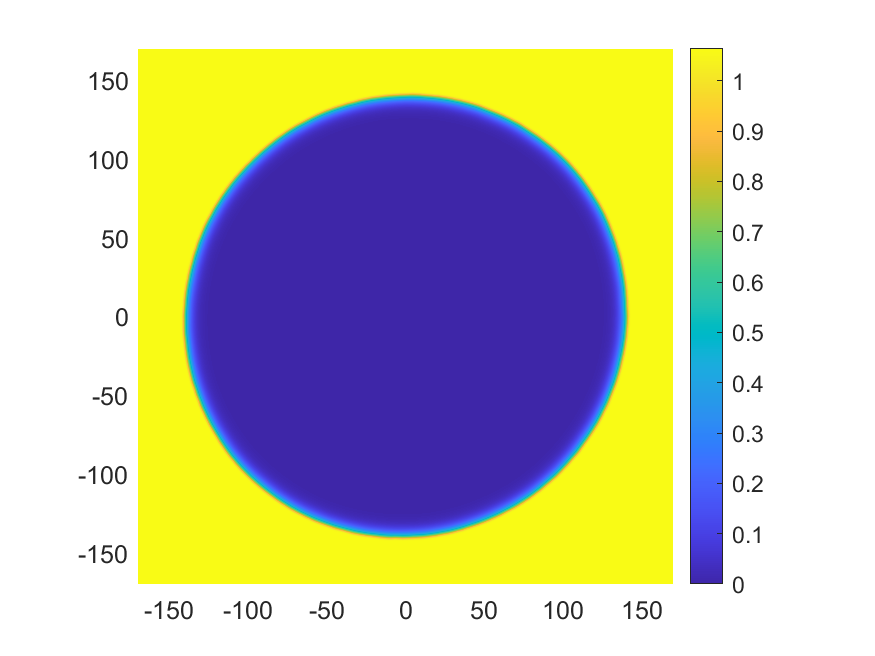}}\\[-1mm]
	\subfigure[$t=10$]{\includegraphics[width=0.245\textwidth,trim=18 8 18 8,clip]{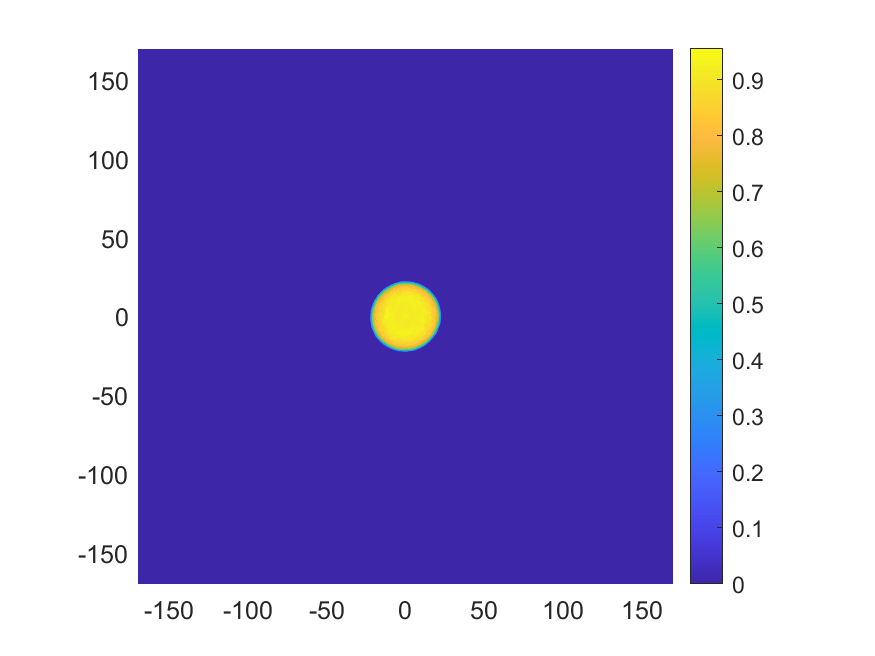}}%
	\subfigure[$t=30$]{\includegraphics[width=0.245\textwidth,trim=18 8 18 8,clip]{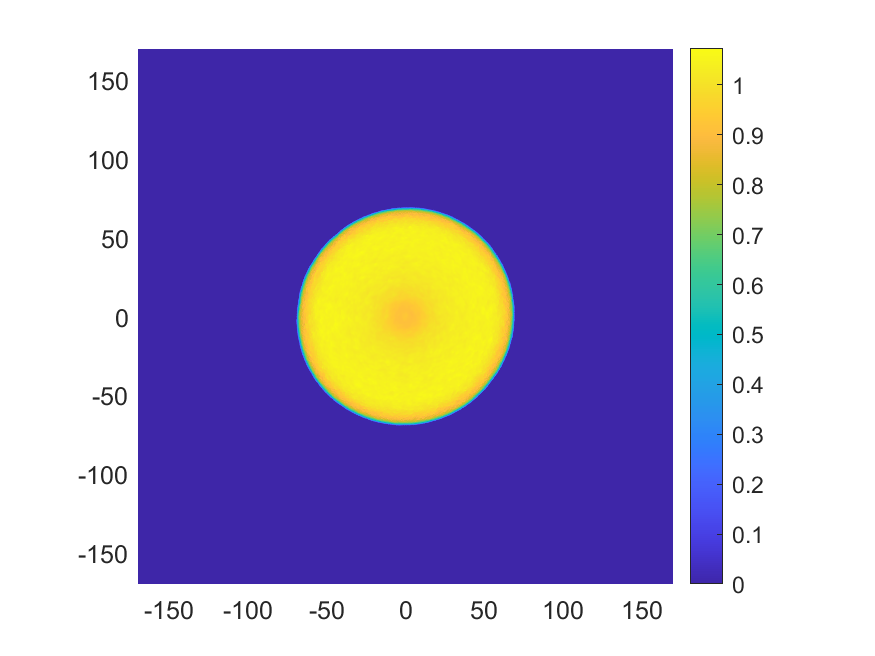}}%
	\subfigure[$t=50$]{\includegraphics[width=0.245\textwidth,trim=18 8 18 8,clip]{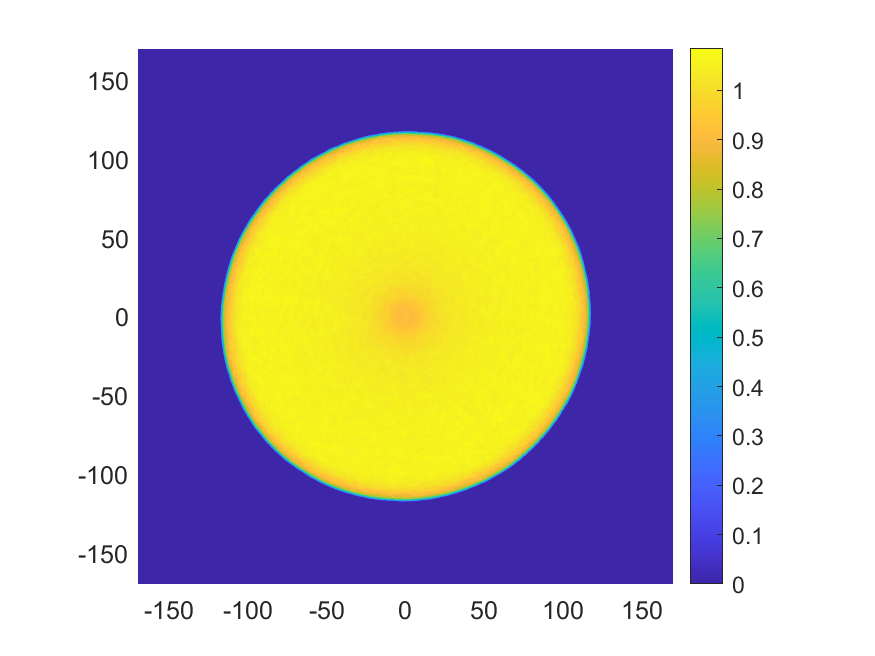}}%
	\subfigure[$t=60$]{\includegraphics[width=0.245\textwidth,trim=18 8 18 8,clip]{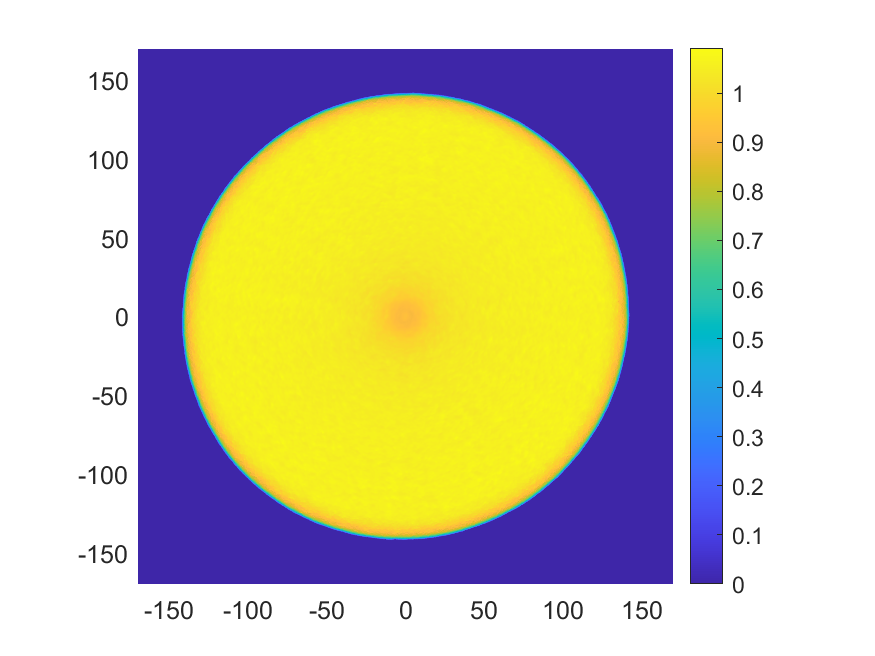}}
	\caption{The evolution of $n_h^m$ (top) and $b_h^m$ (bottom), $n_0=1.07$,  $\sigma_0=4$,  $\chi_0=7.5$.}
	\label{n_107_sigma_4_chi_75}
\end{figure}

Then we take $\sigma_0=4$, corresponding to a lower agar concentration and thus enhanced bacterial diffusivity. In this regime, sufficiently high nutrient levels lead to colony morphologies within region D of the morphological diagram, while lower values of $n_0$ may still give rise to region-E-type patterns (see \cite{KAMT1997}). To examine the transition from region E to region D, we first fix the initial nutrient concentration at $n_0=0.71$ and investigate the effect of chemotaxis by varying $\chi_0$. Comparing Figures \ref{n_071_sigma_4_chi_0} and
\ref{n_071_sigma_4_chi_5}, we observe that increasing $\chi_0$ not only accelerates colony growth but also qualitatively changes the colony morphology. A plausible explanation is that, when $\chi_0$ is sufficiently large, the propagation speed of the leading front exceeds the characteristic diffusive spreading rate. In this case, the mechanism responsible for fjord formation is suppressed, so that fine branches begin to merge, the pattern gradually loses its branched structure, and the colony eventually approaches a more homogeneous disk-like morphology. For $n_0=1.07$, a more pronounced chemotactic effect is observed; see Figures \ref{n_107_sigma_4_chi_0} and \ref{n_107_sigma_4_chi_75}. The results indicate that, near the transition between regions E and D in the morphological diagram, chemotaxis stabilizes the colony interface and thereby suppresses the instabilities that lead to branch formation. As the chemotactic effect becomes stronger, the emergence of new branches is reduced and the colony develops a more uniform morphology, while its overall propagation is accelerated. The numerical simulations successfully reproduce these characteristic behaviors, providing further evidence of the reliability of the present scheme for capturing the morphological evolution of bacterial colonies.


\section{Conclusions}
\label{5_conclusion}
In this paper, we develop high-order structure-preserving fully discrete schemes for a chemotaxis model describing the growth of bacterial colonies. We first establish two important structural properties of its continuous system, namely, total mass conservation and the nonnegativity of the solutions. Then, combining the mass-lumped finite element method with the KKT correction and the Lagrange multiplier method, we construct numerical schemes that preserve these properties at the discrete level. Moreover, we prove the existence and unique solvability of the nonlinear algebraic equation arising from the Lagrange multiplier correction and the local convergence of the corresponding secant iteration. Numerical experiments confirm the expected convergence orders and structure-preserving properties, and successfully capture the evolution of bacterial colony patterns. However, deriving rigorous error estimates remains challenging because of the strongly nonlinear cross-diffusive coupling in the bacterial-density equation. The error analysis strategy in \cite{Wangliufeng}, which separates the error analysis into temporal and spatial parts and derives suitable regularity estimates for the temporal semi-discrete solutions, provides a possible approach to handling such nonlinear terms. Moreover, linear systems with variable coefficients need to be solved at each time step to ensure the approximation accurately in the growth of the bacterial colony. How to design more efficient scheme is still an interesting topic. All of these will be studied in future.

\section*{Statements and Declarations}
\textbf{Competing Interests}
The authors have no relevant financial or non-financial interests to disclose.\\
\\
\textbf{Author Contributions}
Shuyan Ren: Method, Implementation, Numerical analysis, Writing,  Editing.  Kun Wang:  Method, Numerical analysis, Writing, Editing. Nianyu Yi:  Method, Numerical analysis, Writing, Editing.\\
\\
\textbf{Data Availability}
Enquiries about data availability should be directed to
the authors.
%\end{linenumbers}
\end{document}